\documentclass[12pt]{article}

\usepackage{amsfonts,amssymb,subeqnarray,eqsection,indent}
\usepackage{bm}    
\usepackage{url}

\begin{document}

\title{\vspace*{-2cm}
       Complete monotonicity for inverse powers \\
       of some combinatorially defined polynomials}

\date{January 9, 2013 \\
      revised November 10, 2013 \\[5mm]
      \parbox{5in}{{\bf Expanded version:}  Contains Appendices~A and B
            that are not included, due to space constraints, in the version
            ({\tt arXiv:1301.2449v2}) that will be published in
            {\em Acta Mathematica}\/.}
     }

\author{
     {\small Alexander D.~Scott} \\[-2mm]
     {\small\it Mathematical Institute}  \\[-2mm]
     {\small\it University of Oxford} \\[-2mm]
     {\small\it 24--29 St.~Giles'}\/ \\[-2mm]
     {\small\it Oxford OX1 3LB, UK}                         \\[-2mm]
     {\small\tt scott@maths.ox.ac.uk}                        \\[3mm]
     {\small Alan D.~Sokal\thanks{Also at Department of Mathematics,
           University College London, London WC1E 6BT, UK.}}  \\[-2mm]
     {\small\it Department of Physics}       \\[-2mm]
     {\small\it New York University}         \\[-2mm]
     {\small\it 4 Washington Place}          \\[-2mm]
     {\small\it New York, NY 10003 USA}      \\[-2mm]
     {\small\tt sokal@nyu.edu}               \\[-2mm]
     {\protect\makebox[5in]{\quad}}  
     \\[-4mm]
}

\maketitle
\thispagestyle{empty}   

\vspace*{-8mm}

\begin{abstract}
We prove the complete monotonicity on $(0,\infty)^n$
for suitable inverse powers of the spanning-tree polynomials of graphs
and, more generally, of the basis generating polynomials of
certain classes of matroids.
This generalizes a result of Szeg\H{o}
and answers, among other things,
a long-standing question of Lewy and Askey
concerning the positivity of Taylor coefficients for certain rational functions.
Our proofs are based on two {\em ab initio}\/ methods for proving that
$P^{-\beta}$ is completely monotone on a convex cone $C$:
the determinantal method and the quadratic-form method.
These methods are closely connected with
harmonic analysis on Euclidean Jordan algebras
(or equivalently on symmetric cones).
We furthermore have a variety of constructions that,
given such polynomials, can create other ones with the same property:
among these are algebraic analogues of the matroid operations
of deletion, contraction, direct sum, parallel connection,
series connection and 2-sum.
The complete monotonicity of $P^{-\beta}$ for some $\beta > 0$
can be viewed as a strong quantitative version of
the half-plane property (Hurwitz stability) for $P$,
and is also related to the Rayleigh property for matroids.
\end{abstract}


\bigskip
\noindent
{\bf Mathematics Subject Classification (MSC 2000) codes:}
05C31 (Primary);
05A15, 05A20, 05B35, 05C05, 05C50, 05E99, 15A15, 15B33, 15B57, 17C99,
26A48, 26B25, 26C05, 32A99,
43A85, 44A10, 60C05,
82B20 (Secondary).

\clearpage

\newtheorem{theorem}{Theorem}[section]
\newtheorem{proposition}[theorem]{Proposition}
\newtheorem{lemma}[theorem]{Lemma}
\newtheorem{corollary}[theorem]{Corollary}
\newtheorem{definition}[theorem]{Definition}
\newtheorem{conjecture}[theorem]{Conjecture}
\newtheorem{question}[theorem]{Question}
\newtheorem{problem}[theorem]{Problem}
\newtheorem{example}[theorem]{Example}

\renewcommand{\theenumi}{\alph{enumi}}
\renewcommand{\labelenumi}{(\theenumi)}
\def\eop{\hbox{\kern1pt\vrule height6pt width4pt
depth1pt\kern1pt}\medskip}
\def\prf{\par\noindent{\bf Proof.\enspace}\rm}
\def\rmk{\par\medskip\noindent{\bf Remark\enspace}\rm}

\newcommand{\be}{\begin{equation}}
\newcommand{\ee}{\end{equation}}
\newcommand{\<}{\langle}
\renewcommand{\>}{\rangle}
\newcommand{\widebar}{\overline}
\def\reff#1{(\protect\ref{#1})}
\def\spose#1{\hbox to 0pt{#1\hss}}
\def\ltapprox{\mathrel{\spose{\lower 3pt\hbox{$\mathchar"218$}}
    \raise 2.0pt\hbox{$\mathchar"13C$}}}
\def\gtapprox{\mathrel{\spose{\lower 3pt\hbox{$\mathchar"218$}}
    \raise 2.0pt\hbox{$\mathchar"13E$}}}
\def\textprime{${}^\prime$}
\def\proof{\par\medskip\noindent{\sc Proof.\ }}
\def\firstproof{\par\medskip\noindent{\sc First Proof.\ }}
\def\secondproof{\par\medskip\noindent{\sc Second Proof.\ }}
\def\qed{ $\square$ \bigskip}
\def\proofof#1{\bigskip\noindent{\sc Proof of #1.\ }}
\def\firstproofof#1{\bigskip\noindent{\sc First Proof of #1.\ }}
\def\secondproofof#1{\bigskip\noindent{\sc Second Proof of #1.\ }}
\def\altproofof#1{\bigskip\noindent{\sc Alternate Proof of #1.\ }}
\def\half{ {1 \over 2} }
\def\third{ {1 \over 3} }
\def\twothird{ {2 \over 3} }
\def\smfrac#1#2{{\textstyle{#1\over #2}}}
\def\smsmfrac#1#2{{\scriptstyle{#1\over #2}}}
\def\smhalf{ \smfrac{1}{2} }
\newcommand{\real}{\mathop{\rm Re}\nolimits}
\renewcommand{\Re}{\mathop{\rm Re}\nolimits}
\newcommand{\imag}{\mathop{\rm Im}\nolimits}
\renewcommand{\Im}{\mathop{\rm Im}\nolimits}
\newcommand{\sgn}{\mathop{\rm sgn}\nolimits}
\newcommand{\tr}{\mathop{\rm tr}\nolimits}
\newcommand{\per}{\mathop{\rm per}\nolimits}
\newcommand{\supp}{\mathop{\rm supp}\nolimits}
\newcommand{\diag}{\mathop{\rm diag}\nolimits}
\newcommand{\pf}{\mathop{\rm pf}\nolimits}
\newcommand{\Det}{\mathop{\rm Det}\nolimits}
\newcommand{\Mdet}{\mathop{\rm Mdet}\nolimits}
\newcommand{\Jdet}{\mathop{\rm Jdet}\nolimits}
\def\hboxscript#1{ {\hbox{\scriptsize\em #1}} }
\renewcommand{\emptyset}{\varnothing}

\newcommand{\restrict}{\upharpoonright}
\newcommand{\implies}{\;\Longrightarrow\;}

\newcommand{\scra}{{\mathcal{A}}}
\newcommand{\scrb}{{\mathcal{B}}}
\newcommand{\scrc}{{\mathcal{C}}}
\newcommand{\scrchat}{{\mathcal{C}^{\,\widehat{\,}}\,}}
\newcommand{\scrd}{{\mathcal{D}}}
\newcommand{\scrf}{{\mathcal{F}}}
\newcommand{\scrg}{{\mathcal{G}}}
\newcommand{\scrh}{{\mathcal{H}}}
\newcommand{\scrk}{{\mathcal{K}}}
\newcommand{\scrl}{{\mathcal{L}}}
\newcommand{\scrm}{{\mathcal{M}}}
\newcommand{\scro}{{\mathcal{O}}}
\newcommand{\scrp}{{\mathcal{P}}}
\newcommand{\scrr}{{\mathcal{R}}}
\newcommand{\scrs}{{\mathcal{S}}}
\newcommand{\scrt}{{\mathcal{T}}}
\newcommand{\scrv}{{\mathcal{V}}}
\newcommand{\scrw}{{\mathcal{W}}}
\newcommand{\scrz}{{\mathcal{Z}}}

\newcommand{\ahat}{{\widehat{a}}}
\newcommand{\Zhat}{{\widehat{Z}}}
\newcommand{\cc}{{\mathbf{c}}}
\renewcommand{\k}{{\mathbf{k}}}
\newcommand{\m}{{\mathbf{m}}}
\newcommand{\n}{{\mathbf{n}}}
\newcommand{\vv}{{\mathbf{v}}}
\newcommand{\bv}{{\mathbf{v}}}
\newcommand{\bz}{{\mathbf{z}}}
\newcommand{\w}{{\mathbf{w}}}
\newcommand{\x}{{\mathbf{x}}}
\newcommand{\y}{{\mathbf{y}}}
\newcommand{\z}{{\bf z}}
\newcommand{\zbar}{{\bar{\mathbf{z}}}}
\newcommand{\zero}{{\mathbf{0}}}
\newcommand{\one}{{\mathbf{1}}}
\newcommand{\bR}{{\bf R}}  
\newcommand{\bRtilde}{{\widetilde{\bf R}}}
\newcommand{\bRhat}{{\widehat{\bf R}}}

\newcommand{\ofo}{ {{}_1 \! F_1} }
\newcommand{\oft}{ {{}_1 \! F_2} }

\newcommand{\mF}{{\mathcal F}}
\newcommand{\mG}{{\mathcal G}}
\newcommand{\mH}{{\mathcal H}}

\newcommand{\C}{{\mathbb C}}
\newcommand{\D}{{\mathbb D}}
\newcommand{\Z}{{\mathbb Z}}
\newcommand{\N}{{\mathbb N}}
\newcommand{\Q}{{\mathbb Q}}
\newcommand{\R}{{\mathbb R}}
\newcommand{\HH}{{\mathbb H}}
\newcommand{\OO}{{\mathbb O}}
\newcommand{\RR}{{\mathbb R}}

\newcommand{\varphibar}{{\bar{\varphi}}}
\newcommand{\bvarphi}{{\boldsymbol{\varphi}}}
\newcommand{\bvarphibar}{{\bar{\boldsymbol{\varphi}}}}
\newcommand{\psibar}{{\bar{\psi}}}

\def\cbar{{\overline{C}}}

\newcommand{\bigdash}{%
\smallskip\begin{center} \rule{5cm}{0.1mm} \end{center}\smallskip}


\newenvironment{sarray}{
             \textfont0=\scriptfont0
             \scriptfont0=\scriptscriptfont0
             \textfont1=\scriptfont1
             \scriptfont1=\scriptscriptfont1
             \textfont2=\scriptfont2
             \scriptfont2=\scriptscriptfont2
             \textfont3=\scriptfont3
             \scriptfont3=\scriptscriptfont3
           \renewcommand{\arraystretch}{0.7}
           \begin{array}{l}}{\end{array}}

\newenvironment{scarray}{
             \textfont0=\scriptfont0
             \scriptfont0=\scriptscriptfont0
             \textfont1=\scriptfont1
             \scriptfont1=\scriptscriptfont1
             \textfont2=\scriptfont2
             \scriptfont2=\scriptscriptfont2
             \textfont3=\scriptfont3
             \scriptfont3=\scriptscriptfont3
           \renewcommand{\arraystretch}{0.7}
           \begin{array}{c}}{\end{array}}

\clearpage

\tableofcontents

\clearpage

\section{Introduction}  \label{sec1}

If $P$ is a univariate or multivariate polynomial with
real coefficients and strictly positive constant term,
and $\beta$ is a positive real number,
it is sometimes of interest to know whether $P^{-\beta}$
has all nonnegative (or even strictly positive) Taylor coefficients.
A problem of this type arose in the late 1920s in
Friedrichs and Lewy's study of the discretized time-dependent wave equation
in two space dimensions:
they needed the answer for the case
$P(y_1,y_2,y_3) = (1-y_1)(1-y_2) + (1-y_1)(1-y_3) + (1-y_2)(1-y_3)$
and $\beta = 1$.
Lewy contacted Gabor Szeg\H{o},
who proceeded to solve a generalization of this problem:
Szeg\H{o} \cite{Szego_33} showed that for any $n \ge 1$, the polynomial
\be
   P_n(y_1,\ldots,y_n)  \;=\; \sum_{i=1}^n \prod_{j \neq i} (1-y_j)
 \label{eq1.1}
\ee
has the property that $P_n^{-\beta}$ has nonnegative Taylor coefficients
for all $\beta \ge 1/2$.
(The cases $n=1,2$ are of course trivial;
 the interesting problem is for $n \ge 3$.)
Szeg\H{o}'s proof was surprisingly indirect,
and exploited the Gegenbauer--Sonine addition theorem for Bessel functions
together with Weber's first exponential integral.\footnote{
   These formulae for Bessel functions can be found in
   \cite[p.~367, eq.~11.41(17)]{Watson_44}
   and \cite[p.~394, eq.~13.3(4)]{Watson_44}, respectively.
   For the special case $n=3$,
   Szeg\H{o} also gave a version of the proof using
   Sonine's integral for the product of three Bessel functions
   \cite[p.~411, eq.~13.46(3)]{Watson_44}.
   Szeg\H{o} commented in his introduction \cite[p.~674]{Szego_33} that
   ``Die angewendeten Hilfsmittel stehen allerdings in keinem Verh\"altnis
     zu der Einfachheit des Satzes.''
   (``The tools used are, however, disproportionate to the
      simplicity of the result.'')

   Szeg\H{o} in fact proved the {\em strict}\/ positivity
   of the Taylor coefficients for all $n$ when $\beta = 1$,
   and for $n > 4\beta/(2\beta-1)$ when $\beta > 1/2$
   \cite[S\"atze~I--III]{Szego_33}.
   In this paper we shall concentrate on nonnegativity
   and shall not worry about whether strict positivity holds or not.
   But see Remark~1 after Theorem~\ref{thm.multidim.HBW}.
}
Shortly thereafter, Kaluza \cite{Kaluza_33} provided an elementary
(albeit rather intricate) proof, but only for $n=3$ and $\beta=1$.
In the early 1970s, Askey and Gasper \cite{Askey_72}
gave a partially alternate proof,
using Jacobi polynomials in place of Bessel functions.
Finally, Straub \cite{Straub_08} has very recently produced
simple and elegant proofs for the cases $n=3,4$ and $\beta=1$,
based on applying a positivity-preserving operator to another
rational function whose Taylor coefficients are known
(by a different elementary argument) to be nonnegative
(indeed strictly positive).

Askey and Gasper, in discussing both Szeg\H{o}'s problem and a related
unsolved problem of Lewy and Askey,
expressed the hope that
``there should be a combinatorial interpretation of these results''
and observed that
``this might suggest new methods'' \cite[p.~340]{Askey_72}.
The purpose of the present paper is to provide
such a combinatorial interpretation,
together with new and elementary (but we think powerful) methods of proof.
As a consequence we are able to prove a far-reaching generalization of
Szeg\H{o}'s original result, which includes as a special case
an affirmative solution to the problem of Lewy and Askey.
Indeed, we give {\em two different}\/ proofs for the Lewy--Askey problem,
based on viewing it as a member of two different families of generalizations
of the $n=3$ Szeg\H{o} problem.
Our methods turn out to be closely connected with
harmonic analysis on Euclidean Jordan algebras
(or equivalently on symmetric cones) \cite{Faraut_94}.

\subsection{Spanning-tree polynomials and series-parallel graphs}

{}From a combinatorial point of view,
one can see that Szeg\H{o}'s polynomial \reff{eq1.1}
is simply the spanning-tree generating polynomial
$T_G({\bf x})$ for the $n$-cycle $G=C_n$,
\be
   T_{C_n}(x_1,\ldots,x_n)  \;=\; \sum_{i=1}^n \prod_{j \neq i} x_j
   \;,
 \label{eq.TG.Cn}
\ee
after the change of variables $x_i = 1-y_i$.
This suggests that an analogous result might hold for
the spanning-tree polynomials of some wider class of graphs.\footnote{
   See \reff{def.TG.0}/\reff{def.TG} below for the general definition
   of the spanning-tree polynomial $T_G({\bf x})$
   for a connected graph $G$ \cite{Choe_hpp,Sokal_bcc2005}.
}
This conjecture is indeed true, as we shall show.
Moreover (and this will turn out to be quite important in what follows),
the change of variables $x_i = 1 -y_i$
can be generalized to $x_i = c_i - y_i$
for constants $c_i > 0$ that are not necessarily equal.
We shall prove:

\begin{theorem}
  \label{thm1.serpar}
Let $G=(V,E)$ be a connected series-parallel graph,
and let $T_G({\bf x})$ be its spanning-tree polynomial
in the variables ${\bf x} = \{x_e\}_{e \in E}$.
Then, for all $\beta \ge 1/2$ and
all choices of strictly positive constants ${\bf c} = \{c_e\}_{e \in E}$,
the function $T_G({\bf c} - {\bf y})^{-\beta}$
has nonnegative Taylor coefficients in the variables ${\bf y}$.

Conversely, if $G$ is a connected graph and
there exists $\beta \in (0,1) \setminus \{\smhalf\}$
such that $T_G({\bf c} - {\bf y})^{-\beta}$
has nonnegative Taylor coefficients (in the variables ${\bf y}$)
for all ${\bf c} > \zero$,
then $G$ is series-parallel.
\end{theorem}

The proof of the direct half of Theorem~\ref{thm1.serpar}
is completely elementary (and indeed quite simple).
The converse relies on a deep result from harmonic analysis on
Euclidean Jordan algebras
\cite{Gindikin_75,Ishi_00} \cite[Chapter~VII]{Faraut_94},
for which, however, there now exist two different elementary proofs
\cite{Shanbhag_88,Casalis_94} \cite{Sokal_riesz}.

Let us recall that a $C^\infty$ function $f(x_1,\ldots,x_n)$
defined on $(0,\infty)^n$ is termed {\em completely monotone}\/
if its partial derivatives of all orders alternate in sign, i.e.
\be
   (-1)^k {\partial^k f \over \partial x_{i_1} \,\cdots\, \partial x_{i_k}}
   \;\ge\;0 
 \label{eq.compmono}
\ee
everywhere on $(0,\infty)^n$,
for all $k \ge 0$ and all choices of indices $i_1,\ldots,i_k$.
Theorem \ref{thm1.serpar} can then be rephrased as follows:

\addtocounter{theorem}{-1}
\begin{theorem}
\hspace*{-3mm} ${}^{\bf\prime}$ \hspace{1mm}
Let $G=(V,E)$ be a connected series-parallel graph,
and let $T_G({\bf x})$ be its spanning-tree polynomial.
Then $T_G^{-\beta}$ is completely monotone on $(0,\infty)^E$
for all $\beta \ge 1/2$.

Conversely, if $G=(V,E)$ is a connected graph and
there exists $\beta \in (0,1) \setminus \{\smhalf\}$
such that $T_G^{-\beta}$ is completely monotone on $(0,\infty)^E$,
then $G$ is series-parallel.
\end{theorem}

Allowing arbitrary constants ${\bf c} > \zero$
thus allows the result to be formulated
in terms of complete monotonicity,
and leads to a characterization that is both necessary and sufficient.
Szeg\H{o}'s result (or rather, its generalization to arbitrary ${\bf c}$)
extends to series-parallel graphs and no farther.

\subsection{Determinants}

But this is not the end of the matter:
we can go far beyond series-parallel graphs
if we relax our demands about
the set of $\beta$ for which $T_G^{-\beta}$
is asserted to be completely monotone.
The key here is Kirchhoff's matrix-tree theorem
\cite{Kirchhoff_1847,Brooks_40,Nerode_61,Moon_70,Chen_76,%
Chaiken_78,Chaiken_82,Zeilberger_85,Moon_94,Abdesselam_03},
which shows how spanning-tree polynomials can be written as determinants.
This line of thought suggests that complete monotonicity of $P^{-\beta}$
might hold more generally for the homogeneous multiaffine polynomials
arising from determinants of the type studied in
\cite[Section~8.1]{Choe_hpp}.
This too is true;
in fact, such a result holds for a slightly more general class
of polynomials that need not be multiaffine.
We shall prove, once again by elementary methods:

\begin{theorem}
  \label{thm1.det}
Let $A_1,\ldots,A_n$ $(n \ge 1)$ be $m \times m$ real or complex matrices
or hermitian quaternionic matrices,
and let us form the polynomial
\be
   P(x_1,\ldots,x_n)  \;=\; \det\!\left( \sum_{i=1}^n x_i A_i \right)
 \label{def.P.det}
\ee
in the variables ${\bf x} = (x_1,\ldots,x_n)$.
[In the quaternionic case, $\det$ denotes the Moore determinant:
 see Appendix~\ref{sec.moore}.]
Assume further that there exists a linear combination of $A_1,\ldots,A_n$
that has rank $m$ (so that $P \not\equiv 0$).
\begin{itemize}
   \item[(a)]  If $A_1,\ldots,A_n$ are real symmetric positive-semidefinite
      matrices, then $P^{-\beta}$ is completely monotone on $(0,\infty)^n$
      for $\beta = 0,\smhalf,1,\smfrac{3}{2},\ldots$
      and for all real $\mbox{$\beta \ge (m-1)/2$}$.
   \item[(b)]  If $A_1,\ldots,A_n$ are complex hermitian positive-semidefinite
      matrices, then $P^{-\beta}$ is completely monotone on $(0,\infty)^n$
      for $\beta=0,1,2,3,\ldots$
      and for all real $\beta \ge m-1$.
   \item[(c)]  If $A_1,\ldots,A_n$ are quaternionic hermitian
      positive-semidefinite
      matrices, then $P^{-\beta}$ is completely monotone on $(0,\infty)^n$
      for $\beta=0,2,4,6,\ldots$
      and for all real $\beta \ge 2m-2$.
\end{itemize}
\end{theorem}

These curious conditions on $\beta$ are not just an artifact
of our method of proof; they really are best possible.
They can be better understood
if we take a slightly more general perspective,
and define complete monotonicity for functions
on an arbitrary open convex cone $C$
in a finite-dimensional real vector space $V$
(see Section~\ref{sec.compmono}).
We then have the following result that ``explains''
Theorem~\ref{thm1.det}:

\begin{theorem}
  \label{thm1.det.cones}
\begin{itemize}
   \item[(a)]  Let $V$ be the real vector space ${\rm Sym}(m,\R)$
      of real symmetric $m \times m$ matrices,
      and let $C \subset V$ be the cone $\Pi_m(\R)$
      of positive-definite matrices.
      Then the map $A \mapsto (\det A)^{-\beta}$
      is completely monotone on $C$
      if and only if
      $\beta \in \{0,\smhalf,1,\smfrac{3}{2},\ldots\} \cup [(m-1)/2,\infty)$.
      Indeed, if
      $\beta \notin \{0,\smhalf,1,\smfrac{3}{2},\ldots\} \cup [(m-1)/2,\infty)$,
      then the map $A \mapsto (\det A)^{-\beta}$
      is not completely monotone on any
      nonempty open convex subcone $C' \subseteq C$.
   \item[(b)]  Let $V$ be the real vector space ${\rm Herm}(m,\C)$
      of complex hermitian $m \times m$ matrices,
      and let $C \subset V$ be the cone $\Pi_m(\C)$
      of positive-definite matrices.
      Then the map $A \mapsto (\det A)^{-\beta}$
      is completely monotone on $C$
      if and only if
      $\beta \in \{0,1,2,3,\ldots\} \cup [m-1,\infty)$.
      Indeed, if
      $\beta \notin \{0,1,2,3,\ldots\} \cup [m-1,\infty)$,
      then the map $A \mapsto (\det A)^{-\beta}$
      is not completely monotone on any
      nonempty open convex subcone $C' \subseteq C$.
   \item[(c)]  Let $V$ be the real vector space ${\rm Herm}(m,\HH)$
      of quaternionic hermitian $m \times m$ matrices,
      and let $C \subset V$ be the cone $\Pi_m(\HH)$
      of positive-definite matrices.
      Then the map $A \mapsto (\det A)^{-\beta}$
      is completely monotone on $C$
      if and only if
      $\beta \in \{0,2,4,6,\ldots\} \cup [2m-2,\infty)$.
      Indeed, if
      $\beta \notin \{0,2,4,6,\ldots\} \cup [2m-2,\infty)$,
      then the map $A \mapsto (\det A)^{-\beta}$
      is not completely monotone on any
      nonempty open convex subcone $C' \subseteq C$.
\end{itemize}
\end{theorem}

\noindent
In particular, if the matrices $A_1,\ldots,A_n$
together span ${\rm Sym}(m,\R)$, ${\rm Herm}(m,\C)$ or  ${\rm Herm}(m,\HH)$
[so that the convex cone they generate has nonempty interior],
then the determinantal polynomial \reff{def.P.det}
has $P^{-\beta}$ completely monotone on $(0,\infty)^n$
if and {\em only~if}\/ $\beta$ belongs to the set enumerated
in Theorem~\ref{thm1.det}.\footnote{
   Br\"and\'en \cite{Branden_12} has very recently used this latter fact
   to determine the exact set of $\alpha \in \R$ for which the
   $\alpha$-permanent \cite{Vere-Jones_88}
   is nonnegative on real symmetric (resp.\ complex hermitian)
   positive-semidefinite matrices.
}

The proof of the ``if'' part of Theorem~\ref{thm1.det.cones}
is completely elementary, but the ``only if'' part
again relies on a deep result from harmonic analysis on
Euclidean Jordan algebras,
namely, the characterization of parameters for which the
Riesz distribution is a positive measure
(Theorem~\ref{thm.riesz1} below;
 but see \cite{Shanbhag_88,Casalis_94} \cite{Sokal_riesz}
 and Appendix~\ref{sec.gindikin} below
 for elementary proofs).
In fact, when Theorem~\ref{thm1.det.cones} is rephrased
in this latter context it takes on a unified form:

\begin{theorem}
  \label{thm1.det.cones.Jordan}
Let $V$ be a simple Euclidean Jordan algebra of dimension~$n$ and rank~$r$,
with $n = r + \frac{d}{2} r(r-1)$;
let $\Omega \subset V$ be the positive cone;
and let $\Delta \colon\, V \to \R$ be the Jordan determinant.
Then the map $x \mapsto \Delta(x)^{-\beta}$ is completely monotone on $\Omega$
if and only if
$\beta \in \{0,\frac{d}{2},\ldots,(r-1)\frac{d}{2}\}$
or $\beta > (r-1)\frac{d}{2}$.
Indeed, if $\beta \notin \{0,\frac{d}{2},\ldots,(r-1)\frac{d}{2}\} \cup
   ((r-1)\frac{d}{2},\infty)$,
then the map $x \mapsto \Delta(x)^{-\beta}$
is not completely monotone on any
nonempty open convex subcone $\Omega' \subseteq \Omega$.
\end{theorem}

\noindent
We shall see that
Theorem~\ref{thm1.det.cones.Jordan} is essentially {\em equivalent}\/
to the characterization of parameters for which the
Riesz distribution is a positive measure.
The set of values of $\beta$ described in Theorem~\ref{thm1.det.cones.Jordan}
is known as the {\em Gindikin--Wallach set}\/
and arises in a number of contexts in representation theory
\cite{Berezin_75,Gindikin_75,Rossi_76,Wallach_79,Lassalle_87,%
Faraut_88,Faraut_90,Faraut_94}.

A special case of the construction \reff{def.P.det}
arises \cite[Section~8.1]{Choe_hpp}
when $B$ is an $m \times n$ real or complex matrix of rank $m$,
and we set $P({\bf x}) = \det(BXB^*)$,
where $X = \diag(x_1,\ldots,x_n)$ and ${}^*$ denotes hermitian conjugate.
Then the matrix $A_i$ in \reff{def.P.det} is simply the outer product of the
$i$th column of $B$ with its complex conjugate,
and so is of rank at most 1;
as a consequence, the polynomial $P$ is multiaffine
(i.e., of degree at most 1 in each variable separately).\footnote{
   See Proposition~\ref{prop.detpoly2} below.
}

In particular, let $G=(V,E)$ be a connected graph,
and define its spanning-tree polynomial $T_G({\bf x})$ by
\be
   T_G({\bf x})  \;=\;   \sum_{T \in \scrt(G)} \, \prod_{e \in T} x_e
   \;,
 \label{def.TG.0}
\ee
where ${\bf x} = \{x_e\}_{e \in E}$
is a family of indeterminates indexed by the edges of $G$,
and $\scrt(G)$ denotes the family of edge sets of spanning trees in $G$.
Now let $B$ be the directed vertex-edge incidence matrix
for an arbitrarily chosen orientation of $G$,
with one row (corresponding to an arbitrarily chosen vertex of $G$) deleted;
then the matrix-tree theorem
\cite{Kirchhoff_1847,Brooks_40,Nerode_61,Moon_70,Chen_76,%
Chaiken_78,Chaiken_82,Zeilberger_85,Moon_94,Abdesselam_03,Choe_hpp}
tells us that $T_G({\bf x}) = \det(BXB^{\rm T})$.
Applying Theorem~\ref{thm1.det}(a), we obtain:

\begin{corollary}
   \label{cor.TG}
Let $G=(V,E)$ be a connected graph with $p$ vertices,
and let $T_G({\bf x})$ be its spanning-tree polynomial.
Then $T_G^{-\beta}$ is completely monotone on $(0,\infty)^E$
for $\beta = 0,\smhalf,1,\smfrac{3}{2},\ldots$
and for all real $\beta \ge (p-2)/2$.
\end{corollary}

Likewise, we can apply Theorem~\ref{thm1.det}(b)
to the elementary symmetric polynomial
\be
   E_{2,4}(x_1,x_2,x_3,x_4)  \;=\;
   x_1 x_2 + x_1 x_3 + x_1 x_4 + x_2 x_3 + x_2 x_4 + x_3 x_4
   \;,
 \label{eq.E24}
\ee
which can be represented in the form \reff{def.P.det} with
\be
   A_1 = \left(\!\! \begin{array}{cc}
                       1  & 0  \\
                       0  & 0
                    \end{array}
         \!\!\right)
   ,\quad
   A_2 = \left(\!\! \begin{array}{cc}
                       0  & 0  \\
                       0  & 1
                    \end{array}
         \!\!\right)
   ,\quad
   A_3 = \left(\!\! \begin{array}{cc}
                       1  & 1  \\
                       1  & 1
                    \end{array}
         \!\!\right)
   ,\quad
   A_4 = \left(\!\! \begin{array}{cc}
                       1           & e^{-i\pi/3}  \\
                       e^{i\pi/3}  & 1
                    \end{array}
         \!\!\right)
 \label{eq.matrices.E24}
\ee
or equivalently as $E_{2,4}({\bf x}) = \det(BXB^*)$ with
$B = \displaystyle{ \left(\!\! \begin{array}{cccc}
                                   1  & 0  &  1  &  1  \\
                                   0  & 1  &  1  & e^{i\pi/3}
                               \end{array}
                    \!\!\!\right)}$.
We obtain:

\begin{corollary}
   \label{cor.E24}
The function $E_{2,4}^{-\beta}$ is completely monotone on $(0,\infty)^4$
if and only if $\beta = 0$ or $\beta \ge 1$.
In particular, the function
\be
   \left( \sum\limits_{1 \le i < j \le 4} (1-y_i) (1-y_j) \right) ^{\! -\beta}
 \label{eq.lewy}
\ee
has nonnegative Taylor coefficients for all $\beta \ge 1$.
\end{corollary}

\noindent
Indeed, the ``if'' part can be corroborated by
an explicit Laplace-transform formula for $E_{2,4}^{-\beta}$ for $\beta > 1$:
see \reff{eq.explicit.E24} below.
The ``only if'' follows from the observation made after
Theorem~\ref{thm1.det.cones},
since the matrices $A_1,\ldots,A_4$ in \reff{eq.matrices.E24}
span ${\rm Herm}(2,\C)$.

The second sentence of Corollary~\ref{cor.E24} answers in the affirmative
a question posed long ago by Lewy \cite[p.~340]{Askey_72},
of which Askey remarks that it ``has caused me many hours of frustration''
\cite[p.~56]{Askey_75}.\footnote{
   Askey \cite[p.~56]{Askey_75} comments that, in his view,
   \begin{quote}
   So far the most powerful method of treating problems of this type
   is to translate them into another problem involving special functions
   and then use the results and methods which have been developed
   for the last two hundred years to solve the special function problem.
   So far I have been unable to make a reduction in [Lewy's problem]
   and so have no place to start.
   \end{quote}
   But he immediately adds, wisely, that
   ``it is possible to solve some problems without using special functions,
     so others should not give up on [Lewy's problem].''
}
(See also the recent discussion in \cite{Kauers_08}.)
Indeed, Lewy's question concerned only $\beta=1$,
and made the weaker conjecture that the function \reff{eq.lewy}
multiplied by $(4-y_1-y_2-y_3-y_4)^{-1}$ has nonnegative Taylor coefficients.
This latter factor is now seen to be unnecessary.\footnote{
   Ismail and Tamhankar \cite[p.~483]{Ismail_79}
   mistakenly asserted that
   ``the early coefficients in the power series expansion of
   $$
   \{(1-r) (1-s) + (1-r) (1-t) + (1-r) (1-u) +
     (1-s) (1-t) + (1-s) (1-u) + (1-t) (1-u)\}^{-1}
   $$
   are positive but the later coefficients do change sign'',
   arguing that this is ``because Huygen's [{\em sic}\/] principle
   holds in three-space.''
   Huygens' principle indeed suggests that the coefficients
   approach zero, as Askey and Gasper \cite[p.~340]{Askey_72} observed;
   but this in no way contradicts the nonnegativity of those coefficients.
}

Similarly, Theorem~\ref{thm1.det}(c) applied to the
quaternionic determinant
$\det\Biggl(\!\! \begin{array}{cc}
                      a & q \\
                      \bar{q} & b
                 \end{array}
            \!\!\Biggr)    = ab - q \bar{q}$
for $a,b \in \R$ and $q \in \HH$,
with $A_1,\ldots,A_4$ as above and
\begin{subeqnarray}
   & &
   A_5 = \left(\!\! \begin{array}{cc}
                       1  &  q_5 \\
                       \bar{q}_5 & 1
                    \end{array}
         \!\!\right)
   ,\quad
   A_6 = \left(\!\! \begin{array}{cc}
                       1  & q_6 \\
                       \bar{q}_6 & 1
                    \end{array}
         \!\!\right)
   \\[3mm]
   & &
   q_5 \;=\; \frac{1}{2} - \frac{\sqrt{3}}{6} i
                         - \frac{\sqrt{6}}{3} j
   ,\quad
   q_6 \;=\; \frac{1}{2} - \frac{\sqrt{3}}{6} i
                         - \frac{\sqrt{6}}{12} j
                         - \frac{\sqrt{10}}{4} k
   \;,
\end{subeqnarray}
yields
an analogous result for the elementary symmetric polynomial
\be
   E_{2,6}(x_1,\ldots,x_6)  \;=\;
   x_1 x_2 + x_1 x_3 + \ldots + x_5 x_6
\ee
of degree 2 in six variables\footnote{
   If we define $q_3 = 1$ and $q_4 =  e^{-i\pi/3}$
   [cf.\ \reff{eq.matrices.E24}],
   then $q_3,q_4,q_5,q_6$ are quaternions satisfying
   $$ \real(q_i \bar{q}_j)  \;=\;
      \cases{1  &  if $i=j$  \cr
             \noalign{\vskip 3pt}
             1/2  & if $i \neq j$ \cr
            }
   $$
   or equivalently $|q_i|^2 = 1$ and $|q_i-q_j|^2 = 1$ for all $i \neq j$.
   From this it easily follows that
   $$ \det\left(\!
         \begin{array}{cc}
             x_1 + x_3 + x_4 + x_5 + x_6 &
                x_3 q_3 + x_4 q_4 + x_5 q_5 + x_6 q_6 \\
             x_3 \bar{q}_3 + x_4 \bar{q}_4 + x_5 \bar{q}_5 + x_6 \bar{q}_6 &
                x_2 + x_3 + x_4 + x_5 + x_6
         \end{array}
         \!\!\right)
      \;=\; 
      E_{2,6}(x_1,\ldots,x_6)  \;.
   $$
}:

\begin{corollary}
   \label{cor.E26}
The function $E_{2,6}^{-\beta}$ is completely monotone on $(0,\infty)^6$
if and only if $\beta = 0$ or $\beta \ge 2$.
\end{corollary}

Corollaries~\ref{cor.TG} and \ref{cor.E24} are in fact
special cases of a much more general result concerning
the basis generating polynomials $B_M({\bf x})$ of certain classes of matroids.
(We stress that no knowledge of matroid theory is needed to understand
the main arguments of this paper;  readers allergic to matroids,
or simply unfamiliar with them, can skip all references to them
without loss of logical continuity.
Still, we think that the matroidal perspective is fruitful
and we would like to make some modest propaganda for it.\footnote{
   See \cite{Oxley_11} for background on matroid theory,
   and \cite{Choe_hpp} for background on basis generating polynomials.
   In interpreting Corollary~\ref{cor1.det} below,
   please note that if $G=(V,E)$ is a graph with $k$ connected components,
   then the graphic matroid $M(G)$ has rank $|V|-k$,
   while the cographic matroid $M^*(G)$ has rank $|E| - |V| + k$.
   Note also that the equivalence of ``complex-unimodular matroid'' with
   ``sixth-root-of-unity matroid'' is proven in
   \cite[Theorem~8.9]{Choe_hpp}.
})
So let $M$ be a matroid with ground set $E$,
and let $\scrb(M)$ be its set of bases;
then the {\em basis generating polynomial}\/ of $M$ is, by definition,
$B_M({\bf x}) = \sum\limits_{S \in \scrb(M)} x^S$,
where ${\bf x} = \{x_e\}_{e \in E}$ is a family of indeterminates
indexed by the elements of $M$,
and we have used the shorthand $x^S = \prod\limits_{e \in S} x_e$.

Now let $B$ be an arbitrary $m \times n$ real or complex matrix of rank $m$,
and define $P({\bf x}) = \det(BXB^*)$.
Then, as discussed previously, Theorem~\ref{thm1.det}(a or b)
applies to~$P$ and gives a sufficient condition for $P^{-\beta}$
to be completely monotone.
On the other hand, the Cauchy--Binet formula gives
\be
   P({\bf x}) \;=\;
   \det(BXB^*)  \;=\;
   \sum\limits_{\begin{scarray}
                    S \subseteq [n] \\
                    |S| = m
                 \end{scarray}}
   |\det B_{\star S}|^2 \, x^S
 \label{eq.cb}
\ee
where $B_{\star S}$ denotes the submatrix of $B$ with columns $S$.
Since $\det B_{\star S} \neq 0$ if and only if the columns $S$ of $B$
are linearly independent, we see that $P$ is a {\em weighted}\/ version
of the basis generating polynomial for the matroid $M = M[B]$
that is represented by $B$ (this matroid has rank $m$).
In particular, a matroid is said to be
{\em real-unimodular}\/ (resp.\ {\em complex-unimodular}\/)
if it has a real (resp.\ complex) representing matrix $B$,
with a number of rows equal to its rank,
such that $|\det B_{\star S}|^2 \in \{0,1\}$ for all $S$.\footnote{
   This is not the usual definition of real-unimodular/complex-unimodular,
   but it is equivalent to the usual definition by virtue of
   \cite[Proposition~8.6]{Choe_hpp}.
}
In this case the basis generating polynomial is precisely
$B_M({\bf x}) = \det(BXB^*)$.
We thereby obtain from Theorem~\ref{thm1.det}(a,b) the following result:

\begin{corollary}
   \label{cor1.det}
Let $M$ be a matroid of rank~$r$ on the ground set $E$,
and let $B_M({\bf x})$ be its basis generating polynomial.
\begin{itemize}
   \item[(a)]  If $M$ is a regular [= real-unimodular] matroid,
      then $B_M^{-\beta}$ is completely monotone on $(0,\infty)^E$
      for $\beta = 0,\smhalf,1,\smfrac{3}{2},\ldots$
      and for all real $\beta \ge (r-1)/2$.
      (This holds in particular if $M$ is a graphic or cographic matroid,
       i.e.\ for the spanning-tree or complementary-spanning-tree
       polynomial of a connected graph.)
   \item[(b)]  If $M$ is a complex-unimodular matroid
      [= sixth-root-of-unity matroid],
      then $B_M^{-\beta}$ is completely monotone on $(0,\infty)^E$
      for $\beta=0,1,2,3,\ldots$
      and for all real $\beta \ge r-1$.
\end{itemize}
\end{corollary}

\noindent
In particular, by specializing (a) to a graphic matroid $M(G)$
we recover Corollary~\ref{cor.TG},
and by specializing (b) to the uniform matroid $U_{2,4}$
we recover Corollary~\ref{cor.E24}.

We have also proven a (very) partial converse to Corollary~\ref{cor1.det},
which concerns the cases of rank-$r$ $n$-element simple matroids
in which the matrices $A_1,\ldots,A_n$
together span ${\rm Sym}(r,\R)$ or ${\rm Herm}(r,\C)$:
see Proposition~\ref{cor2.det} below.

\bigskip

{\bf Remark.}  There is also an analogue of Corollary~\ref{cor1.det}
in the quaternionic case.
Recall first the quaternionic analogue of the
Cauchy--Binet formula \reff{eq.cb}
[Proposition~\ref{prop.properties.Mdet}(g) below]:
if $B$ is an $m \times n$ quaternionic matrix, then
$P(x_1,\ldots,x_n) = \det(BXB^*) =
                     \det\Bigl( \sum\limits_{i=1}^n x_i A_i \Bigr)$
is well-defined for $X=\diag(x_1,\ldots,x_n)$
with $x_1,\ldots,x_n$ {\em real}\/
and equals the polynomial
$\sum\limits_{\begin{scarray}
                 S \subseteq [n] \\
                 |S| = m
              \end{scarray}}
 \det[B_{\star S} (B_{\star S})^*] \, x^S$.
(Note, by contrast, that $\det B_{\star S}$ is in general meaningless
 because $B_{\star S}$ need not be hermitian.)
We can then define a matroid $M$ to be {\em quaternionic-unimodular}\/
if its basis generating polynomial can be represented in this way,
i.e.\ if it has a quaternionic representing matrix $B$,
with a number of rows equal to its rank,
such that $\det[B_{\star S} (B_{\star S})^*] \in \{0,1\}$ for all $S$.
For such matroids $M$, Theorem~\ref{thm1.det}(c) implies that
$B_M^{-\beta}$ is completely monotone on $(0,\infty)^E$
for $\beta=0,2,4,6,\ldots$ and for all real $\beta \ge 2r-2$.

A deeper study of quaternionic-unimodular matroids would be of interest.
For instance, is the class of quaternionic-unimodular matroids
closed under duality?
Or even under contraction?
(The class is obviously closed under deletion.)
Which uniform matroids $U_{r,n}$ are quaternionic-unimodular?

A different notion of ``quaternionic-unimodular matroid''
has been introduced recently by Pendavingh and van Zwam \cite{Pendavingh_13}.
It is not clear to us what is the relation between their notion and ours.

\subsection{Quadratic forms}

Of course, $E_{2,4}$ and $E_{2,6}$ are quadratic forms
in the variables ${\bf x}$,
as is the polynomial $E_{2,3}$ arising in the $n=3$ Szeg\H{o} problem.
This suggests that it might be fruitful to study more general
quadratic forms.
We shall prove, by elementary methods:

\begin{theorem}
   \label{thm1.quadratic.cones}
Let $V$ be a finite-dimensional real vector space,
let $B$ be a symmetric bilinear form on $V$
having inertia $(n_+,n_-,n_0)$,
and define the quadratic form $Q(x) = B(x,x)$.
Let $C \subset V$ be a nonempty open convex cone
with the property that $Q(x) > 0$ for all $x \in C$.
Then $n_+ \ge 1$, and moreover:
\begin{itemize}
   \item[(a)]  If $n_+ = 1$ and $n_- = 0$, then
       $Q^{-\beta}$
       is completely monotone on $C$ for all $\beta \ge 0$.
       For all other values of $\beta$,
       $Q^{-\beta}$
       is not completely monotone on any nonempty open convex subcone
       $C' \subseteq C$.
   \item[(b)]  If $n_+ = 1$ and $n_- \ge 1$, then
       $Q^{-\beta}$
       is completely monotone on $C$ for $\beta = 0$
       and for all $\beta \ge (n_- -1)/2$.
       For all other values of $\beta$,
       $Q^{-\beta}$
       is not completely monotone on any nonempty open convex subcone
       $C' \subseteq C$.
   \item[(c)]  If $n_+ > 1$, then
       $Q^{-\beta}$
       is not completely monotone on any nonempty open convex subcone
       $C' \subseteq C$ for any $\beta \neq 0$.
\end{itemize}
\end{theorem}

Theorem~\ref{thm1.quadratic.cones} follows
fairly easily from the classic work
of Marcel Riesz \cite{Riesz_49}
(see also \cite{Duistermaat_91} and \cite[Chapter~VII]{Faraut_94})
treating the case in which $B$ is the Lorentz form on $\R^n$,
\be
   B(x,y)  \;=\;  x_1 y_1 - x_2 y_2 - \ldots - x_n y_n
   \;,
\ee
and $C$ is the Lorentz cone (= forward light cone)
$\{x \in \R^n \colon\, x_1 > \sqrt{x_2^2 + \ldots + x_n^2} \}$.
We are able to give a completely elementary proof
of both the sufficiency and the necessity;
and we are able to give in case (b)
an explicit Laplace-transform formula for $Q^{-\beta}$
(see Proposition~\ref{prop.quadratic.A}).

Specializing Theorem~\ref{thm1.quadratic.cones}
with $V = \R^n$ and $C = (0,\infty)^n$
to the degree-2 elementary symmetric polynomials
\be
   E_{2,n}(x_1,\ldots,x_n)  \;=\;  \sum_{1 \le i < j \le n} x_i x_j
   \;,
\ee
we obtain:

\begin{corollary}
   \label{cor.quadratic}
The function $E_{2,n}^{-\beta}$
is completely monotone on $(0,\infty)^n$ if and only if
$\beta = 0$ or $\beta \ge (n-2)/2$.
\end{corollary}

\noindent
By this method we obtain
an alternate proof of Corollaries~\ref{cor.E24} and \ref{cor.E26} ---
hence in particular a second solution to the Lewy--Askey problem ---
as well as of Szeg\H{o}'s \cite{Szego_33} original result
in the case $n=3$.\footnote{
   In fancy language --- which is, however, completely unnecessary
   for understanding our proofs --- our ``determinantal'' proof of
   Corollary~\ref{cor.E24} is based on harmonic analysis
   on the cone of positive-definite $m \times m$ complex hermitian matrices
   specialized to $m=2$, while our ``quadratic form'' proof is
   based on harmonic analysis on the Lorentz cone in $\R^n$
   specialized to $n=4$.
   The point here is that the Jordan algebra
   ${\rm Herm}(2,\C) \simeq \R \times \R^3$
   can be viewed as a member of two different families of Jordan algebras:
   ${\rm Herm}(m,\C)$ and $\R \times \R^{n-1}$ \cite[p.~98]{Faraut_94}.
   Likewise, our ``determinantal'' proof of
   Corollary~\ref{cor.E26} is based on harmonic analysis
   on the cone of positive-definite $m \times m$ quaternionic hermitian matrices
   specialized to $m=2$, while our ``quadratic form'' proof is
   based on harmonic analysis on the Lorentz cone in $\R^n$
   specialized to $n=6$;
   and we have the isomorphism of Jordan algebras
   ${\rm Herm}(2,\HH) \simeq \R \times \R^5$ \cite[p.~98]{Faraut_94}.
   And finally, our ``determinantal'' proof of
   the $n=3$ Szeg\H{o} result is based on harmonic analysis
   on the cone of positive-definite $m \times m$ real symmetric matrices
   specialized to $m=2$, while our ``quadratic form'' proof is
   based on harmonic analysis on the Lorentz cone in $\R^n$
   specialized to $n=3$;
   and we have ${\rm Sym}(2,\R) \simeq \R \times \R^2$ \cite[p.~98]{Faraut_94}.
}
We also obtain an explicit Laplace-transform formula for $E_{2,n}^{-\beta}$
(see Corollary~\ref{cor.quadratic.E2n}).


\bigskip

{\bf Remark.}
It is easy to see that $E_{2,n}$ is the spanning-tree polynomial
of a graph only if $n=2$ or 3:
a connected graph $G$ whose spanning-tree polynomial is of degree~2
must have precisely three vertices;
if $G$ has multiple edges, then $T_G \neq E_{2,n}$
because monomials corresponding to pairs of parallel edges are absent
from $T_G$;
so $G$ must be either the 3-vertex path or the 3-cycle,
corresponding to $E_{2,2}$ or $E_{2,3}$, respectively.
But this fact can also be seen from our results:
Corollary~\ref{cor.TG} says that $T_G^{-1/2}$ is
completely monotone for all graphs $G$,
while Corollary~\ref{cor.quadratic}
says that $E_{2,n}^{-1/2}$ is not completely monotone when $n > 3$.

\bigskip

Corollaries~\ref{cor.E24}, \ref{cor.E26} and \ref{cor.quadratic}
lead naturally to the following question:
If we write $E_{r,n}$ for the elementary symmetric polynomial of degree $r$
in $n$ variables,
\be
   E_{r,n}(x_1,\ldots,x_n)   \;=\;
   \sum_{1 \le i_1 < i_2 < \ldots < i_r \le n}  x_{i_1} x_{i_2} \cdots x_{i_r}
\ee
(where we set $E_{0,n} \equiv 1$),
then for which $\beta > 0$ is $E_{r,n}^{-\beta}$ completely monotone
on $(0,\infty)^n$?
The cases $r=0$, 1 and $n$ are trivial:
we have complete monotonicity for all $\beta \ge 0$.
Our results for the cases $r=n-1$
(Theorem~\ref{thm1.serpar}${}'$ specialized to cycles $C_n$)
and $r=2$ (Corollary~\ref{cor.quadratic}),
as well as numerical experiments for
$(r,n) = (3,5)$, (3,6) and (4,6),
%
%
%
%
%
%
%
%
%
%
%
%
%
%
%
%
%
%
suggest the following conjecture:

\begin{conjecture}
  \label{conj.Ern}
Let $2 \le r \le n$.
Then $E_{r,n}^{-\beta}$ is completely monotone on $(0,\infty)^n$
if and only if $\beta = 0$ or $\beta \ge (n-r)/2$.
\end{conjecture}

\noindent
However, we have been unable to find a proof
of either the necessity or the sufficiency.

We remark that the elementary symmetric polynomial $E_{r,n}$
is the basis generating polynomial of the uniform matroid $U_{r,n}$.
So Corollary~\ref{cor.quadratic} and Conjecture~\ref{conj.Ern}
concern the same general subject as Corollary~\ref{cor1.det},
namely, complete monotonicity for inverse powers
of the basis generating polynomials of matroids.

%
%
%
%

\subsection{Discussion}

In summary, we have two {\em ab initio}\/ methods for proving,
given a polynomial $P$ and a positive real number $\beta$,
that $P^{-\beta}$ is completely monotone on $(0,\infty)^n$
[or more generally on a convex cone $C$]:
\begin{itemize}
   \item[(a)]  The determinantal method
     (Theorems~\ref{thm1.det}, \ref{thm1.det.cones}
      and \ref{thm1.det.cones.Jordan}:
      see Section~\ref{sec.det}).
   \item[(b)]  The quadratic-form method
     (Theorem~\ref{thm1.quadratic.cones}:
      see Section~\ref{sec.quadratic}).
\end{itemize}
Interestingly, these two methods can be viewed as versions of
the {\em same}\/ construction, involving the determinant
on a Euclidean Jordan algebra and
the Laplace-transform representation of its inverse powers
\cite[Chapters~II--VII]{Faraut_94}.
We discuss this connection in
Sections~\ref{subsec.det.major} and \ref{subsec.quadratic.major}.

In addition to these two {\em ab initio}\/ methods,
we have a variety of constructions that,
given such polynomials, can create other ones with the same property
(see Section~\ref{sec.constructions}).
Among these are algebraic analogues of the graph (or matroid) operations
of deletion, contraction, direct sum\footnote{
   By ``direct sum'' of graphs we mean either
   disjoint union (``0-sum'') or gluing at a cut vertex (``1-sum'').
   Both of these operations correspond to the direct sum of matroids.
},
parallel connection, series connection and 2-sum (but {\em not}\/ duality).
By combining these operations with our {\em ab initio}\/ proofs,
we are able to prove the complete monotonicity of $T_G^{-\beta}$
for some values of $\beta$ beyond those covered by Corollary~\ref{cor.TG}:

\begin{proposition}
   \label{prop.TG.extended}
Fix $p \ge 2$,
and let $G=(V,E)$ be any graph that can be obtained from copies
of the complete graph $K_p$ by
parallel connection, series connection, direct sum, deletion and contraction.
Then $T_G^{-\beta}$ is completely monotone on $(0,\infty)^E$
for $\beta = 0,\smhalf,1,\smfrac{3}{2},\ldots$
and for all real $\beta \ge (p-2)/2$.
\end{proposition}

\noindent
In particular, the case $p=3$ covers series-parallel graphs;
this is essentially our proof of the direct half of
Theorem~\ref{thm1.serpar}${}'$.
We also have versions of this proposition for matroids:
see Propositions~\ref{prop.BM.extended1} and \ref{prop.BM.extended2} below.
Finally, in Propositions~\ref{prop.min3conn.graphs}
and \ref{prop.min3conn.matroids} we give excluded-minors
characterizations of the class of graphs/matroids handled by
Propositions~\ref{prop.TG.extended} and \ref{prop.BM.extended1},
respectively.

But even in the graphic case,
we are still far from having a complete answer
to the following fundamental problem:

\begin{problem}
   \label{problem.Gbeta}
Given a graph $G=(V,E)$, for which real numbers $\beta > 0$
is the function $T_G^{-\beta}$ completely monotone on $(0,\infty)^E$?
\end{problem}

\noindent
This question can be rephrased usefully as follows:

\addtocounter{theorem}{-1}
\begin{problem}
\hspace*{-3mm} ${}^{\bf\prime}$ \hspace{1mm}
For each $\beta > 0$, characterize the class $\scrg_\beta$
of graphs for which
$T_G^{-\beta}$ is completely monotone on $(0,\infty)^E$.
\end{problem}

\noindent
We will show in Section~\ref{sec.application.graphs}
that the class $\scrg_\beta$
is closed under minors --- so that it can be characterized by listing
the excluded minors --- and under parallel connection.
Furthermore, it is closed under series connection
when (but only when) $\beta \ge 1/2$.

In this paper we have solved Problem~\ref{problem.Gbeta}${}'$ in a few cases:
\begin{itemize}
   \item For $\beta \in \{{1 \over 2},1,{3 \over 2},\ldots\}$,
      $\scrg_\beta =$ all graphs.
      See Corollary~\ref{cor1.det}(a).
   \item For $\beta \in (0,{1 \over 2})$,
      $\scrg_\beta =$ graphs obtained from forests by parallel extension
      of edges (= graphs with no $K_3$ minor).
      See Theorem~\ref{thm.g0half}.
   \item For $\beta \in ({1 \over 2},1)$,
      $\scrg_\beta =$ series-parallel graphs (= graphs with no $K_4$ minor).
      See Theorems~\ref{thm1.serpar}${}'$ and \ref{thm.ghalf1}.
\end{itemize}
So the first unsolved cases are $\beta \in (1,{3 \over 2})$:
Might it be true that $\scrg_\beta =$ all graphs with no $K_5$ minor?
Or might there exist, alternatively, other excluded minors?
We have been thus far unable to determine the complete monotonicity
of $T_G^{-\beta}$ for the cases $G=W_4$ (the wheel with four spokes)
and $G = K_5 - e$ (the complete graph $K_5$ with one edge deleted).
Indeed, for $k \ge 2$ we do not even know the answer to the following question:

\begin{question}
Fix an integer $k \ge 0$.  Must we have $\scrg_\beta = \scrg_{\beta'}$
whenever $\beta,\beta' \in ({k \over 2}, {k+1 \over 2})$?
\end{question}

Let us mention, finally, an alternative approach to ``converse'' results
that we have not pursued, for lack of competence.
When $P^{-\beta}$ does {\em not}\/ have all nonnegative Taylor coefficients,
this fact should in most cases be provable either by explicit computation
of low-order coefficients or by asymptotic computation of suitable families
of high-order coefficients (computer experiments can usually suggest
which families to focus on).  This type of multivariate asymptotic
calculation has been pioneered recently by Pemantle and collaborators
\cite{Pemantle_02,Pemantle_04,Pemantle_08,Pemantle_10,Baryshnikov_11}
and involves some rather nontrivial algebraic geometry/topology.
In fact, Baryshnikov and Pemantle \cite[Section~4.4]{Baryshnikov_11}
have recently used their method to study the asymptotics of the
Taylor coefficients of $P_n^{-\beta}$
for the Szeg\H{o} polynomial \reff{eq1.1} with $n=3$,
but thus far only for $\beta > 1/2$.\footnote{
   The formula in their Theorem~4.4 has a misprint:
   the power $-1/2$ should be $\beta - \smfrac{3}{2}$.
}
It would be interesting to know whether this analysis can be extended
to the case $\beta < 1/2$, thereby providing an explicit proof
that some of the Taylor coefficients are asymptotically negative.
More generally, one might try to study
the elementary symmetric polynomials $E_{r,n}$:
after the $n=3$ Szeg\H{o} case $E_{2,3}$,
the next simplest would probably be the Lewy--Askey case $E_{2,4}$
[i.e., \reff{eq.lewy}].


\subsection{Some further remarks}

{\bf The half-plane property.}
Let us recall that a polynomial $P$ with complex coefficients
is said to have the {\em half-plane property}\/
\cite{Choe_hpp,Choe_05,Wagner_05,Branden_07,Borcea_09,Wagner_09,Branden_10,%
Wagner_11}
if either $P \equiv 0$ or else $P(x_1,\ldots,x_n) \neq 0$
whenever $x_1,\ldots,x_n$ are complex numbers
with strictly positive real part.\footnote{
   A polynomial $P \not\equiv 0$ with the half-plane property
   is also termed {\em Hurwitz stable}\/.
}
We shall show (Corollary~\ref{cor.multidim.HBW} below)
that if $P$ is a polynomial with real coefficients
that is strictly positive on $(0,\infty)^n$
and such that $P^{-\beta}$ is completely monotone on $(0,\infty)^n$
for at least one $\beta > 0$,
then $P$ necessarily has the half-plane property
(but not conversely).
The complete monotonicity of $P^{-\beta}$
can therefore be thought of as a strong quantitative form
of the half-plane property.
In particular, it follows that the determinantal polynomials considered
in Theorem~\ref{thm1.det} have the half-plane property ---
a fact that can easily be proven directly
(Corollary~\ref{cor.det.hpp} below).
The same is true for the quadratic polynomials considered
in Theorem~\ref{thm1.quadratic.cones}:
see \cite[Theorem~5.3]{Choe_hpp}
and Theorem~\ref{thm.quadratic.hpp} below.

\bigskip

{\bf The Rayleigh property.}
Complete monotonicity is also connected with the
Rayleigh property \cite{Choe_Rayleigh} for matroids and,
more generally, for multiaffine polynomials.
Let us say that a function $f$ is
{\em completely monotone of order $K$}\/
if the inequalities \reff{eq.compmono} hold for $0 \le k \le K$.
Thus, a function is completely monotone of order 0 (resp.\ 1)
if and only if it is nonnegative (resp.\ nonnegative and decreasing).
A function is completely monotone of order 2 if, in addition,
$\partial^2 f / \partial x_i \partial x_j \ge 0$ for all $i,j$.
Specializing this to $f = P^{-\beta}$ where $P$ is a polynomial,
we obtain
\be
   P \, {\partial^2 P \over \partial x_i \partial x_j}
   \;\le\;
   (\beta+1) \, {\partial P \over \partial x_i} \,
                {\partial P \over \partial x_j}
   \qquad\hbox{for all } i,j \;.
 \label{eq.C-Rayleigh}
\ee
If $P$ is multiaffine, then $\partial^2 P / \partial x_i^2 = 0$,
so it suffices to consider the cases $i \neq j$.
The inequality \reff{eq.C-Rayleigh} is then a generalization
of the Rayleigh (or negative-correlation) inequality
in which an extra constant $C=\beta+1$ is inserted on the right-hand side.
(The ordinary Rayleigh property corresponds to $\beta \downarrow 0$,
 hence to taking $f = -\log P$ and omitting the $k=0$ condition.)
It would be interesting to know whether the combinatorial consequences
of the Rayleigh property ---
such as the matroidal-support property \cite{Wagner_08} ---
extend to the $C$-Rayleigh property for arbitrary $C < \infty$.
It would also be interesting to extend the results of the present paper,
which address complete monotonicity of order $\infty$,
to complete monotonicity of finite orders $K$.
In what way do the conditions on $\beta$ become $K$-dependent?

\bigskip

{\bf Connected-spanning-subgraph polynomials.}
Let us remark that the literature contains some other examples
of multivariate polynomials $P$
for which $P^{-\beta}$ has all nonnegative Taylor coefficients,
for some specified set of numbers $\beta$.
For instance, Askey and Gasper \cite{Askey_77} showed
that this is the case for
\be
   P(x,y,z)  \;=\;  1 - \half (x+y+z) +  \half xyz
 \label{eq.askey-gasper}
\ee
whenever $\beta \ge (\sqrt{17}-3)/2 \approx 0.561553$;
Gillis, Reznick and Zeilberger \cite{Gillis_83}
later gave an elementary proof.
Likewise, Koornwinder \cite{Koornwinder_78} proved this for
\be
   P(x,y,z,u)  \;=\;  1 - \half (x+y+z+u) + \half (xyz + xyu + xzu + yzu)
                        - xyzu
 \label{eq.koornwinder}
\ee
whenever $\beta \ge 1$;
an elementary proof later emerged from the combined work of
Ismail and Tamhankar \cite{Ismail_79}
and Gillis--Reznick--Zeilberger \cite{Gillis_83}.

It turns out that these two examples also have a combinatorial interpretation:
not in terms of the spanning-tree polynomial $T_G({\bf x})$,
but rather in terms of the connected-spanning-subgraph polynomial
\cite{Royle-Sokal,Sokal_bcc2005}
\be
   C_G({\bf v})  \;=\;
   \sum_{\begin{scarray}
           A \subseteq E \\
           (V,A) \, {\rm connected}
         \end{scarray}}
   \prod_{e \in A} v_e
   \;,
\ee
which has $T_G({\bf x})$ as a limiting case:
\be
   T_G({\bf x})  \;=\;
  \lim_{\lambda \to 0} \lambda^{-(|V|-1)} C_G(\lambda {\bf x})
  \;.
 \label{eq.treelimit}
\ee
If we specialize to $G=C_n$ and make the change of variables
$v_i = -\lambda (1-z_i)$ with $0 < \lambda < n$
--- thus defining
 $P_{G,\lambda}(\bz) = T_G(-\lambda(\one-\bz)) / T_G(-\lambda\one)$ ---
it then turns out that the Askey--Gasper polynomial \reff{eq.askey-gasper}
corresponds to the case $n=3$, $\lambda=1$,
while the Koornwinder polynomial \reff{eq.koornwinder}
corresponds to the case $n=4$, $\lambda=2$.
On the other hand, in the limit $\lambda \to 0$
we recover a multiple of the Szeg\H{o} polynomial \reff{eq1.1};
this is simply a special case of \reff{eq.treelimit}.

In the same way that the complete monotonicity of $T_G^{-\beta}$
is a strong quantitative form of the half-plane property,
it turns out that the nonnegativity of Taylor coefficients of
$P_{G,\lambda}^{-\beta}$ in these examples 
is a strong quantitative form of the multivariate Brown--Colbourn property
(or more precisely, the multivariate property BC${}_\lambda$)
discussed in \cite{Royle-Sokal,Sokal_bcc2005}.
But it seems to be a difficult problem to determine the set of pairs
$(\lambda,\beta)$ for which $P_{G,\lambda}^{-\beta}$ has nonnegative
Taylor coefficients, even in the simplest case $G=C_3$.
We have some partial results on this problem,
but we leave these for a future paper.


\subsection{Plan of this paper}

The plan of this paper is as follows:
In Section~\ref{sec.compmono} we define complete monotonicity on cones
and recall the Bernstein--Hausdorff--Widder--Choquet theorem;
we also prove a general result showing that
complete monotonicity of $P^{-\beta}$ on a cone $C \subset V$
implies the nonvanishing of $P$ in the complex tube $C+iV$.
In Section~\ref{sec.constructions}
we discuss some general constructions
by which new polynomials $P$ with $P^{-\beta}$ completely monotone
can be obtained from old ones.
In Section~\ref{sec.det} we present the determinantal construction
and prove Theorems~\ref{thm1.det}, \ref{thm1.det.cones}
and \ref{thm1.det.cones.Jordan}.
In Section~\ref{sec.quadratic} we present the quadratic-form construction
and prove Theorem~\ref{thm1.quadratic.cones}.
In Section~\ref{sec.posdef} we present briefly the theory of
positive-definite functions (in the semigroup sense) on convex cones
--- which is a close relative of the theory of completely monotone
functions --- and its application to the class of cones treated here.
Finally, in Section~\ref{sec.graphs+matroids}
we apply the results of Sections~\ref{sec.compmono}--\ref{sec.quadratic}
to the spanning-tree polynomials of graphs
and the basis generating polynomials of matroids;
in particular we analyze the series-parallel case
and prove Theorem~\ref{thm1.serpar}.

In this arXiv version of the paper we include two appendices
that will be omitted from the journal version due to space constraints:
Appendix~\ref{sec.moore} reviewing the definition and main properties
of the Moore determinant for hermitian quaternionic matrices,
and Appendix~\ref{sec.gindikin} explaining an elementary proof
of Gindikin's characterization of parameters for which the
Riesz distribution is a positive measure.

We have tried hard to make this paper comprehensible
to the union (not the intersection!)\
of combinatorialists and analysts.
We apologize in advance to experts in each of these fields
for boring them every now and then with overly detailed
explanations of elementary facts.

\section{Complete monotonicity on cones} \label{sec.compmono}

In the Introduction we defined complete monotonicity for
functions on $(0,\infty)^n$.
For our later needs (see Sections~\ref{sec.det} and \ref{sec.quadratic}),
it turns out to be natural to consider complete monotonicity
on more general open convex cones $C \subset \R^n$.
This is a genuine generalization,
because for $n \ge 3$, an open convex cone is not necessarily
the affine image of a (possibly higher-dimensional) orthant,
i.e.\ it need not have ``flat sides'':
an example is the Lorentz cone
$\{{\bf x} \in \R^n \colon\, x_1 > \sqrt{x_2^2 + \ldots + x_n^2} \}$
in dimension $n \ge 3$.

\begin{definition}
   \label{def.compmono.cones}
Let $V$ be a finite-dimensional real vector space,
and let $C$ be an open convex cone in $V$.
Then a $C^\infty$ function $f \colon\, C \to \R$
is termed {\em completely monotone}\/
if for all $k \ge 0$,
all choices of vectors ${\bf u}_1, \ldots, {\bf u}_k \in C$,
and all ${\bf x} \in C$,
we have
\be
   (-1)^k D_{{\bf u}_1} \cdots D_{{\bf u}_k} f(x)  \;\ge\;  0
 \label{eq.def.compmono.cones}
\ee
where $D_{\bf u}$ denotes a directional derivative.
A function $f$
is termed {\em conditionally completely monotone}\/
if the inequality \reff{eq.def.compmono.cones} holds
for all $k \ge 1$ but not necessarily for $k=0$.\footnote{
   The terminology ``conditionally completely monotone'' is new,
   but we think it felicitous:
   it is chosen by analogy with ``conditionally positive definite matrix''
   \cite{Berg_84,Bapat_97},
   with which this concept is in fact closely related
   \cite{Horn_67,Horn_69a,Horn_69b,Berg_08,Berg_84}.
   ({\em Warning}\/:  The book \cite{Berg_84} uses the term
    ``negative definite'' for what we would call
    ``conditionally negative definite''.)

   Please note that if $f$ is bounded below,
   then $f$ is conditionally completely monotone
   if and only if there exists a constant $c$ such that $f+c$ is
   completely monotone.
   But there also exist conditionally completely monotone functions
   that are unbounded below
   (and hence for which such a constant $c$ cannot exist):
   examples on $(0,\infty)$ are
   $f(x) = -(a+x)^\alpha$ and $f(x) = -\log(a+x)$
   with $a \ge 0$ and $0 < \alpha \le 1$.

   Of course, it follows immediately from the definition
   that $f$ is conditionally completely monotone
   if and only if $-D_{\bf u} f$ is completely monotone
   for all vectors ${\bf u} \in C$.
   In the multidimensional case this seems rather difficult to work with;
   but in the one-dimensional case $C=(0,\infty)$
   it says that $f$ is conditionally completely monotone
   if and only if $-f'$ is completely monotone.
}
\end{definition}

\noindent
Of course, if the inequality \reff{eq.def.compmono.cones}
holds for all ${\bf u}_1, \ldots, {\bf u}_k$ in some set $S$,
then by linearity and continuity it holds also for all
${\bf u}_1, \ldots, {\bf u}_k$ in the closed convex cone generated by $S$.
This observation also shows the equivalence of
Definition~\ref{def.compmono.cones},
specialized to the case $V = \R^n$ and $C = (0,\infty)^n$,
with the definition given in the Introduction.

If $T\colon\, (V_1,C_1) \to (V_2,C_2)$ is a positive linear map
(i.e., a linear map $T \colon\, V_1 \to V_2$ satisfying $T[C_1] \subseteq C_2$)
and $f\colon\, C_2 \to \R$ is completely monotone,
then it is easily seen that
$f \circ T \colon\, C_1 \to \R$ is completely monotone.
Conversely, if $f \circ T$ is completely monotone
for {\em all}\/ positive linear maps
$T \colon\, (\R^n, (0,\infty)^n) \to (V_2,C_2)$
{\em for arbitrarily large $n$}\/,
then $f$ is completely monotone:
for if \reff{eq.def.compmono.cones} fails for some $k$,
then we can take $n=k$ and $T {\bf e}_i = {\bf u}_i$
(where ${\bf e}_i$  is the $i$th coordinate unit vector in $\R^n$)
and $f \circ T$ will fail one of the $k$th-order
complete-monotonicity inequalities.

Let us next recall some elementary facts.
If $f$ is completely monotone,
then $f(0^+) = \lim_{x \to 0, x \in C} f(x)$ exists
and equals $\sup_{x \in C} f(x)$,
but it might be $+\infty$.
The product of two completely monotone functions is completely monotone.
If $f$ is completely monotone and
$\Phi\colon\, [0,\infty) \to [0,\infty)$ is absolutely monotone
(i.e.\ its derivatives of all orders are everywhere nonnegative),
then $\Phi \circ f$ is completely monotone.
If $f$ is conditionally completely monotone and
$\Phi\colon\, (-\infty,\infty) \to [0,\infty)$ is absolutely monotone,
then $\Phi \circ f$ is completely monotone.
(In particular, this occurs when $\Phi$ is the exponential function.)
Finally, a locally uniform limit of a sequence of
completely monotone functions is completely monotone.

The fundamental fact in the theory of completely monotone functions
on $(0,\infty)$ is the Bernstein--Hausdorff--Widder theorem \cite{Widder_46}:
A function $f$ defined on $(0,\infty)$ is completely monotone
if and only if it can be written in the form
\be
   f(x)  \;=\;  \int_0^\infty e^{-tx} \, d\mu(t)
 \label{eq.BHW}
\ee
where $\mu$ is a nonnegative Borel measure on $[0,\infty)$.
We shall need a multidimensional version of the
Bernstein--Hausdorff--Widder theorem, valid for arbitrary cones.
Such a result was proven by Choquet \cite{Choquet_69}\footnote{
   See also Nussbaum \cite{Nussbaum_55},
   Devinatz and Nussbaum \cite{Devinatz_61},
   Hirsch \cite[Section~VII.2]{Hirsch_72},
   Gl\"ockner \cite[Section~16]{Glockner_03}
   and Thomas \cite{Thomas_04}.
}:


\begin{theorem}[Bernstein--Hausdorff--Widder--Choquet theorem]
   \label{thm.multidim.HBW}
Let $V$ be a finite-dimensional real vector space,
let $C$ be an open convex cone in $V$,
and let $C^* = \{ \ell \in V^* \colon\: \langle \ell,x \rangle \ge 0
                       \hbox{ for all } x \in C \}$
be the closed dual cone.
Then a function $f \colon\, C \to \R$
is completely monotone if and only if
there exists a positive measure $\mu$ on $C^*$ satisfying
\be
   f(x)  \;=\;  \int\limits_{C^*} e^{-\langle \ell,x \rangle} \, d\mu(\ell)
   \;.
 \label{eq.thm.multidim.HBW}
\ee
In this case, $\mu(C^*) = f(0^+)$;
in particular, $\mu$ is finite if and only if $f$ is bounded.

In particular, if $f$ is completely monotone on $C$,
then it is extendible [using \reff{eq.thm.multidim.HBW}]
to an analytic function on the complex tube $C+iV$ satisfying
\be
   |D_{{\bf u}_1} \cdots D_{{\bf u}_k} f(x+iy)|
   \;\le\;
   (-1)^k D_{{\bf u}_1} \cdots D_{{\bf u}_k} f(x)
 \label{eq.thm.multidim.HBW.2}
\ee
for all $k \ge 0$, $x \in C$, $y \in V$
and ${\bf u}_1, \ldots, {\bf u}_k \in C$.
\end{theorem}

\medskip\noindent
{\bf Remarks.}
1.  Since $C$ is nonempty and open, it is not hard to see that
$\ell \in C^* \setminus \{0\}$ implies $\langle \ell,x \rangle > 0$
for all $x \in C$.  It then follows from \reff{eq.thm.multidim.HBW}
that either
\begin{itemize}
   \item[(a)] $\mu$ is supported on $\{0\}$, in which case $f$ is constant,
\end{itemize}
or else
\begin{itemize}
   \item[(b)] we have the {\em strict}\/ inequality
\be
   (-1)^k D_{{\bf u}_1} \cdots D_{{\bf u}_k} f(x)  \;>\;  0
 \label{eq.def.compmono.cones.strict}
\ee
for all $k \ge 0$, all ${\bf u}_1, \ldots, {\bf u}_k \in C$,
and all ${\bf x} \in C$.
\end{itemize}
Furthermore, for ${\bf u}_1, \ldots, {\bf u}_k$ in the {\em closure}\/ of $C$,
the left-hand side of \reff{eq.def.compmono.cones.strict}
is either strictly positive for all $x \in C$
or else identically zero on $C$.
Of course, the latter case can occur:
e.g.\ $f(x_1,x_2) = e^{-x_1}$ on $(0,\infty)^2$.

2.  In our definition of complete monotonicity,
the cone $C$ plays two distinct roles:
it is the domain on which $f$ is defined,
and it provides the direction vectors ${\bf u}_i$
for which the inequalities \reff{eq.def.compmono.cones} hold.
Choquet \cite{Choquet_69} elegantly separates these roles,
and considers functions on an arbitrary open set $\Omega \subseteq V$
that are completely monotone with respect to the cone $C$.
He then proves the integral representation \reff{eq.thm.multidim.HBW}
under the hypothesis $\Omega + C \subseteq \Omega$.
This is a beautiful generalization, but we shall not need it.
\medskip

\bigskip

By virtue of Theorem~\ref{thm.multidim.HBW},
one way to test a function $f$ for complete monotonicity
is to compute its inverse Laplace transform
and ask whether it is nonnegative and supported on $C^*$.
Of course, this procedure is not necessarily well-defined,
because the inverse Laplace transform need not exist;
moreover, if it does exist, it may need to be understood
as a distribution in the sense of Schwartz \cite{Schwartz_66}
rather than as a pointwise-defined function.
But we can say this:  If $f \colon\, C \to \R$
is the Laplace transform of a distribution $T$ on $V^*$,
then $f$ is completely monotone if and only if
$T$ is positive (hence a positive measure) and supported on $C^*$.
This follows from the injectivity of the Laplace transform
on the space $\scrd'(\R^n)$ of distributions \cite[p.~306]{Schwartz_66}.
Note that the complete monotonicity of $f$ can fail
either because $T$ fails to be positive
or because $T$ fails to be supported on $C^*$.
In the former case, we can conclude that
$f$ is not completely monotone on
{\em any}\/ nonempty open convex subcone $C' \subseteq C$.
In the latter case, $T$ might possibly be supported on some larger proper cone;
but if it isn't (e.g.\ if the smallest convex cone containing the support
of $T$ is all of $V^*$), then once again we can conclude that
$f$ is not completely monotone on
any nonempty open convex subcone $C' \subseteq C$.
And finally, if $f$ is {\em not}\/ the Laplace transform of any
distribution on $V^*$, then it is certainly not the Laplace transform of a
positive measure, and hence is not completely monotone on
any nonempty open convex subcone $C' \subseteq C$.

\bigskip

In the applications to be made in this paper,
the functions $f$ will typically be of the form $F^{-\beta}$,
where $F$ is a function (usually a polynomial)
that is strictly positive on the cone $C$
and has an analytic continuation to the tube $C+iV$
(for polynomials this latter condition of course holds trivially).
The following corollary of Theorem~\ref{thm.multidim.HBW}
shows that the complete monotonicity of $F^{-\beta}$ on the real cone $C$
implies the absence of zeros of $F$ in the complex tube $C+iV$:

\begin{corollary}
   \label{cor.multidim.HBW}
Let $V$ be a finite-dimensional real vector space and
let $C$ be an open convex cone in $V$.
Let $F$ be an analytic function on the tube $C+iV$
that is real and strictly positive on $C$.
If $F^{-\beta}$ is completely monotone on $C$ for at least one $\beta > 0$,
then $F$ is nonvanishing on $C+iV$.
[In particular, when $V = \R^n$ and $C = (0,\infty)^n$,
 the function $F$ has the half-plane property.]
\end{corollary}

\proof
Suppose that $G = F^{-\beta}$ is completely monotone on $C$;
then by Theorem~\ref{thm.multidim.HBW}
it has an analytic continuation to $C+iV$ (call it also $G$).
Now suppose that $S = \{z \in C+iV \colon\, F(z)=0 \}$ is nonempty.
Choose a simply connected domain $D \subset (C+iV) \setminus S$
such that $D \cap C \neq \emptyset$ and
$\bar{D} \cap S \neq \emptyset$.\footnote{
   For instance, let $\Omega$ be a simply connected open subset of $C$
   whose closure is a compact subset of $C$
   and which satisfies $(\Omega +iV) \cap S \neq \emptyset$;
   fix a norm $\| \,\cdot\, \|$ on $V$;
   and let $$R  \;=\; \inf\limits_{\begin{scarray}
                                      x \in \Omega,\, y \in V \\
                                      x+iy \in S
                                   \end{scarray}}  \|y\|  \;.$$
   By compactness we must have $R > 0$
   (for otherwise we would have $\bar{\Omega} \cap S \neq \emptyset$,
    contrary to the hypothesis that $F>0$ on $C$).
   Now take $D = \Omega + iB_R$,
   where $B_R = \{y \in V \colon\, \|y\| < R \}$.
}
Then $H = F^{-\beta}$ is a well-defined analytic function on $D$
(we take the branch that is real and positive on $D \cap C$).
On the other hand, $H$ coincides with $G$ on the real environment $D \cap C$,
so it must coincide with $G$ everywhere in $D$.
But $\lim\limits_{z \to z_0} |H(z)| = +\infty$
for all $z_0 \in \bar{D} \cap S$,
which contradicts the analyticity of $G$ on $C+iV$.
\qed

\medskip\noindent
{\bf Remarks.}
1.  It also follows that the analytic function $G = F^{-\beta}$
defined on $C+iV$ is nonvanishing there.

2.  The hypothesis that $F$ have an analytic continuation to $C+iV$
is essential;  it cannot be derived as a consequence of the
complete monotonicity of $F^{-\beta}$ on $C$.
To see this, take $V=\R$ and $C=(0,\infty)$
and consider $F(x) = (1 + \half e^{-x})^{-1/\beta}$ with any $\beta > 0$.

3.  The converse of Corollary~\ref{cor.multidim.HBW} is
easily seen to be false:
for instance, the univariate polynomial $P(x) = 1+x^2$ has the
half-plane property (i.e.\ is nonvanishing for $\real x > 0$),
but $P^{-\beta}$ is not completely monotone
on $(0,\infty)$ for any $\beta > 0$.
The same holds for the bivariate multiaffine polynomial
$P(x_1,x_2) = 1 + x_1 x_2$.
So the complete monotonicity of $P^{-\beta}$ for some $\beta > 0$
is strictly stronger than the half-plane property.
It would be interesting to know whether similar counterexamples can be found
if $P$ is required to be homogeneous, or homogeneous and multiaffine.
\medskip

In this paper we will typically consider a polynomial $P$
that is strictly positive on an open convex cone $C$,
and we will ask for which values of $\beta$
the function $P^{-\beta}$ is completely monotone.
We begin with a trivial observation:
If $P$ is a nonconstant polynomial,
then $P^{-\beta}$ cannot be completely monotone
on any nonempty open convex cone for any $\beta < 0$
(because $P$ grows at infinity in all directions
except at most a variety of codimension 1);
and $P^{-\beta}$ is trivially completely monotone for $\beta=0$.
So we can restrict attention to $\beta > 0$.

%

Given a function $F \colon\, C \to (0,\infty)$
--- for instance, a polynomial ---
we can ask about the set
\be
   \scrb_F  \;=\;  \{ \beta > 0 \colon\, F^{-\beta}
                    \hbox{ is completely monotone on } C  \}  \;.
\ee
Clearly $\scrb_F$ is a closed additive subset of $(0,\infty)$.
In particular, we either have $\scrb_F \subseteq [\epsilon,\infty)$
for some $\epsilon > 0$ or else $\scrb_F = (0,\infty)$.
The following easy lemma \cite{Horn_67,Berg_08} characterizes the latter case:

\begin{lemma}
   \label{lemma.BF.0infty}
Let $V$ be a finite-dimensional real vector space,
let $C$ be an open convex cone in $V$,
and let $F \colon\, C \to (0,\infty)$.
Then the following are equivalent:
\begin{itemize}
   \item[(a)]  $F^{-\beta}$ is completely monotone on $C$ for all $\beta > 0$.
   \item[(b)]  $F^{-\beta_i}$ is completely monotone on $C$ for a sequence
      $\{\beta_i\}$ of strictly positive numbers converging to 0.
   \item[(c)]  $-\log F$ is conditionally completely monotone on $C$.
\end{itemize}
\end{lemma}

\proof
(a)$\implies$(b) is trivial, and (b)$\implies$(c) follows from
\be
   -\log F  \;=\;  \lim\limits_{\beta \downarrow 0}
                           {F^{-\beta} - 1  \over \beta}
\ee
and its derivatives with respect to $x$.
Finally, (c)$\implies$(a) follows from
$F^{-\beta} = \exp(-\beta \log F)$
and the fact that $\exp$ is absolutely monotone
(i.e.\ has all derivatives nonnegative) on $(-\infty,\infty)$.
\qed

Already for $C = (0,\infty)$ it seems to be a difficult problem
to characterize in a useful way the functions $F$ described in
Lemma~\ref{lemma.BF.0infty},
or even the subclass consisting of polynomials $P$.\footnote{
   Functions $f = F^{-1}$ for which $f^\beta$
   is completely monotone for all $\beta > 0$
   are sometimes called {\em logarithmically completely monotone}\/
   \cite{Berg_08}.
}
For polynomials $P(x) = \prod (1+x/x_i)$,
a {\em necessary}\/ condition from Corollary~\ref{cor.multidim.HBW}
is that $P$ have the half-plane property,
i.e.\ $\real x_i \ge 0$ for all $i$.
A {\em sufficient}\/ condition is that all $x_i$ be real and positive;
and for quadratic polynomials this condition is necessary as well.
But already for quartic polynomials the situation becomes more complicated:
for instance, we can take $x_1 = a + b i$, $x_2 = a - b i$,
$x_3 = x_4 = c$ with $0 < c \le a$ and $b \in \R$,
and it is not hard to see that
$-\log P$ is conditionally completely monotone on $(0,\infty)$.\footnote{
   The function $-\log P$ is conditionally completely monotone on $(0,\infty)$
   if and only if $(\log P)' = \sum\limits_i (x+x_i)^{-1}$
   is completely monotone on $(0,\infty)$;
   and this happens if and only if its inverse Laplace transform,
   which is $g(t) = \sum\limits_i e^{-tx_i}$,
   is nonnegative on $[0,\infty)$.
}

It also seems to be a difficult problem to characterize the
closed additive subsets $S \subseteq (0,\infty)$
that can arise as $S = \scrb_F$.

\begin{example}
\rm
Fix $a>0$, and consider $F(x) = (1 + a e^{-x})^{-1}$.
Then the function $F(x)^{-\beta} = (1 + a e^{-x})^{\beta}$
is obviously completely monotone on $(0,\infty)$
whenever $\beta \in \{0,1,2,3,\ldots\}$.
On the other hand, if $\beta \notin \{0,1,2,3,\ldots\}$
we claim that $F^{-\beta}$ is {\em not}\/ completely monotone.
Indeed, for $0 < a \le 1$ the convergent binomial expansion
\be
   F(x)^{-\beta} \;=\; \sum_{k=0}^\infty   a^k {\beta \choose k} e^{-kx}
\ee
shows that $F^{-\beta}$ is the Laplace transform of the signed measure
$\sum_{k=0}^\infty  a^k {\beta \choose k} \delta_k$,
which is nonnegative if and only if $\beta \in \{0,1,2,3,\ldots\}$.
On the other hand, for $a > 1$ the function $F^{-\beta}$
has singularities in the right half-plane at $x = \log a \pm i\pi$
whenever $\beta \notin \{0,1,2,3,\ldots\}$,
so it is not the Laplace transform of {\em any}\/ distribution.
%
\medskip
\end{example}

%
%
%
%
%

\begin{example}
\rm
It is an interesting problem
\cite{Askey_73,Askey_74b,Fields_75,Gasper_75,Moak_87,Zastavnyi_00}
to determine the pairs $(\mu,\lambda) \in \R^2$ for which the function
\be
   F_{\mu,\lambda}(x)  \;=\; x^{-\mu} (x^2+1)^{-\lambda}
\ee
is completely monotone on $(0,\infty)$.
It is easy to show that there is a function $\mu_\star(\lambda)$
such that $F_{\mu,\lambda}$ is completely monotone
if and only if $\mu \ge \mu_\star(\lambda)$;
furthermore, the function $\mu_\star$ is subadditive.
The state of current knowledge about $\mu_\star$ seems to be:
\begin{subeqnarray}
   \mu_\star(\lambda) \;=\; -2\lambda
           & \quad & \hbox{for } \lambda \le 0 \\
   \lambda \;<\; \mu_\star(\lambda) \;\le\; \min(2\lambda,1)
           & \quad & \hbox{for } 0 < \lambda < 1  \\
   \mu_\star(\lambda) \;=\; 2\lambda + o(\lambda)
           & \quad & \hbox{for } \lambda \downarrow 0  \\
   \mu_\star(\lambda) \;=\; \lambda
           & \quad & \hbox{for } \lambda \ge 1
\end{subeqnarray}
It seems to be an open problem even to prove that $\mu_\star$ is continuous.
\medskip
\end{example}

\section{Constructions}
    \label{sec.constructions}

In this section we discuss some general constructions
by which new polynomials $P$ with $P^{-\beta}$ completely monotone
can be obtained from old ones.
In the situations we have in mind,
the vector space $V$ decomposes as a direct sum $V = V_1 \oplus V_2$
and the cone $C$ is a product cone $C = C_1 \times C_2$
(with $C_1 \subset V_1$ and $C_2 \subset V_2$).
Since we shall be using the letters $A,B,C,D$ in this section
to denote functions, we shall write our cones as $\scrc$.

Let us begin with a trivial fact:
a function $f(x,y)$ that is completely monotone on $\scrc_1 \times \scrc_2$
can be specialized by fixing $y$ to a specific point in $\scrc_2$,
and the resulting function will be completely monotone on $\scrc_1$.
In particular, this fixed value can then be taken to zero or infinity,
and if the limit exists --- possibly with some rescaling ---
then the limiting function is also completely monotone on $\scrc_1$.
Rather than stating a general theorem of this kind,
let us just give the special case that we will need,
which concerns functions of the form
\be
   f(x,y)  \;=\; [A(x) \,+\, B(x) \, y]^{-\beta}
 \label{eq.A+Bx}
\ee
with $V_2 = \R$ and $\scrc_2 = (0,\infty)$.

\begin{lemma}
  \label{lemma.compmono.delcon}
Let $V$ be a finite-dimensional real vector space and
let $\scrc$ be an open convex cone in $V$.
Fix $\beta > 0$, and let $A,B \colon\; \scrc \to (0,\infty)$.
If $(A+By)^{-\beta}$ is completely monotone on $\scrc \times (0,\infty)$,
then $A^{-\beta}$ and $B^{-\beta}$ are completely monotone on $\scrc$.
\end{lemma}

\proof
Restrict to fixed $y \in (0,\infty)$ and then take $y \downarrow 0$;
this proves that $A^{-\beta}$ is completely monotone.
Restrict to fixed $y \in (0,\infty)$, multiply by $y^{\beta}$
and then take $y \uparrow \infty$;
this proves that $B^{-\beta}$ is completely monotone.
\qed

\noindent
As we shall see later, this trivial lemma is an analytic version of
deletion ($y \to 0$) or contraction ($y \to \infty$)
for graphs or matroids.

Let us also observe a simple but important fact about complete monotonicity
for functions defined on a product cone $\scrc_1 \times \scrc_2$:

\begin{lemma}
  \label{lemma.compmono.trivial}
For $i=1,2$, let $V_i$ be a finite-dimensional real vector space and
let $\scrc_i$ be an open convex cone in $V_i$.
Let $f \colon\, \scrc_1 \times \scrc_2 \to \R$.
Then the following are equivalent:
\begin{itemize}
   \item[(a)] $f$ is completely monotone on $\scrc_1 \times \scrc_2$.
   \item[(b)] For all $k \ge 0$, all $y \in \scrc_2$,
       and all choices of vectors ${\bf u}_1, \ldots, {\bf u}_k \in \scrc_2$,
       the function
\be
   F_{k,y,{\bf u}_1, \ldots, {\bf u}_k}(x)
   \;\equiv\;
   \left.
   (-1)^k \, {\partial^k \over \partial t_1 \dots \partial t_k}
   \, f(x, y + t_1 {\bf u}_1 + \ldots + t_k {\bf u}_k)
   \right| _{t_1 = \,\ldots\, = t_k = 0}
\ee
is completely monotone on $\scrc_1$.
   \item[(c)] For all $k \ge 0$, all $y \in \scrc_2$,
       and all choices of vectors ${\bf u}_1, \ldots, {\bf u}_k \in \scrc_2$,
       there exists a positive measure $\mu_{k,y,{\bf u}_1, \ldots, {\bf u}_k}$
       on $\scrc_1^*$ such that
\be
   F_{k,y,{\bf u}_1, \ldots, {\bf u}_k}(x)
   \;=\;
   \int\limits_{\scrc_1^*}  e^{-\langle \ell,x \rangle} \,
        d\mu_{k,y,{\bf u}_1, \ldots, {\bf u}_k}(\ell)  \;.
\ee
\end{itemize}
%
\end{lemma}

\noindent
In particular, when $V_2 = \R$ and $\scrc_2 = (0,\infty)$,
(b) reduces to the statement that the functions
$F_{k,y}(x) = (-1)^k \partial^k f/\partial y^k$
are completely monotone on $\scrc_1$ for all $k \ge 0$ and all $y > 0$,
and (c) reduces analogously.

\proofof{Lemma~\ref{lemma.compmono.trivial}}
(a) $\iff$ (b) is a trivial consequence
of the definition \reff{eq.def.compmono.cones},
while (b) $\iff$ (c) follows immediately from the
the Bernstein--Hausdorff--Widder--Choquet theorem
(Theorem~\ref{thm.multidim.HBW}).
\qed

Statement (c) can be rephrased loosely as saying that
the inverse Laplace transform of $f(x,y)$ with respect to $x$
is a completely monotone function of $y \in \scrc_2$.
(To make this more precise, one should add the same qualifications
 as in the paragraph after Theorem~\ref{thm.multidim.HBW}.)

One important application of this lemma
concerns functions of the form \reff{eq.A+Bx}:

\begin{lemma}
  \label{lemma.compmono.invLT.affine}
Let $V$ be a finite-dimensional real vector space and
let $\scrc$ be an open convex cone in $V$.
Fix $\beta > 0$, and let $A,B \colon\; \scrc \to (0,\infty)$.
Then the following are equivalent:
\begin{itemize}
   \item[(a)]  $(A+By)^{-\beta}$ is completely monotone on
       $\scrc \times (0,\infty)$.
   \item[(b)]  $B^{-\beta} \exp(-tA/B)$ is completely monotone
      on $\scrc$ for all $t \ge 0$.
   \item[(c)]  $B^{\kappa} (A+zB)^{-(\beta+\kappa)}$ is completely monotone
      on $\scrc$ for all $\kappa \ge -\beta$ and all $z \ge 0$.
   \item[(d)]  $B^{k} (A+zB)^{-(\beta+k)}$ is completely monotone
      on $\scrc$ for all integers $k \ge 0$ and all $z \ge 0$.
\end{itemize}
\end{lemma}

\proof
We have the Laplace-transform formula
\be
   (A+By)^{-\beta}
   \;=\;
   \int\limits_0^\infty e^{-ty} \;
   {t^{\beta-1} \over \Gamma(\beta)} \, B^{-\beta} \, \exp(-tA/B)
   \; dt
   \;.
 \label{eq.lemma.compmono.invLT.affine.proof}
\ee
Therefore,
Lemma~\ref{lemma.compmono.trivial}(a) $\iff$ (c)
with $\scrc_1 = (0,\infty)$ and $\scrc_2 = \scrc$
proves the equivalence of (a) and (b).

Now assume that $B^{-\beta} \exp(-tA/B)$ is completely monotone
on $\scrc$ for all $t \ge 0$.
Then we can multiply by $e^{-zt} t^{p-1}/\Gamma(p)$
for any $p > 0$ and integrate over $t \in (0,\infty)$,
and the result will be completely monotone.
This (together with a trivial evaluation at $t=0$ to handle $p=0$)
shows that (b) $\implies$ (c).

(c) $\implies$ (d) is trivial.

The equivalence of (a) and (d) follows from
Lemma~\ref{lemma.compmono.trivial}(a) $\iff$ (b),
used with $\scrc_1 = \scrc$ and $\scrc_2 = (0,\infty)$.
\qed


\bigskip

\begin{corollary}
  \label{cor.compmono.invLT.affine}
Let $V$ be a finite-dimensional real vector space,
let $\scrc$ be an open convex cone in $V$,
and let $A,B \colon\; \scrc \to (0,\infty)$.
Define
\begin{eqnarray}
   \scrb_B  & = &  \{ \beta > 0 \colon\, B^{-\beta}
                    \hbox{ is completely monotone on } \scrc  \}
      \\[1mm]
   \scrb_{A+By}  & = &  \{ \beta > 0 \colon\, (A+By)^{-\beta}
                 \hbox{ is completely monotone on } \scrc \times (0,\infty)  \}
     \qquad
\end{eqnarray}
Then $\scrb_{A+By} + \scrb_B \subseteq \scrb_{A+By}$.

In particular, if $\scrb_B = (0,\infty)$,
then $\scrb_{A+By}$ is either the empty set or all of $(0,\infty)$
or a closed interval $[\beta_0,\infty)$ with $\beta_0 > 0$.
\end{corollary}

\proof
This follows immediately from
Lemma~\ref{lemma.compmono.invLT.affine}(a)$\iff$(b):
for if $\beta \in \scrb_{A+By}$ and $\lambda \in \scrb_B$,
then $B^{-\beta} \exp(-tA/B)$ and $B^{-\lambda}$
are both completely monotone on $\scrc$,
hence so is their product, hence $\beta+\lambda \in \scrb_{A+By}$.
\qed

Lemma~\ref{lemma.compmono.invLT.affine}
leads to the following extremely important result,
which (as we shall see later) is an analytic version of
parallel connection for graphs or matroids:

\begin{proposition}
   \label{prop.compmono.parallel}
Let $V$ be a finite-dimensional real vector space and
let $\scrc$ be an open convex cone in $V$.
Fix $\beta > 0$, and let $A,B,C,D \colon\; \scrc \to (0,\infty)$.
Suppose that $(A+By)^{-\beta}$ and $(C+Dy)^{-\beta}$
are completely monotone on $\scrc \times (0,\infty)$.
Then the same is true of $(AD+BC+BDy)^{-\beta}$.
\end{proposition}

\proof
By Lemma~\ref{lemma.compmono.invLT.affine},
$B^{-\beta} \exp(-tA/B)$ and $D^{-\beta} \exp(-tC/D)$
are completely monotone on $\scrc$ for all $t \ge 0$.
Hence the same is true of their product, which is \linebreak
$(BD)^{-\beta} \exp[-t(AD+BC)/(BD)]$.
But then using Lemma~\ref{lemma.compmono.invLT.affine} again
(this time in the reverse direction),
we conclude that $(AD+BC+BDy)^{-\beta}$
is completely monotone on $\scrc \times (0,\infty)$.
\qed

We also have an analytic version of series connection for graphs or matroids,
but only for $\beta \ge 1/2$:

\begin{proposition}
   \label{prop.compmono.series}
Let $V$ be a finite-dimensional real vector space and
let $\scrc$ be an open convex cone in $V$.
Fix $\beta \ge 1/2$, and let $A,B,C,D \colon\; \scrc \to (0,\infty)$.
Suppose that $(A+By)^{-\beta}$ and $(C+Dy)^{-\beta}$
are completely monotone on $\scrc \times (0,\infty)$.
Then the same is true of $[AC+ (AD+BC)y]^{-\beta}$.
\end{proposition}

To prove Proposition~\ref{prop.compmono.series},
we begin with a lemma that we think is of independent interest;
both the sufficiency and the necessity will play important roles for us.

\begin{lemma}
   \label{lemma.compmono.series}
For $\beta \in \R$ and $\lambda > 0$, the function
\be
   F_{\beta,\lambda}(u,v)  \;=\;
     (u+v)^{-\beta} \, \exp\biggl( - \lambda {uv \over u+v} \biggr)
 \label{lemma.compmono.series.eq1}
\ee
is completely monotone on $(0,\infty)^2$ if and only if $\beta \ge 1/2$.

In particular, for $\beta \ge 1/2$
there exists a positive measure $\mu_{\beta,\lambda}$
on $[0,\infty)^2$ such that
\be
   (u+v)^{-\beta} \, \exp\biggl( - \lambda {uv \over u+v} \biggr)
   \;=\;
   \int\limits_{[0,\infty)^2}
       e^{-t_1 u - t_2 v} \, d\mu_{\beta,\lambda}(t_1,t_2)
   \;.
 \label{lemma.compmono.series.eq2}
\ee
%
\end{lemma}

\proof
``If'':
Since $(u+v)^{-(\beta-\!\smhalf)}$ is completely monotone when $\beta \ge 1/2$,
it suffices to prove the complete monotonicity for $\beta = 1/2$.
But this follows immediately from the identity
\be
   (u+v)^{-1/2} \, \exp\biggl( - \lambda {uv \over u+v} \biggr)
   \;=\;
   {1 \over \sqrt{\pi}}
   \int_{-\infty}^\infty 
   \exp\!\left[
     - \Bigl( s + {\sqrt{\lambda} \over 2} \Bigr)^2 u
     - \Bigl( s - {\sqrt{\lambda} \over 2} \Bigr)^2 v
         \right]
   \, ds
   \;, \qquad
 \label{lemma.compmono.series.eq3}
\ee
which is easily verified by completing the square in the Gaussian integral.
The statement about the measure $\mu_{\beta,\lambda}$ then follows from
the Bernstein--Hausdorff--Widder--Choquet theorem
(Theorem~\ref{thm.multidim.HBW}).

``Only if'':
If $F_{\beta,\lambda}$ is completely monotone,
then so is $F_{\beta',\lambda}$ for all $\beta' > \beta$;
so it suffices to prove the failure of complete monotonicity
for $0 < \beta < 1/2$.  Now, by Lemma~\ref{lemma.compmono.trivial},
$F_{\beta,\lambda}$ is completely monotone on $(0,\infty)^2$
if and only if the functions
\be
   F_{\beta,\lambda;k,v}(u)  \;=\;
   (-1)^k {\partial^k \over \partial v^k} F_{\beta,\lambda}(u,v)
\ee
are completely monotone on $(0,\infty)$ for all $k \ge 0$ and all $v > 0$,
or equivalently if their inverse Laplace transforms with respect to $u$,
\be
   G_{\beta,\lambda;k,v}(t)  \;=\;
   (-1)^k {\partial^k \over \partial v^k}
   \left[ (t/\lambda)^{(\beta-1)/2} e^{-(t+\lambda)v} v^{1-\beta}
              I_{\beta-1}(2v \sqrt{\lambda t}) 
   \right]
\ee
(see \cite[p.~245, eq.~5.6(35)]{Erdelyi_54}),
are nonnegative for all $k \ge 0$ and all $t,v > 0$
(here $I_{\beta-1}$ is the modified Bessel function).
For $k=0$ this manifestly holds for all $\beta \ge 0$;
but let us now show that for $k=1$ it holds only for $\beta \ge 1/2$.
We have
\be
   G_{\beta,\lambda;1,v}(t)  \;=\;
   (t/\lambda)^{(\beta-1)/2} e^{-(t+\lambda)v} v^{1-\beta}
   \left[ \Bigl( t+\lambda + {\beta-1 \over v} \Bigr)
              I_{\beta-1}(2v \sqrt{\lambda t})
          \,-\,
          2 \sqrt{\lambda t} I'_{\beta-1}(2v \sqrt{\lambda t})
   \right]  \;,
\ee
and we need the term in square brackets to be nonnegative for all $t,v > 0$.
Write $x = 2v \sqrt{\lambda t}$ and eliminate $v$ in favor of $x$;
we need
\be
   t \,+\, 2 \sqrt{\lambda t} \left[ {\beta-1 \over x} -
                                     {I'_{\beta-1}(x) \over I_{\beta-1}(x)}
                              \right]
     \,+\, \lambda
   \;\ge\; 0
\ee
for all $t,x > 0$.
This quadratic in $\sqrt{t}$ is nonnegative for all $t > 0$
if and only if
\be
   {\beta-1 \over x} - {I'_{\beta-1}(x) \over I_{\beta-1}(x)}
   \;\ge\;  -1  \;.
\ee
But using the large-$x$ asymptotic expansion
\be
   {d \over dx} \, \log I_{\beta-1}(x)
   \;=\;
   1 \,-\, {1 \over 2x} \,+\, O(1/x^2)
   \;,
\ee
we see that
\be
  {\beta-1 \over x} - {I'_{\beta-1}(x) \over I_{\beta-1}(x)}
  \;=\;
  -1 \,+\, {\beta-\smhalf \over x} \,+\, O(1/x^2)
  \;,
\ee
which is $<-1$ for all sufficiently large $x$ whenever $\beta < 1/2$.
\qed

\medskip\noindent
{\bf Remarks.}
1.  It is obvious by rescaling of $u$ and $v$ that, for any given $\beta$,
the functions $F_{\beta,\lambda}$ are either completely monotone
for all $\lambda>0$ or for none.

2.  The appeal to the Bernstein--Hausdorff--Widder--Choquet theorem
can be avoided:
for $\beta=1/2$, the integral representation \reff{lemma.compmono.series.eq3}
already provides the desired measure $\mu_{1/2,\lambda}$;
and for $\beta > 1/2$, \reff{lemma.compmono.series.eq3} together with
\be
   (u+v)^{-(\beta-\!\smhalf)}  \;=\;
   \int\limits_0^\infty 
       {t^{\beta-\!\smfrac{3}{2}} \over \Gamma(\beta-\smhalf)} \,
       e^{-t(u+v)}
       \; dt
 \label{eq.lemma.compmono.series.eq5}
\ee
represents $\mu_{\beta,\lambda}$ as the convolution of two positive measures.
Indeed, multiplying \reff{lemma.compmono.series.eq3} by
\reff{eq.lemma.compmono.series.eq5},
one obtains after a straightforward change of variables the explicit formula
\be
   \mu_{\beta,\lambda}(t_1,t_2)
   \;=\;
   {(4\lambda)^{1-\beta}  \over \Gamma(\smhalf) \, \Gamma(\beta-\smhalf)}
   \:
   P(t_1,t_2,\lambda)^{\! \beta -\!\smfrac{3}{2}}
   \:
   \chi(t_1,t_2,\lambda)
 \label{eq.lemma.compmono.series.eq6}
\ee
where
\begin{eqnarray}
   P(t_1,t_2,\lambda)
   & = &
   2 (t_1 t_2 + \lambda t_1 + \lambda t_2) \,-\, (t_1^2 + t_2^2 + \lambda^2)
 \label{eq.lemma.compmono.series.eq6.bis1}
\end{eqnarray}
and
\begin{eqnarray}
   \chi(t_1,t_2,\lambda)
   & = &
   \cases{ 1  & if $t_1,t_2 \ge 0$ and
                $|\sqrt{t_1} - \sqrt{t_2}| \le \sqrt{\lambda}
                 \le \sqrt{t_1} + \sqrt{t_2}$  \cr
           \noalign{\vskip 6pt}
           0  & otherwise \cr
         }
 \label{eq.lemma.compmono.series.eq6.bis2}
\end{eqnarray}
The constraint $\chi(t_1,t_2,\lambda) \neq 0$ states simply that
$\sqrt{t_1}, \sqrt{t_2}, \sqrt{\lambda}$ form the sides of a triangle;
and $P(t_1,t_2,\lambda)$ is precisely 16 times the square of the area
of this triangle
(Heron's formula \cite[Section~3.2]{Coxeter_67}).\footnote{
   It is curious that similar expressions,
   involving the area of a triangle in terms of its sides,
   arise also in Sonine's integral for the product of three Bessel functions
   \cite[p.~411, eq.~13.46(3)]{Watson_44}
   \cite[p.~36, eq.~(4.39) and p.~40]{Askey_75}
   --- a formula that Szeg\H{o} \cite{Szego_33} employed in one version of
   his nonnegativity proof for \reff{eq1.1} in the case $n=3$.
   Probably this is not an accident;
   it would be interesting to understand the precise connection
   between Sonine's formula and
   \reff{lemma.compmono.series.eq2}/\reff{eq.lemma.compmono.series.eq6}
   [cf.\ also \reff{eq.prop.K3} below].
}
In view of these explicit formulae,
the proof of Lemma~\ref{lemma.compmono.series}
is in fact completely elementary.

3.  It would be interesting to know whether Lemma~\ref{lemma.compmono.series}
can be generalized to other ratios of elementary symmetric polynomials,
e.g.\ $E_{r,n}^{-\beta} \exp(-\lambda E_{r+1,n}/E_{r,n})$.
By Lemma~\ref{lemma.compmono.invLT.affine} this would determine the
complete monotonicity of $E_{r+1,n+1}^{-\beta}$.
\medskip

\bigskip

\proofof{Proposition~\ref{prop.compmono.series}}
By Lemma~\ref{lemma.compmono.invLT.affine},
$B^{-\beta} \exp(-t_1 A/B)$ and $D^{-\beta} \exp(-t_2 C/D)$
are completely monotone on $\scrc$ for all $t_1,t_2 \ge 0$.
Hence the same is true of their product for any choice of $t_1,t_2 \ge 0$.
We now use the identity \reff{lemma.compmono.series.eq2},
multiplied on both sides by $(BD)^{-\beta}$,
with $u=A/B$ and $v=C/D$.
This shows that $(AD+BC)^{-\beta} \exp[-\lambda AC/(AD+BC)]$
is completely monotone on $\scrc$ for all $\lambda  \ge 0$.
But then using Lemma~\ref{lemma.compmono.invLT.affine} again
(this time in the reverse direction),
we conclude that $[AC+ (AD+BC)y]^{-\beta}$
is completely monotone on $\scrc \times (0,\infty)$.
\qed

Using Lemma~\ref{lemma.compmono.series}
we can also show that the spanning-tree polynomial of the 3-cycle
--- or equivalently, the elementary symmetric polynomial of degree 2
in three variables --- has the property that $P^{-\beta}$ is
completely monotone on $(0,\infty)^3$ {\em if and only if}\/ $\beta \ge 1/2$:

\begin{proposition}
   \label{prop.K3}
The function $F(x_1,x_2,x_3) = (x_1 x_2 + x_1 x_3 + x_2 x_3)^{-\beta}$
is completely monotone on $(0,\infty)^3$ if and only if
$\beta = 0$ or $\beta \ge 1/2$.
\end{proposition}

\proof
As always, $F$ is completely monotone for $\beta=0$
and not completely monotone for $\beta < 0$,
so it suffices to consider $\beta > 0$.
Using Lemma~\ref{lemma.compmono.invLT.affine} with
$A = x_1 x_2$, $B = x_1 + x_2$, $y=x_3$,
we see that $F$ is completely monotone on $(0,\infty)^3$ if and only if
the function $F_{\beta,\lambda}$ defined in \reff{lemma.compmono.series.eq1}
is completely monotone for all $\lambda \ge 0$.
But by Lemma~\ref{lemma.compmono.series},
this occurs if and only if $\beta \ge 1/2$.
\qed

\medskip\noindent
{\bf Remarks.}
1.  We shall later give two further independent proofs
of Proposition~\ref{prop.K3}:
one based on harmonic analysis on the cone of positive-definite
$m \times m$ real symmetric matrices specialized to $m=2$
[Corollary~\ref{cor.TG}, which follows from results to be proved
in Section~\ref{sec.det}, together with Proposition~\ref{cor2.det}(a)],
and one based on harmonic analysis on the Lorentz cone in $\R^n$
specialized to $n=3$ [Corollary~\ref{cor.quadratic},
which follows from results to be proved in Section~\ref{sec.quadratic}].
The point here is that the Jordan algebra
    ${\rm Sym}(2,\R) \simeq \R \times \R^2$
    can be viewed as a member of two different families of Jordan algebras:
    ${\rm Sym}(m,\R)$ and $\R \times \R^{n-1}$ \cite[p.~98]{Faraut_94}.

2. Proposition~\ref{prop.K3} implies that the property stated in
Proposition~\ref{prop.compmono.series} does {\em not}\/ hold
for $0 < \beta < 1/2$.
Indeed, it suffices to take $\scrc = (0,\infty)^2$
and $A=x_1$, $B=1$, $C=x_2$, $D=1$,
leading to the function $(x_1 x_2 + x_1 y + x_2 y)^{-\beta}$.

3. Proposition~\ref{prop.K3} is, of course,
just Theorem~\ref{thm1.serpar}${}'$ restricted to the 3-cycle $G=K_3$.
In particular it implies
Szeg\H{o}'s \cite{Szego_33} result (except the strict positivity)
for the polynomial \reff{eq1.1} in the special case $n=3$.

4. Combining \reff{eq.lemma.compmono.invLT.affine.proof} with
\reff{lemma.compmono.series.eq2}/\reff{eq.lemma.compmono.series.eq6},
we obtain for $\beta > 1/2$ the formula
\begin{eqnarray}
   & &
   (x_1 x_2 + x_1 x_3 + x_2 x_3)^{-\beta}
   \;\:=\;\:
   {4^{1-\beta}
    \over
    \Gamma(\smhalf) \, \Gamma(\beta-\smhalf) \, \Gamma(\beta)}
   \: \times
        \nonumber \\
   & & \qquad\qquad
   \int\limits_0^\infty
   \int\limits_0^\infty
   \int\limits_0^\infty
   \:
   e^{-t_1 x_1 - t_2 x_2 - t_3 x_3}
   \;
   P(t_1,t_2,t_3)^{\! \beta -\!\smfrac{3}{2}}_+
   \; dt_1 \, dt_2 \, dt_3
   \;, \qquad\quad
 \label{eq.prop.K3}
\end{eqnarray}
which provides an explicit elementary proof of
the direct (``if'') half of Proposition~\ref{prop.K3}.
See also Remark~1 in Section~\ref{subsec.det.major}
for an alternate derivation of \reff{eq.prop.K3},
and see Corollary~\ref{cor.quadratic.E2n}
for a generalization from $E_{2,3}$ to $E_{2,n}$.
\bigskip

The following generalization of Lemma~\ref{lemma.compmono.series}
is an analytic version of series extension of a single edge:

\begin{lemma}
   \label{lemma.serpar.compmono}
Let $V$ be a finite-dimensional real vector space and
let $\scrc$ be an open convex cone in $V$.
Let $f$ be completely monotone on $\scrc \times (0,\infty)$,
and let $\beta \ge 1/2$.
Then the function
\be
   f_\beta^{\#}(x,y,z)
   \;=\;
   (y+z)^{-\beta} \, f\biggl( x, {yz \over y+z} \biggr)
\ee
is completely monotone on $\scrc \times (0,\infty)^2$.
\end{lemma}

\proof
By the Bernstein--Hausdorff--Widder--Choquet theorem
(Theorem~\ref{thm.multidim.HBW})
and linearity, it suffices to prove the lemma for
$f(x,y) = \exp(- \langle \ell, x \rangle - \lambda y)$
with $\ell \in \scrc^*$ and $\lambda \ge 0$.
The variable $x$ now simply goes for the ride,
so that the claim follows immediately from Lemma~\ref{lemma.compmono.series}.
\qed


\section{Determinantal polynomials}
  \label{sec.det}

In this section we consider polynomials defined by determinants
as in \reff{def.P.det}.
We begin with some preliminary algebraic facts
about such determinantal polynomials.
After this, we turn to the analytic results that are our principal concern.
We first prove a simple abstract version of the half-plane property
for determinantal polynomials.
Then we turn to the main topic of this section,
namely, the proof of
Theorems~\ref{thm1.det}, \ref{thm1.det.cones} and \ref{thm1.det.cones.Jordan}.

\subsection{Algebraic preliminaries}

First, some notation:
We write $[m] = \{1,\ldots,m\}$.
If $A = (a_{ij})_{i,j=1}^m$ is an $m \times m$ matrix
and $I,J \subseteq [m]$,
we denote by $A_{IJ}$ the submatrix of $A$
corresponding to the rows $I$ and the columns $J$,
all kept in their original order.
We write $I^c$ to denote the complement of $I$ in $[m]$.
Then $\epsilon(I,J) = (-1)^{\sum_{i \in I} i + \sum_{j \in J} j}$
is the sign of the permutation that takes $I I^c$ into $J J^c$.

We begin with a simple formula for the determinant of the sum of two matrices,
which ought to be found in every textbook of matrix theory
but seems to be surprisingly little known.\footnote{
   This formula can be found in
   \cite[pp.~162--163, Exercise 6]{Marcus_75}
   and \cite[pp.~221--223]{Korepin_93}.
   It can also be found
   --- albeit in an ugly notation that obscures what is going on ---
   in
   \cite[pp.~145--146 and 163--164]{Marcus_75}
   \cite[pp.~31--33]{Rump_97}
   \cite[pp.~281--282]{Prells_03};
   and in an even more obscure notation in
   \cite[p.~102, item~5]{Aitken_56}.
   We would be grateful to readers who could supply additional references.
   The proof here is taken from \cite[Lemma~A.1]{Caracciolo_07}.

   We remark that, by the same method, one can prove
   a formula analogous to \reff{eq.detsum}
   in which all three occurrences of determinant are replaced by
   permanent and the factor $\epsilon(I,J)$ is omitted.
}

\begin{lemma}
   \label{lemma.detpoly1}
Let $A,B$ be $m \times m$ matrices with elements in a commutative ring $R$.
Then\footnote{
   The determinant of an empty matrix is of course defined to be 1.
   This makes sense in the present context even if the ring $R$
   lacks an identity element:  the term $I=J=\emptyset$ contributes
   $\det B$ to the sum \reff{eq.detsum}, while the term $I=J=[m]$
   contributes $\det A$.
}
\be
   \det(A+B)  \;=    \sum_{\begin{scarray}
                              I,J \subseteq [m] \\
                              |I| = |J|
                           \end{scarray}}
                     \epsilon(I,J) \,
                     (\det A_{IJ}) (\det B_{I^c J^c})
   \;.
 \label{eq.detsum}
\ee
\end{lemma}

\proof
Using the definition of determinant and expanding the products, we have
\begin{equation}
   \det(A+B)  \;=\;
   \sum_{\pi \in \scrs_m} \sgn(\pi)
   \sum_{I \subseteq [m]} \prod_{i \in I} a_{i \pi(i)}
                          \prod_{i' \in I^c} b_{i' \pi(i')}
   \;,
\end{equation}
where the outermost sum runs over all permutations $\pi$ of $[m]$.
Define now $J = \pi[I]$.  Then we can interchange the order of summation:
\begin{equation}
   \det(A+B)  \;=\;
   \sum_{\begin{scarray}
            I,J \subseteq [m] \\
            |I| = |J|
         \end{scarray}}
   \sum_{\begin{scarray}
            \pi \in \scrs_m \\
            \pi[I] = J
         \end{scarray}}
   \sgn(\pi)
   \prod_{i \in I} a_{i \pi(i)}
   \prod_{i' \in I^c} b_{i' \pi(i')}
   \;.
\end{equation}
Suppose now that $|I|=|J|=k$, and let us write
$I = \{i_1,\ldots,i_k\}$ and $J = \{j_1,\ldots,j_k\}$
where the elements are written in increasing order,
and likewise
$I^c = \{i'_1,\ldots,i'_{m-k}\}$ and $J = \{j'_1,\ldots,j'_{m-k}\}$.
Let $\pi' \in \scrs_k$ and $\pi'' \in \scrs_{m-k}$ be the permutations
defined so that
\begin{subeqnarray}
  \pi'(\alpha) = \beta    & \longleftrightarrow &   \pi(i_\alpha) = j_\beta \\
  \pi''(\alpha) = \beta   & \longleftrightarrow &   \pi(i'_\alpha) = j'_\beta
\end{subeqnarray}
It is easy to see that $\sgn(\pi) = \sgn(\pi') \sgn(\pi'') \epsilon(I,J)$.
The formula \reff{eq.detsum} then follows
by using twice again the definition of determinant.
\qed

The following special case is frequently useful:

\begin{corollary}
   \label{cor.detpoly1}
Let $A,B$ be $m \times m$ matrices with elements in a commutative ring $R$,
with at least one of them being a diagonal matrix.  Then
\be
   \det(A+B)  \;=    \sum_{I \subseteq [m]}
                     (\det A_{II}) (\det B_{I^c I^c})
   \;.
\ee
\end{corollary}

\proof
If $A$ is diagonal, then $\det A_{IJ} = 0$ whenever $I \neq J$,
and likewise for $B$.
\qed

{\bf Remark.}  We will see in Appendix~\ref{sec.moore}
that Corollary~\ref{cor.detpoly1} generalizes to quaternions,
but Lemma~\ref{lemma.detpoly1} does not:
see Proposition~\ref{prop.detsum.quaternions} and the remark following it.
\bigskip

Iterating Lemma~\ref{lemma.detpoly1}, we obtain a formula
for the determinant of a sum of $n$ matrices:
\be
   \det\!\left( \sum_{k=1}^n  A_k \right)
   \;=\;
   \sum_{ {\bf I}, {\bf J} }
        \epsilon({\bf I},{\bf J}) \,
        \prod_{k=1}^n  \det[(A_k)_{I_k J_k}]
   \;,
 \label{eq.detsum.iterated}
\ee
where the sum runs over ordered partitions
${\bf I} = (I_1,\ldots,I_n)$ and
${\bf J} = (J_1,\ldots,J_n)$ of $[m]$
into $n$ {\em possibly empty}\/ blocks
satisfying $|I_k| = |J_k|$ for all $k$;
here $\epsilon({\bf I},{\bf J})$
is the sign of the permutation taking
$I_1 I_2 \cdots I_n$ into $J_1 J_2 \cdots J_n$.

We can now say something about determinantal polynomials
of the type \reff{def.P.det}.
Recall that if $A$ is a (not-necessarily-square) matrix
with elements in a commutative ring $R$,
then the {\em (determinantal) rank}\/ of $A$
is defined to be the largest integer $r$
such that $A$ has a nonzero $r \times r$ minor;
if no such minor exists (i.e., $A=0$), we say that $A$ has rank 0.

\begin{proposition}
   \label{prop.detpoly2}
Let $A_1,\ldots,A_n$ be $m \times m$ matrices with elements
in a commutative ring $R$,
and let $x_1,\ldots,x_n$ be indeterminates.
Then
\be
   P(x_1,\ldots,x_n)  \;=\; \det\!\left( \sum_{i=1}^n x_i A_i \right)
 \label{def.P.det.bis}
\ee
is a homogeneous polynomial of degree $m$ with coefficients in $R$.
Furthermore, the degree of $P$ in the variable $x_i$ is $\le {\rm rank}(A_i)$.
[In particular, if each $A_i$ is of rank at most 1,
 then $P$ is multiaffine.]
\end{proposition}

\proof
Both assertions about $P$ are immediate consequences of
\reff{eq.detsum.iterated}.
\qed

We are grateful to Andrea Sportiello for drawing our attention
to Lemma~\ref{lemma.detpoly1} and its proof,
and for showing us this elegant proof
of Proposition~\ref{prop.detpoly2}.

An analogue of Proposition~\ref{prop.detpoly2} holds also for
hermitian quaternionic matrices, albeit with a different proof:
here the determinant is the Moore determinant,
and ``rank'' means left row rank (= right column rank);
moreover, the polynomial $P(x_1,\ldots,x_n)$
is defined initially by letting $x_1,\ldots,x_n$ be {\em real}\/ numbers.
See Proposition~\ref{prop.detpoly2.quaternions}.

\subsection{The half-plane property}

Now we take an analytic point of view,
so that the commutative ring $R$ will be either $\R$ or $\C$.

In this subsection we make a slight digression from our main theme,
by showing that
if the $A_i$ are complex hermitian positive-semidefinite matrices
(this of course includes the special case of
 real symmetric positive-semidefinite matrices),
then the determinantal polynomial \reff{def.P.det.bis}
has the half-plane property.
This turns out to be an easy extension of the proof
of \cite[Theorem~8.1(a)]{Choe_hpp}.
Indeed, we can go farther,
by first stating the result in a clean abstract way,
and then deducing the half-plane property for \reff{def.P.det.bis}
as an immediate corollary.

\begin{proposition}
   \label{prop.det.hpp}
Let $V$ be the real vector space ${\rm Herm}(m,\C)$
of complex hermitian $m \times m$ matrices,
and let $C \subset V$ be the cone $\Pi_m(\C)$
of positive-definite matrices.
Then the polynomial $P(A) = \det A$
is nonvanishing on the tube $C+iV$.

[Of course, by restriction the same result holds when
$V$ is the real vector space ${\rm Sym}(m,\R)$
of real symmetric $m \times m$ matrices
and $C \subset V$ is the cone $\Pi_m(\R)$
of positive-definite matrices.]
\end{proposition}

\begin{corollary}
   \label{cor.det.hpp}
Let $A_1,\ldots,A_n$ be complex hermitian positive-semidefinite
$m \times m$ matrices.
Then the polynomial
\be
   P(x_1,\ldots,x_n)  \;=\; \det\!\left( \sum_{i=1}^n x_i A_i \right)
 \label{def.P.det.bis2}
\ee
has the half-plane property, that is,
either $P \equiv 0$ or else $P(x_1,\ldots,x_n) \neq 0$
whenever $\real x_i > 0$ for all $i$.
\end{corollary}

\firstproofof{Proposition~\ref{prop.det.hpp}}
Let $A \in C+iV$, and let $\psi$ be a nonzero vector in $\C^m$.
Then the Hermitian form
$\psi^* A \psi = \sum_{i,j=1}^m \psibar_i A_{ij} \psi_j$
has strictly positive real part,
and in particular is nonzero;
it follows that $A\psi \neq 0$.
Since this is true for every nonzero $\psi \in \C^m$,
we conclude that $\ker A = \{0\}$,
i.e.\ $A$ is nonsingular;
and this implies that (and is in fact equivalent to) $\det A \neq 0$.
\qed

\bigskip\noindent
{\sc Second Proof of Proposition~\ref{prop.det.hpp}}
\cite[Lemma~4.1]{Branden_07}.
Write $A = P+iQ$ with $P$ positive-definite and $Q$ hermitian.
Then $P$ has a positive-definite square root $P^{1/2}$,
and we have
\begin{subeqnarray}
   \det(P+iQ)  & = &
   \det[P^{1/2} (I + i P^{-1/2} Q P^{-1/2}) P^{1/2}]
        \\[1mm]
   & = &
   (\det P) \, \det(I + i P^{-1/2} Q P^{-1/2})
   \;.
\end{subeqnarray}
This is nonzero because all the eigenvalues of
$I + i P^{-1/2} Q P^{-1/2}$ have real part equal to 1.
\qed

\proofof{Corollary~\ref{cor.det.hpp}}
If at least one of the matrices $A_1,\ldots,A_n$
is strictly positive definite,
then $\sum_{i=1}^n x_i A_i$ lies in $C+iV$
whenever $\real x_i > 0$ for all $i$.
Proposition~\ref{prop.det.hpp} then implies that
$P(x_1,\ldots,x_n) \neq 0$.

The general case can be reduced to this one by replacing
$A_i \to A_i + \epsilon I$ ($\epsilon > 0$)
and taking $\epsilon \downarrow 0$.
By Hurwitz's theorem,
the limiting function is either
nonvanishing on the product of open right half-planes
or else is identically zero.
%
\qed

%
%

\medskip\noindent
{\bf Remarks.}
1.  The first proof of Proposition~\ref{prop.det.hpp}
is an abstraction of the proof of
Choe {\em et al.}\/ \cite[Theorem~8.1(a)]{Choe_hpp};
it can also be found in Faraut and Kor\'anyi \cite[Lemma~X.1.2]{Faraut_94}
and, for $V = {\rm Sym}(m,\R)$, in H\"ormander \cite[p.~85]{Hormander_90}.
Faraut and Kor\'anyi furthermore observe that
$A \mapsto A^{-1}$ is an involutive holomorphic automorphism
of the tube $C+iV$, having $I$ as its unique fixed point;
and these facts hold not only for ${\rm Herm}(m,\C)$
but in fact for any simple Euclidean Jordan algebra
\cite[Theorem~X.1.1]{Faraut_94}.

2.  It is natural to ask whether there is an analogous result for permanents,
which would extend \cite[Theorem~10.2]{Choe_hpp}
in the same way that Proposition~\ref{prop.det.hpp}
extends \cite[Theorem~8.1(a)]{Choe_hpp}.
The most obvious such extension would be:
Let $V$ be the vector space $\R^{m \times n}$
of real $m \times n$ matrices ($m \le n$)
and let $C \subset V$ be the cone $(0,\infty)^{m \times n}$
of matrices with strictly positive entries.
Then the polynomial $Q(A) = \per A$
is nonvanishing on the tube $C+iV$.
But this is simply {\em false}\/:
for instance, for $m=n=2$ we have
$\per \left(\! \begin{array}{cc}
                   1+i & 1-i \\
                   1-i & 1+i
               \end{array}  \!\right) = 0$.

\subsection{Proof of Theorems~\ref{thm1.det}, \ref{thm1.det.cones}
   and \ref{thm1.det.cones.Jordan}}
  \label{subsec.det.major}

Let us now turn to the main results of this section,
namely, Theorems~\ref{thm1.det.cones} and \ref{thm1.det.cones.Jordan}
and their consequence, Theorem~\ref{thm1.det}.
The proof of all these results turns out to be surprisingly easy;
and all but the ``only if'' part of Theorems~\ref{thm1.det.cones}
and \ref{thm1.det.cones.Jordan} is completely elementary.

In order to make the elementary nature of the proof as transparent as possible,
we proceed as follows:
First we prove the direct (``if'') half of 
Theorem~\ref{thm1.det.cones}(a,b) by completely elementary methods,
without reference to Euclidean Jordan algebras.
Then we prove Theorem~\ref{thm1.det.cones.Jordan}:
we will see that the proof of the direct half of this theorem
is a straightforward abstraction of the preceding elementary proof
in the concrete cases;
only the converse (``only if'') half is really deep.

\proofof{the direct half of Theorem~\ref{thm1.det.cones}(a,b)}
Let us begin with the real case and $\beta = 1/2$.
We use the Gaussian integral representation
\be
   (\det A)^{-1/2}
   \;=\;
   \int\limits_{\R^m}
   \exp( -\x^{\rm T} A \x )
   \,
   \prod\limits_{j=1}^m \displaystyle {d x_j \over \sqrt{\pi}}
 \label{eq.gaussian.real}
\ee
where $A$ is a real symmetric positive-definite $m \times m$ matrix
and we have written $\x = (x_1,\ldots,x_m) \in \R^m$.
This exhibits $(\det A)^{-1/2}$ as the Laplace transform
of a positive measure on $\Pi_m(\R)^* = \overline{\Pi_m(\R)}$,
namely, the push-forward of Lebesgue measure $d\x/\pi^{m/2}$ on $\R^m$
by the map $\x \mapsto \x \x^{\rm T}$.
[We remark that this measure is supported on positive-semidefinite matrices
 of rank 1.]
Alternatively, one can see directly,
by differentiating under the integral sign in \reff{eq.gaussian.real},
that the $k$-fold directional derivative of $(\det A)^{-1/2}$
in directions $B_1,\ldots,B_k \in \Pi_m(\R)$
has sign $(-1)^k$, because each derivative brings down a factor
$-\x^{\rm T} B_i \x \le 0$.

Since a product of completely monotone functions is completely monotone,
it follows immediately that $(\det A)^{-N/2}$ is completely monotone
for all positive integers $N$.
We remark that this can alternatively be seen
from the Gaussian integral representation
\be
   (\det A)^{-N/2}
   \;=\;
   \int\limits_{(\R^m)^N}
   \exp\biggl( - \sum\limits_{\alpha=1}^N \x^{(\alpha)\rm T}
                          A \x^{(\alpha)}
       \biggr)
   \,
   \prod\limits_{\alpha=1}^N
   \prod\limits_{j=1}^m \displaystyle {d x^{(\alpha)}_j \over \sqrt{\pi}}
 \label{eq.gaussian.real.Nspecies}
\ee
where we have introduced vectors $\x^{(\alpha)} \in \R^m$
for $\alpha=1,\ldots,N$.
If we assemble these vectors into an $m \times N$ matrix $X$,
then we have exhibited $(\det A)^{-N/2}$ as the Laplace transform
of a positive measure on $\overline{\Pi_m(\R)}$,
namely, the push-forward of Lebesgue measure $dX/\pi^{mN/2}$
on $\R^{m \times N}$
by the map $X \mapsto X X^{\rm T}$.
This measure is supported on positive-semidefinite matrices of rank
$\min(N,m)$.

Finally, for real values of $\beta > (m-1)/2$,
we use the integral representation\footnote{
   See e.g.\ \cite{Ingham_33},
       \cite[pp.~585--586]{Siegel_35},
       \cite[p.~41]{Terras_88},
       \cite[Lemma~1]{Etingof_99},
       or \cite[Theorem~7.2.2 and Corollary~7.2.4]{Anderson_03}.
}
\be
   (\det A)^{-\beta}
   \;=\:
   \left[ \pi^{m(m-1)/4} \prod\limits_{j=0}^{m-1}
              \Gamma\Bigl( \beta - {j \over 2} \Bigr)
   \right] ^{\! -1}
   \times \int\limits_{B > 0}
       e^{-\tr(AB)} \, (\det B)^{\beta -\! \smfrac{m+1}{2}} \, dB
   \;,
   \quad
 \label{eq.traceintegral.real}
\ee
where the integration runs over
real symmetric positive-definite $m \times m$ matrices $B$,
with measure $dB = \prod\limits_{1 \le i \le j \le m} dB_{ij}$.
This exhibits $(\det A)^{-\beta}$ as the Laplace transform
of a positive measure on $\Pi_m(\R)$.

The proof in the complex case is completely analogous.
For $\beta=1$ we use the Gaussian integral representation
\be
   (\det A)^{-1}
   \;=\;
   \int\limits_{\C^m}
   \exp( - \zbar^{\rm T} A \z )
   \,
   \prod\limits_{j=1}^m
      \displaystyle {(d \real z_j)(d \imag z_j) \over \pi}
 \label{eq.gaussian.complex}
\ee
where $A$ is a complex hermitian positive-definite matrix,
$\z = (z_1,\ldots,z_m) \in \C^m$
and $\bar{\hphantom{\z}}$ denotes complex conjugation.
For real values of $\beta > m-1$, we use the integral representation\footnote{
   See e.g.\ \cite[eqns.~(5.13)--(5.40)]{Goodman_63}.
   Please note that \cite[Problem~7.9]{Anderson_03} and
   \cite[eq.~(1.2)]{Graczyk_03} appear to contain errors in the
   normalization constant.
}
\begin{eqnarray}
   (\det A)^{-\beta}
   \;=\;
   \left[ \pi^{m(m-1)/2} \prod\limits_{j=0}^{m-1}
              \Gamma( \beta - j )
   \right] ^{\! -1}
   \times \int\limits_{B > 0}
       e^{-\tr(AB)} \, (\det B)^{\beta - m} \, dB
   \;,
   \qquad
 \label{eq.traceintegral.complex}
\end{eqnarray}
where the integration runs over
complex hermitian positive-definite $m \times m$ matrices $B$,
with measure
\be
   dB  \;=\;  \prod_{i=1}^m dB_{ii} \,
              \prod_{1 \le i < j \le m}  (d \real B_{ij})(d \imag B_{ij})
   \;.
  \vspace*{-3mm}
 \label{def.DB.complex}
\ee
\qed

{\bf Remarks.}
1.  By applying \reff{eq.traceintegral.real} for the case $m=2$ to
$A = \displaystyle  \left(\!\!\! \begin{array}{cc}
                                    x_1 + x_3   & x_3  \\
                                    x_3         & x_2 + x_3
                                 \end{array}
                    \!\!\!\right)$,
we obtain after a bit of algebra an alternate derivation
of the formula \reff{eq.prop.K3} for $(x_1 x_2 + x_1 x_3 + x_2 x_3)^{-\beta}$,
valid for $\beta > 1/2$.
In particular it implies the direct (``if'') half
of Proposition~\ref{prop.K3}, and hence also
Szeg\H{o}'s \cite{Szego_33} result (except the strict positivity)
for the polynomial \reff{eq1.1} in the special case $n=3$.

2. 
Similarly, by combining \reff{eq.traceintegral.complex}/\reff{def.DB.complex}
for the case $m=2$ with \reff{eq.E24}/\reff{eq.matrices.E24},
we obtain after some algebra an explicit formula for $E_{2,4}(\x)^{-\beta}$,
valid for $\beta > 1$:

\vspace*{-9mm}
\begin{eqnarray}
   & &
   (x_1 x_2 + x_1 x_3 + x_1 x_4 + x_2 x_3 + x_2 x_4 + x_3 x_4)^{-\beta}
   \;\:=\;\:
   {3^{\smfrac{3}{2}-\beta\vphantom{\biggl(}}
    \over
    2\pi \, \Gamma(\beta-1) \, \Gamma(\beta)}
   \: \times
        \nonumber \\
   & & \qquad
   \int\limits_0^\infty
   \int\limits_0^\infty
   \int\limits_0^\infty
   \int\limits_0^\infty
   \;
   e^{-t_1 x_1 - t_2 x_2 - t_3 x_3 - t_4 x_4}
   \;
   Q(t_1,t_2,t_3,t_4)_+^{\! \beta-2}
   \; dt_1 \, dt_2 \, dt_3 \, dt_4
 \label{eq.explicit.E24}
\end{eqnarray}
where
\be
   Q(t_1,t_2,t_3,t_4)
   \;=\;
   (t_1 t_2 + t_1 t_3 + t_1 t_4 + t_2 t_3 + t_2 t_4 + t_3 t_4)
   \,-\,
   (t_1^2 + t_2^2 + t_3^2 + t_4^2)
   \;.
 \label{def.Q.t14}
\ee
This formula provides an explicit elementary proof of
the direct half of Corollary~\ref{cor.E24}
--- and in particular solves the Lewy--Askey problem ---
in the same way that \reff{eq.prop.K3}
provides an explicit elementary proof of
the direct half of Proposition~\ref{prop.K3}.
Note also that by setting $x_4 = 0$ in \reff{eq.explicit.E24}
and performing the integral over $t_4$,
we obtain \reff{eq.prop.K3}.
Finally, see Corollary~\ref{cor.quadratic.E2n}
for a generalization from $E_{2,3}$ and $E_{2,4}$ to $E_{2,n}$.
\medskip

\bigskip

As preparation for the proof of Theorem~\ref{thm1.det.cones.Jordan},
let us review some facts from the theory of
harmonic analysis on Euclidean Jordan algebras
(see \cite{Faraut_94} for definitions and further background).

Let $V$ be a simple Euclidean (real) Jordan algebra
of dimension~$n$ and rank~$r$,
with Peirce subspaces $V_{ij}$ of dimension $d$;
recall that $n = r + \frac{d}{2} r(r-1)$.
It is illuminating (though not logically necessary)
to know that there are precisely five cases \cite[p.~97]{Faraut_94}:
\begin{itemize}
   \item[(a)] $V = {\rm Sym}(m,\R)$, the space of $m \times m$
      real symmetric matrices ($d=1$, $r=m$);
   \item[(b)] $V = {\rm Herm}(m,\C)$, the space of $m \times m$
      complex hermitian matrices ($d=2$, $r=m$);
   \item[(c)] $V = {\rm Herm}(m,\HH)$, the space of $m \times m$
      quaternionic hermitian matrices ($d=4$, $r=m$);
   \item[(d)] $V = {\rm Herm}(3,\OO)$, the space of $3 \times 3$
      octonionic hermitian matrices ($d=8$, $r=3$); and
   \item[(e)] $V = \R \times \R^{n-1}$ ($d=n-2$, $r=2$).
\end{itemize}
We denote by $(x|y) = \tr(xy)$ the inner product on~$V$,
where $\tr$ is the Jordan trace and $xy$ is the Jordan product.\footnote{
   In cases~(a) and (b), we have $(A|B) = \tr(AB)$;
   in case~(c) we have $(A|B) = \smhalf \tr(AB+BA)$;
   and in case~(e) we have
   $\big( (x_0,{\bf x}) | (y_0,{\bf y}) \big) =
    2(x_0 y_0 + {\bf x} \cdot {\bf y})$.
}

Let $\Omega \subset V$ be the positive cone
(i.e.\ the interior of the set of squares in $V$,
 or equivalently the set of invertible squares in $V$);
it is open, convex and self-dual.\footnote{
   In cases~(a)--(d) the positive cone $\Omega$ is the cone of
   positive-definite matrices;  in case~(e) it is the Lorentz cone
   $\{(x_0,{\bf x}) \colon\, x_0 > \sqrt{{\bf x}^2} \}$.
}
We denote by $\Delta(x)$ the Jordan determinant on~$V$:
it is a homogeneous polynomial of degree $r$ on~$V$,
which is strictly positive on $\Omega$
and vanishes on $\partial\Omega$.\footnote{
   In cases~(a) and (b), the Jordan determinant is the ordinary determinant;
   in case~(c) it is the Moore determinant (see Appendix~\ref{sec.moore});
   in case~(d) it is the Freudenthal determinant
   \cite{Freudenthal_54,Dray_98,Baez_02,Moldovan_09};
   and in case~(e) it is the Lorentz quadratic form
   $\Delta(x_0,{\bf x}) = x_0^2 - {\bf x}^2$.
}
We have the following fundamental Laplace-transform formula
\cite[Corollary~VII.1.3]{Faraut_94}:
for $y \in \Omega$ and
$\real\alpha > (r-1) \frac{d}{2} = \frac{n}{r} - 1$,
\be
   \int\limits_\Omega e^{-(x|y)} \, \Delta(x)^{\alpha - \frac{n}{r}} \, dx
   \;=\;
   \Gamma_\Omega(\alpha) \, \Delta(y)^{-\alpha}
 \label{eq.laplace}
\ee
where\footnote{
   Here $dx$ is Lebesgue measure on the Euclidean space $V$
   with inner product $(\,\cdot\,|\,\cdot\,)$.
   Thus in case~(a) we have
   $(B|B) = \sum\limits_{i=1}^m B_{ii}^2 +
            2 \!\!\!\sum\limits_{1 \le i < j \le m} B_{ij}^2$,
   so that
   $$ dx
      \;=\;
      \prod\limits_{i=1}^m dB_{ii} \,
      \prod\limits_{1 \le i < j \le m} \! \sqrt{2} dB_{ij}
      \;=\;
      2^{m(m-1)/4} \, dB \;,$$
   showing that \reff{eq.laplace}/\reff{def.Gamma_Omega}
   agrees with \reff{eq.traceintegral.real}.
   Likewise, in case~(b) we have
   $(B|B) = \sum\limits_{i=1}^m B_{ii}^2 +
            2 \!\!\!\sum\limits_{1 \le i < j \le m} |B_{ij}|^2$,
   so that
   $$ dx
      \;=\;
      \prod\limits_{i=1}^m dB_{ii} \,
      \prod_{1 \le i < j \le m}  \! (\sqrt{2} d \real B_{ij})
                                 (\sqrt{2} d \imag B_{ij})
      \;=\;
      2^{m(m-1)/2} \, dB \;,$$
   showing that \reff{eq.laplace}/\reff{def.Gamma_Omega}
   agrees with \reff{eq.traceintegral.complex}.
}
\be
   \Gamma_\Omega(\alpha)
   \;=\;
   (2\pi)^{(n-r)/2} \,
     \prod\limits_{j=0}^{r-1} \Gamma\Bigl( \alpha - j \frac{d}{2} \Bigr)
   \;.
 \label{def.Gamma_Omega}
\ee

Thus, for $\real\alpha > (r-1) \frac{d}{2}$,
the function $\Delta(x)^{\alpha - \frac{n}{r}} / \Gamma_\Omega(\alpha)$
is locally integrable on $\overline{\Omega}$ and polynomially bounded;
it therefore defines a tempered distribution $\scrr_\alpha$ on~$V$
by the usual formula
\be
   \scrr_\alpha(\varphi)
   \;=\;
   {1 \over \Gamma_\Omega(\alpha)} \,
   \int\limits_\Omega \varphi(x) \, \Delta(x)^{\alpha - \frac{n}{r}} \, dx
      \qquad \hbox{for } \varphi \in \scrs(V)
   \;.
 \label{def.scrr}
\ee
Using \reff{eq.laplace}, a beautiful argument
--- which is a special case of I.~Bernstein's general method
for analytically continuing distributions of the form $\scrp_\Omega^\lambda$
\cite{Bernstein_72,Bjork_79} ---
shows that the distributions $\scrr_\alpha$
can be analytically continued to the whole complex $\alpha$-plane:

\begin{theorem} {$\!\!\!$ \bf \protect\cite[Theorem~VII.2.2 et seq.]{Faraut_94}
                 \ }
   \label{thm.basicproperties}
The distributions $\scrr_\alpha$ can be analytically continued
to the whole complex $\alpha$-plane
as a tempered-distribution-valued entire function of $\alpha$.
The distributions $\scrr_\alpha$ have support contained in $\overline{\Omega}$
and have the following properties:
\begin{subeqnarray}
   \scrr_0  & = & \delta   \\[1mm]
   \scrr_\alpha * \scrr_\beta  & = & \scrr_{\alpha+\beta}
      \slabel{eq.riesz.convolution}   \\[1mm]
   \Delta(\partial/\partial x) \, \scrr_\alpha & = & \scrr_{\alpha-1}
      \slabel{eq.riesz.bernstein}  \\[0mm]
   \Delta(x) \, \scrr_\alpha & = &
      \!\left( \prod\limits_{j=0}^{r-1} \Bigl( \alpha - j \frac{d}{2} \Bigr)
      \!\right) \scrr_{\alpha+1}
      \slabel{eq.riesz.product}
      \label{eq.riesz}
\end{subeqnarray}
(here $\delta$ denotes the Dirac measure at 0).
Finally, the Laplace transform of $\scrr_\alpha$ is
\be
   (\scrl \scrr_\alpha)(y)  \;=\; \Delta(y)^{-\alpha}
 \label{eq.laplace.analcont}
\ee
for $y$ in the complex tube $\Omega + iV$.\footnote{
   The property \reff{eq.riesz.product} is not explicitly stated in
   \cite{Faraut_94}, but for $\real\alpha > (r-1) \frac{d}{2}$
   it is an immediate consequence of \reff{def.Gamma_Omega}/\reff{def.scrr},
   and then for other values of $\alpha$
   it follows by analytic continuation
   (see also \cite[Proposition~3.1(iii) and Remark~3.2]{Hilgert_01}).
}
\end{theorem}

\noindent
The distributions $\{\scrr_\alpha\}_{\alpha \in \C}$
constructed in Theorem~\ref{thm.basicproperties}
are called the {\em Riesz distributions}\/
on the Euclidean Jordan algebra $V$.

It is fairly easy to find a {\em sufficient}\/ condition for a
Riesz distribution to be a positive measure:

\begin{proposition}
  {$\!\!\!$ \bf \protect\cite[Proposition~VII.2.3]{Faraut_94}
   (see also \protect\cite[Section~3.2]{Hilgert_01}
      \protect\cite{Lassalle_87,Bonnefoy-Casalis_90})\hspace*{-3mm}}
 \label{prop.positive}
\begin{itemize}
   \item[(a)]  For $\alpha = k \frac{d}{2}$ with $k=0,1,\ldots,r-1$,
        the Riesz distribution $\scrr_\alpha$ is a positive measure
        that is supported on the set of elements of $\overline{\Omega}$
        of rank exactly $k$ (which is a subset of $\partial\Omega$).
   \item[(b)]  For $\alpha > (r-1) \frac{d}{2}$,
        the Riesz distribution $\scrr_\alpha$ is a positive measure
        that is supported on $\Omega$
        and given there by a density (with respect to Lebesgue measure)
        that lies in $L^1_{\rm loc}(\overline{\Omega})$.
\end{itemize}
\end{proposition}

\noindent
Indeed, part (b) is immediate from the definition \reff{def.scrr},
while part (a) follows by reasoning that abstracts the constructions
given in \reff{eq.gaussian.real}/\reff{eq.gaussian.real.Nspecies} and
\reff{eq.gaussian.complex} above for the special cases of
real symmetric and complex hermitian matrices.\footnote{
   Thus, for integer $N \ge 0$
   in the real symmetric (resp.\ complex hermitian) case,
   the positive measure $R_{N/2}$ is supported
   on the positive-semidefinite matrices of rank $\min(N,m)$
   \cite[Proposition~VII.2.3]{Faraut_94}
   and is nothing other than the push-forward of Lebesgue measure
   on $\R^{m \times N}$ (resp.\ $\C^{m \times N}$)
   by the map $X \mapsto X X^{\rm T}$ (resp.\ $X \mapsto X X^*$),
   as discussed above during the proof of the direct half
   of Theorem~\ref{thm1.det.cones}(a).
   For $N \ge m$ this is a straightforward calculation
   \cite{Speer_76,Blekher_82,Etingof_99,Muirhead_82,Anderson_03},
   which shows the equivalence of
   \reff{eq.gaussian.real.Nspecies} and \reff{eq.traceintegral.real}
   [or the corresponding formulae in the complex case];
   for $0 \le N \le m-1$ it follows by comparing
   \reff{eq.gaussian.real.Nspecies} with \reff{eq.laplace.analcont}
   and invoking the injectivity of the Laplace transform
   on the space of distributions \cite[p.~306]{Schwartz_66}:
   see \cite[Proposition~VII.2.4]{Faraut_94}.
}

It is a highly nontrivial fact that
the converse of Proposition~\ref{prop.positive} also holds:

\begin{theorem} {$\!\!\!$ \bf \protect\cite[Theorem~VII.3.1]{Faraut_94} \ }
   \label{thm.riesz1}
The Riesz distribution $\scrr_\alpha$
is a positive measure if and only if
$\alpha = 0,\frac{d}{2},\ldots,(r-1)\frac{d}{2}$
or $\alpha > (r-1)\frac{d}{2}$.
\end{theorem}

\noindent
This fundamental fact was first proven by Gindkin \cite{Gindikin_75}
(see also \cite{Berezin_75,Ishi_00}) and is generally considered to be deep.
However, there now exist two elementary proofs:
one that is a fairly simple but clever application of
Theorem~\ref{thm.basicproperties} and Proposition~\ref{prop.positive}
\cite{Shanbhag_88,Casalis_94} \cite[Appendix]{Sokal_riesz},
and another that analyzes the integrability of
$\Delta(x)^{\alpha - \frac{n}{r}}$ near $\partial\Omega$
and characterizes those $\alpha \in \C$ for which $\scrr_\alpha$
is a locally finite complex measure \cite{Sokal_riesz}.
In Appendix~\ref{sec.gindikin} below we give the first of these proofs,
thereby making the present paper nearly self-contained:
if one grants the elementary properties of the Riesz distributions
(Theorem~\ref{thm.basicproperties} and Proposition~\ref{prop.positive}),
then everything else is explicitly proven.

Using Proposition~\ref{prop.positive} for the case $d=4$,
we can prove the direct half of Theorem~\ref{thm1.det.cones}(c)
[i.e.\ the quaternionic case] by complete analogy with the
elementary proofs given above for the real and complex cases.
Moreover, with Theorem~\ref{thm.riesz1} in hand,
the proof of Theorem~\ref{thm1.det.cones.Jordan}
(and hence also of the converse half of Theorem~\ref{thm1.det.cones})
becomes utterly trivial:

\proofof{Theorem~\ref{thm1.det.cones.Jordan}}
Equation~\reff{eq.laplace.analcont} writes $\Delta(y)^{-\alpha}$
as the Laplace transform of the distribution $\scrr_\alpha$,
which is supported on $\overline{\Omega}$
(this is the closed dual cone of $\Omega$ since $\Omega$ is self-dual).
By the Bernstein--Hausdorff--Widder--Choquet theorem
(Theorem~\ref{thm.multidim.HBW}), it follows that
the map $y \mapsto \Delta(y)^{-\alpha}$ is completely monotone on $\Omega$
if and only if $\scrr_\alpha$ is a positive measure;
moreover, if $\scrr_\alpha$ is not a positive measure,
then this map is not completely monotone on
any nonempty open convex subcone $\Omega' \subseteq \Omega$.
So Theorem~\ref{thm1.det.cones.Jordan} is an immediate consequence
of Theorem~\ref{thm.riesz1}.
\qed

\bigskip

{\bf Remark.}
The formulae in this section arise in multivariate statistics
in connection with the Wishart distribution \cite{Muirhead_82,Anderson_03};
in recent decades some statisticians have introduced
the formalism of Euclidean Jordan algebras as a unifying device
\cite{Bonnefoy-Casalis_90,Casalis_94,Casalis_96,Massam_94,Massam_97}.
These formulae also arise in quantum field theory in studying the analytic
continuation of Feynman integrals to ``complex space-time dimension''
\cite{Speer_76,Blekher_82,Etingof_99,CSS_analcont}.
\medskip

\bigskip

Theorem~\ref{thm1.det.cones} is an immediate consequence of
Theorem~\ref{thm1.det.cones.Jordan},
once we use the fact that the Jordan determinant coincides
with the ordinary determinant on ${\rm Sym}(m,\R)$ and ${\rm Herm}(m,\C)$
and with the Moore determinant (see Appendix~\ref{sec.moore})
on ${\rm Herm}(m,\HH)$.
Theorem~\ref{thm1.det} is in turn an easy consequence
of Theorem~\ref{thm1.det.cones}:

\proofof{Theorem~\ref{thm1.det}}
If $A_1,\ldots,A_n$ are strictly positive definite,
Theorem~\ref{thm1.det} is an immediate corollary
of Theorem~\ref{thm1.det.cones}
and the definition of complete monotonicity on cones
(Definition~\ref{def.compmono.cones}).
The general case is reduced to this one by replacing
$A_i \to A_i + \epsilon I$ ($\epsilon > 0$)
and taking $\epsilon \downarrow 0$.
\qed


Conversely, if $\beta$ does {\em not}\/ lie in the set described in
Theorem~\ref{thm1.det.cones}, then the map $A \mapsto (\det A)^{-\beta}$
is {\em not}\/ completely monotone on {\em any}\/ nonempty open convex subcone
of the cone of positive-definite matrices.
This means, in particular, that if the matrices $A_1,\ldots,A_n$
span ${\rm Sym}(m,\R)$ or ${\rm Herm}(m,\C)$
[so that the convex cone they generate is open],
then the determinantal polynomial \reff{def.P.det.bis2}
does {\em not}\/ have $P^{-\beta}$ completely monotone on $(0,\infty)^n$.
See Section~\ref{subsec.matroids.converse}
for an application of this idea to graphs and matroids.

\section{Quadratic forms} 
       \label{sec.quadratic}

In this section we consider quadratic forms 
(= homogeneous polynomials of degree~2).
We begin by proving an abstract theorem giving
a necessary and sufficient condition
for such a quadratic form to be nonvanishing in a complex tube $C+iV$;
in the special case $C = (0,\infty)^n$
this corresponds to the half-plane property.
We then employ these results as one ingredient
in our proof of Theorem~\ref{thm1.quadratic.cones}.

\subsection{The half-plane property}

%
%

In this subsection we proceed in three steps.
First we study the analytic geometry associated to a
symmetric bilinear form $B$ on a finite-dimensional real vector space $V$
(Lemma~\ref{lemma.quadratic.1}).
Next we extend the quadratic form $Q(x) = B(x,x)$
to the complexified space $V+iV$ and study the values it takes
(Proposition~\ref{prop.quadratic}).
Finally we introduce the additional structure of an open convex cone
$C \subseteq V$ on which $Q$ is assumed strictly positive
(Corollary~\ref{cor.quadratic.C} and Theorem~\ref{thm.quadratic.hpp}).

So let $V$ be a finite-dimensional real vector space,
and let $B \colon\, V \times V \to \R$
be a symmetric bilinear form having inertia $(n_+,n_-,n_0)$.
Define $\scrs_+ = \{ x\colon\, B(x,x) > 0\}$
and $\scrs_- = \{ x\colon\, B(x,x) < 0\}$.
Clearly $\scrs_+$ and $\scrs_-$ are open cones
(not in general convex or even connected).
Indeed, $\scrs_+$ and $\scrs_-$ are {\em never}\/ convex
(except when they are empty)
because $x \in \scrs_\pm$ implies $-x \in \scrs_\pm$
but manifestly $0 \notin \scrs_\pm$.

Many of our proofs will involve choosing a basis in $V$
(and hence identifying $V$ with $\R^n$)
in such a way that $B$ takes the form
\be
   B(x,y)  \;=\;
   \sum_{i=1}^{n_+} x_i y_i  \,-\, \sum_{i=n_+ +1}^{n_+ + n_-} x_i y_i
   \;.
 \label{eq.B.diagonalized}
\ee
Moreover, whenever $\scrs_+ \neq \emptyset$ we can choose the basis
in this construction such that the first coordinate direction
lies along any desired vector $x \in \scrs_+$:
that this can be done follows from the standard Gram--Schmidt proof
of the canonical form \reff{eq.B.diagonalized}.

Elementary analytic geometry gives:

\begin{lemma}
   \label{lemma.quadratic.1}
Let $V$ be a finite-dimensional real vector space,
and let $B$ be a symmetric bilinear form on $V$ having inertia $(n_+,n_-,n_0)$.
\begin{itemize}
   \item[(a)]  If $n_+ = 0$, then $\scrs_+ = \emptyset$.
   \item[(b)]  If $n_+ = 1$, then $\scrs_+$ is a nonempty disconnected open
      non-convex cone, and there exists a nonempty open convex cone $\scrc$
      such that $\scrs_+ = \scrc \cup (-\scrc)$
      and $\scrc \cap (-\scrc) = \emptyset$.
      [Clearly $\scrc$ is uniquely determined modulo a sign.]
      Moreover, $B(x,y) \ge B(x,x)^{1/2} B(y,y)^{1/2} > 0$
      whenever $x,y \in \scrc$;
      and $B(x,y)^2 \ge B(x,x) \, B(y,y)$
      whenever $x \in \overline{\scrs_+}$ and $y \in V$ (or vice versa).
   \item[(c)]  If $n_+ \ge 2$, then $\scrs_+$ is a nonempty connected open
      non-convex cone.
      Moreover, for each $x \in \scrs_+$, the set
\be
   \scrt_+(x)  \;=\;
   \{ y \in V \colon\; B(x,y)^2 < B(x,x) \, B(y,y) \}
\ee
is a nonempty open non-convex cone that is contained in $\scrs_+$
and has a nonempty intersection with every neighborhood of $x$;
moreover we can write
\be
   \scrt_+(x)  \;=\;
   \{ y \in V \colon\; \R x + \R y \hbox{ is a two-dimensional
        subspace contained within } \scrs_+ \cup \{0\} \, \}
\ee
[that is, it is a two-dimensional subspace on which $B$ is positive-definite].
\end{itemize}
Analogous statements hold for $\scrs_-$
when $n_- = 0$, $n_- = 1$ or $n_- \ge 2$.
\end{lemma}

\proof
%
%
(a) is trivial.

(b) Assume that $B$ takes the canonical form \reff{eq.B.diagonalized}
with $n_+ = 1$,
and define $\scrc$ to be the ``forward light cone''
\be
   \scrc  \;=\;
   \{ x \in \R^n \colon\;  x_1^2 - x_2^2 - \ldots - x_{n_-+1}^2 > 0
                             \hbox{ and } x_1 > 0 \}
   \;.
 \label{def.scrc}
\ee
It is immediate that $\scrs_+ = \scrc \cup (-\scrc)$
and $\scrc \cap (-\scrc) = \emptyset$,
and the statements about $\scrs_+$ follow easily.

Now consider any $x,y \in \scrc$ and define
\be
   g(\alpha) \;=\; B(x+\alpha y, x+\alpha y)  \;=\;
   B(x,x) \,+\, 2\alpha B(x,y) \,+\, \alpha^2 B(y,y)
   \;.
\ee
We have $g(0) = B(x,x) > 0$;
but for the special value $\alpha_\star = -x_1/y_1$
the vector $x+\alpha_\star y$ has its first component equal to zero
and hence $g(\alpha_\star) = B(x+\alpha_\star y, x+\alpha_\star y) \le 0$.
So the quadratic equation $g(\alpha) = 0$ has a real solution,
which implies that its discriminant is nonnegative,
i.e.\ that $B(x,y)^2 \ge B(x,x) \, B(y,y)$.

Next assume that $x \in \overline{\scrs_+}$ and $y \in V$.
If $B(x,x) = 0$ or $B(y,y) \le 0$, the assertion is trivial;
so we can assume that $x,y \in \scrs_+$.
By the replacements $x \to \pm x$ and $y \to \pm y$
(which do not affect the desired conclusion)
we may assume that $x,y \in \scrc$.
But in this case the desired inequality has already been proven.

Finally, using $B(x,y) > 0$ for $x,y \in \scrc$
it is easily checked that $\scrc$ is convex.

(c) Clearly $\scrs_+$ is nonempty;
and as explained earlier it is non-convex.
To prove that $\scrs_+$ is connected,
we can assume that $B$ takes the form \reff{eq.B.diagonalized} with $n_+ \ge 2$.
It is now sufficient to find a path in $\scrs_+$
from an arbitrary vector $x \in \scrs_+$ to the vector $e_1 = (1,0,\ldots,0)$.
But this is easy:  first move coordinates $x_i$ with $i > n_+$
monotonically to zero
[this increases $B(x,x)$ monotonically and hence stays in $\scrs_+$];
then rotate and scale inside the subspace spanned by the
first $n_+$ coordinates to obtain $e_1$.


Now assume that $B$ takes the canonical form \reff{eq.B.diagonalized}
with $n_+ \ge 2$ and with the given vector $x \in \scrs_+$
lying along the first coordinate direction.
Then an easy computation shows that
$y = (y_1,y_2,\ldots,y_n)$ belongs to $\scrt_+(x)$
if and only if $y' \equiv (0,y_2,\ldots,y_n)$ belongs to $\scrs_+$.
Therefore, the preceding results (b,c) applied with
$n_+$ replaced by $n_+ - 1$ show that $\scrt_+(x)$
is a nonempty open non-convex cone,
which is obviously contained in $\scrs_+$;
and by taking $y'$ small we see that $\scrt_+(x)$
meets every neighborhood of $x$.
Moreover, $\R x + \R y = \R x + \R y'$
is a two-dimensional subspace if and only if $y' \neq 0$
(i.e.\ $y \notin \R x$);
and since $B(x,y') = 0$, we have
$\R x + \R y' \subseteq \scrs_+ \cup \{0\}$
if and only if $y' \in \scrs_+ \cup \{0\}$.
\qed

{\bf Remark.}  Note the sharp contrast between parts (b) and (c):
in the latter case, given any $x \in \scrs_+$
there is a nonempty open cone of vectors $y$
satisfying the Schwarz inequality (strictly) with $x$;
while in the former case all vectors $y \in V$ satisfy the
{\em reverse}\/ Schwarz inequality with $x$.
\medskip

\bigskip

We now consider the quadratic form $Q(x) = B(x,x)$,
extended to the complexified space $V+iV$ in the obvious way:
$Q(x+iy) = B(x,x) - B(y,y) +2i B(x,y)$.
We want to study the values taken by $Q$
in the complex tubes $\scrs_+ + iV$ and $\scrs_- + iV$,
and in particular the presence or absence of zeros.
We write $\HH$ to denote the open right half-plane
$\{ \zeta \in \C \colon\, \real\zeta > 0\}$.

\begin{proposition}
   \label{prop.quadratic}
Let $V$ be a finite-dimensional real vector space,
let $B$ be a symmetric bilinear form on $V$ having inertia $(n_+,n_-,n_0)$,
and let $Q$ be the associated quadratic form extended to $V+iV$.
\begin{itemize}
   \item[(a)]  If $n_+ = 1$, then for every $x \in \scrs_+$
      we have $Q[\HH x] = Q[\scrs_+ + iV] = \C \setminus (-\infty,0]$.
      In particular, $Q$ is nonvanishing on $\scrs_+ + iV$.
   \item[(b)]  If $n_+ \ge 2$, then for every $x \in \scrs_+$
      and $y \in \scrt_+(x)$ [recall that $\scrt_+(x)$ is nonempty]
      we have $Q[\,[1,\infty)x + i (\R x + \R y)] = \C$.
      In particular, for each $x \in \scrs_+$ and $y \in \scrt_+(x)$
      there exists $z \in \R x + \R y$ such that $Q(x+iz) = 0$.
\end{itemize}
Analogous statements hold for $\scrs_-$ when $n_- = 1$ or $n_- \ge 2$.
\end{proposition}

\proof
(a) If $x \in \scrs_+$, then $Q(\zeta x) = \zeta^2 Q(x)$
can take any value in $\C \setminus (-\infty,0]$
as $\zeta$ ranges over $\HH$.
On the other hand, if $x \in \scrs_+$ and $y \in V$,
then $Q(x+iy) = B(x,x) - B(y,y) +2i B(x,y)$
cannot take a value in $(-\infty,0]$,
as $B(x,y) = 0$ implies by Lemma~\ref{lemma.quadratic.1}(b)
that $y \notin \scrs_+$, i.e.\ $B(y,y) \le 0$.

(b) Given $x \in \scrs_+$ and $y \in \scrt_+(x)$,
the vector $y' = y + \mu x$ with $\mu = -B(x,y)/B(x,x)$
satisfies $B(x,y') = 0$ and $Q(y') > 0$.
Therefore
\be
   Q(\lambda x + i(\alpha x + \beta y'))
   \;=\;
   (\lambda+i\alpha)^2 Q(x) \,-\, \beta^2 Q(y')
   \;,
\ee
and this is easily seen to take all complex values
as $\lambda$ ranges over $[1,\infty)$,
$\alpha$ over $\R$, and $\beta$ over $(0,\infty)$.
\qed

\bigskip

Now we introduce the additional structure of an open convex cone
$C \subseteq V$ on which $Q$ is assumed strictly positive
(i.e.\ $C \subseteq \scrs_+$).
The hypotheses in the following result are identical
to those of Theorem~\ref{thm1.quadratic.cones}.

\begin{corollary}
   \label{cor.quadratic.C}
Let $V$ be a finite-dimensional real vector space,
let $B$ be a symmetric bilinear form on $V$
having inertia $(n_+,n_-,n_0)$,
and define the quadratic form $Q(x) = B(x,x)$.
Let $C \subset V$ be a nonempty open convex cone
with the property that $Q(x) > 0$ for all $x \in C$.
Then $n_+ \ge 1$, and moreover:
\begin{itemize}
   \item[(a)]  If $n_+ = 1$, then either $C \subseteq \scrc$
or $C \subseteq -\scrc$ [where $\scrc$ is defined as in
Lemma~\ref{lemma.quadratic.1}(b)],
and we have $B(x,y) > 0$ for all $x,y \in C$.
   \item[(b)]  If $n_+ \ge 2$, then for every $x \in C$
we have $Q[\,[1,\infty)x + iV] = \C$,
and in particular there exists $z \in V$ such that $Q(x+iz) = 0$.
\end{itemize}
\end{corollary}

\proof
(a) Let us use the canonical form \reff{eq.B.diagonalized}
with $n_+ = 1$ and $\scrc$ defined by \reff{def.scrc}.
By hypothesis the cone $C$ is contained within $\scrs_+ = \scrc \cup (-\scrc)$;
but since $C$ is convex we must in fact have either
$C \subseteq \scrc$ or $C \subseteq -\scrc$:
for otherwise $C$ would contain a point with $x_1 = 0$
and hence $B(x,x) \le 0$.\footnote{
   Alternatively, we can argue (assuming for simplicity that $n_0=0$) that
   a convex set $S$ contained in $\overline{\scrs_+}$
   must either lie on a single line through the origin
   (which is obviously impossible if $S$ is open)
   or else be contained in either
   $\overline{\scrc}$ or $-\overline{\scrc}$.
   For if we had $x \in S \cap \overline{\scrc}$ and
   $y \in S \cap -\overline{\scrc}$
   with $x$ and $y$ not on the same line through the origin,
   then there would exist $\lambda \in (0,1)$ such that
   $z = \lambda x + (1-\lambda) y$
   has $z_1 = 0$ and $z \neq 0$,
   hence $z \notin \overline{\scrs_+}$.
}
The remaining statements follow from Lemma~\ref{lemma.quadratic.1}(b)
and Proposition~\ref{prop.quadratic}(a).

(b) Since $C \subseteq \scrs_+$, the statements follow from
Lemma~\ref{lemma.quadratic.1}(c) and Proposition~\ref{prop.quadratic}(b).
\qed

\medskip

We can also summarize our results in a way that extends
\cite[Theorem~5.3]{Choe_hpp} from $C = (0,\infty)^n$
to general open convex cones $C$:

\begin{theorem}
   \label{thm.quadratic.hpp}
Let $V$ be a finite-dimensional real vector space,
let $B$ be a symmetric bilinear form on $V$
having inertia $(n_+,n_-,n_0)$,
and define the quadratic form $Q(x) = B(x,x)$.
Let $C \subset V$ be a nonempty open convex cone
with the property that $Q(x) > 0$ for all $x \in C$.
Then $n_+ \ge 1$, and the following are equivalent:
\begin{itemize}
   \item[(a)]  $n_+ = 1$.
   \item[(b)]  $Q$ is nonvanishing on the tube $C+iV$.
   \item[(c)]  $Q[C+iV] = \C \setminus (-\infty,0]$.
   \item[(d)]  If $x,y \in V$ with $Q(x) \ge 0$,
      then $B(x,y)^2 \ge B(x,x) \, B(y,y)$.
   \item[(e)]  If $x,y \in C$,
      then $B(x,y)^2 \ge B(x,x) \, B(y,y)$.
\end{itemize}
\end{theorem}

\proof
This follows immediately by putting together
the statements from Lemma~\ref{lemma.quadratic.1},
Proposition~\ref{prop.quadratic} and Corollary~\ref{cor.quadratic.C}:
if $n_+ = 1$, then (b,c,d,e) are all true;
and if $n_+ \ge 2$, then (b,c,d,e) are all false.
\qed

{\bf Remark.}  See \cite[Remark~1 after the proof of Theorem~5.3]{Choe_hpp}
for some of the history of this result in the traditional case
$C = (0,\infty)^n$.

\subsection{Proof of Theorem~\ref{thm1.quadratic.cones}}
       \label{subsec.quadratic.major}

\proofof{Theorem~\ref{thm1.quadratic.cones}}
We are concerned with the complete monotonicity of $Q^{-\beta}$,
where $Q(x) = B(x,x)$.
As always, it suffices to consider $\beta > 0$,
because complete monotonicity trivially holds when $\beta=0$
and never holds when $\beta < 0$ (because $Q$ grows at infinity).

We assume that $B$ takes the canonical form \reff{eq.B.diagonalized}
on $V = \R^n$, and we consider separately the three cases:

(a)
The case $n_+ = 1$, $n_- = 0$ is trivial:
we have $Q(x) = x_1^2$,
and the convex cone $C$ must be contained in one of the half-spaces
$\{x_1 > 0\}$ or $\{x_1 < 0\}$.
The map $x \mapsto Q(x)^{-\beta} = x_1^{-2\beta}$
is clearly completely monotone on each of these two half-spaces.

(b)
Next consider the case $n_+ = 1$, $n_- \ge 1$:
here \reff{eq.B.diagonalized} is the Lorentz form in
one ``timelike'' variable and $n_-$ ``spacelike'' variables,
and we have $\scrs_+ = \scrc \cup (-\scrc)$
where $\scrc$ is the forward light cone \reff{def.scrc}.
By Corollary~\ref{cor.quadratic.C}(a) we have either
$C \subseteq \scrc$ or $C \subseteq -\scrc$;
let us suppose the former.

Let us now show that if $\beta \ge (n_- -1)/2$,
then the map $x \mapsto Q(x)^{-\beta}$
is completely monotone on $\scrc$.
The variables $x_{n_- +2},\ldots,x_n$ play no role in this,
so we can assume without loss of generality that $n_0 = 0$, i.e.\ $n=n_- +1$.
For $\beta > (n_- -1)/2 = (n-2)/2$, the desired complete monotonicity
then follows from the integral representation
\cite[pp.~31--34]{Riesz_49}
\cite[eqns.~(24)--(28)]{Duistermaat_91}
\begin{eqnarray}
   & &
   (x_1^2 - x_2^2 - \ldots - x_n^2)^{-\beta}
   \;=\;
   \left[ \pi^{(n-2)/2} \, 2^{2\beta-1} \, \Gamma(\beta)
               \, \Gamma\Bigl( \beta - {n-2 \over 2} \Bigr)
   \right] ^{\! -1}
     \,\times \qquad\qquad   \nonumber \\
   & &
   \qquad\qquad\qquad\qquad
   \int\limits_{\scrc}
      e^{-(x_1 y_1 - x_2 y_2 - \ldots - x_n y_n)}  \,
      (y_1^2 - y_2^2 - \ldots - y_n^2)^{\beta- \! \smfrac{n}{2}} \, dy
 \label{eq.Riesz_integralrep}
\end{eqnarray}
valid for $\beta > (n-2)/2$,
which explicitly represents $Q(x)^{-\beta}$
as the Laplace transform of a positive measure
supported on the closed forward light cone $\scrc^* = \overline{\scrc}$.
For $\beta = (n-2)/2$ the result follows by taking limits.

Conversely, let us show that if $\beta \notin \{0\} \cup [(n_- -1)/2, \infty)$,
then the map $x \mapsto Q(x)^{-\beta}$
is not completely monotone on any nonempty open convex subcone
$C' \subseteq \scrc$.
For suppose that this map {\em is}\/ completely monotone on $C'$
for some such $\beta$:
then by the Bernstein--Hausdorff--Widder--Choquet theorem
(Theorem~\ref{thm.multidim.HBW}),
we must have
(assuming again without loss of generality that $n_0 = 0$)\footnote{
   If the map $x \mapsto Q(x)^{-\beta}$ is completely monotone on $C'$,
   then it is also completely monotone on the cone $C''$
   obtained by projecting $C'$ onto the first $n_+ + n_-$ coordinates
   (since $Q(x)$ is independent of the last $n_0$ coordinates).
}
\be
   (x_1^2 - x_2^2 - \ldots - x_n^2)^{-\beta}
   \;=\;
   \int
      e^{-(x_1 y_1 - x_2 y_2 - \ldots - x_n y_n)}  \,
      d\mu_\beta(y)
 \label{eq.laplace.mubeta}
\ee
for some positive measure $\mu_\beta$ supported on $(C')^*$.
Now, any such measure must clearly be Lorentz-invariant
and homogeneous of degree $2\beta-n$
(this follows from the injectivity of the Laplace transform).
Furthermore, $\mu_\beta$ must be supported on $\overline{\scrc}$,
for otherwise the support would contain a spacelike hyperboloid
$\{ y \in \R^n \colon\, y_1^2 - y_2^2 - \ldots - y_n^2 = \lambda \}$
for some $\lambda < 0$,
whose convex hull is all of $\R^n$
and hence not contained in the proper cone $(C')^*$.
On the other hand, every Lorentz-invariant locally-finite positive measure
on $\R^n$ that is supported on $\overline{\scrc}$ is of the form
\cite[Theorem~IX.33, pp.~70--76]{Reed-Simon_vol2}
\be
   \mu \;=\; c\delta_0 \,+\,  \int\limits_{m \ge 0} d\Omega_m \, d\rho(m)
 \label{eq.lorentz.mu}
\ee
where $c \ge 0$, $\delta_0$ denotes the point mass at the origin,
$\rho$ is a positive measure on $[0,\infty)$,
and $d\Omega_m = \delta(y^2 - m^2) \, dy$
is the unique (up to a constant multiple) Lorentz-invariant measure
on the ``mass hyperboloid''
\be
   H_m \;=\;
   \{ y \in \R^n \colon\, y_1^2 - y_2^2 - \ldots - y_n^2 = m^2
                                \hbox{ and } y_1 > 0 \}
   \;.
\ee
(When $n=2$, we consider $m > 0$ only, as there is no locally-finite
 Lorentz-invariant measure on $\R^2$ that is supported on $H_0$.)
A measure $\mu$ of the form \reff{eq.lorentz.mu}
is homogeneous in precisely three cases:
\begin{itemize}
   \item[(a)]  $c=0$, $d\rho(m) = {\rm const} \times m^{\lambda-1} \, dm$
      with $\lambda > 0$:  here $\mu$ is homogeneous of degree $\lambda-2$.
   \item[(b)]  $c=0$, $\rho = {\rm const} \times \delta_0$
      [for $n > 2$ only]:  here $\mu$ is homogeneous of degree $-2$.
   \item[(c)]  $c \ge 0$, $\rho=0$:
      here $\mu$ is homogeneous of degree $-n$.
\end{itemize}
This proves that a positive measure $\mu_\beta$ can exist {\em only if}\/
$\beta = 0$ or $\beta \ge (n-2)/2$.

(c) Finally, consider the case $n_+ > 1$.
By Corollary~\ref{cor.quadratic.C}(b),
$Q$ has zeros in the tube $C' + iV$
for every nonempty open convex subcone $C' \subseteq \scrs_+$
(indeed, for every nonempty {\em subset}\/ $C' \subseteq \scrs_+$).
We conclude by Corollary~\ref{cor.multidim.HBW} that
$Q^{-\beta}$ cannot be completely monotone on $C'$ for any $\beta > 0$.
\qed

{\bf Remark.}
More can be said about the integral representation \reff{eq.Riesz_integralrep}.
It turns out that the quantity
\be
   R_\beta  \;=\;
   \left[ \pi^{(n-2)/2} \, 2^{2\beta-1} \, \Gamma(\beta)
               \, \Gamma\Bigl( \beta - {n-2 \over 2} \Bigr)
   \right] ^{\! -1}
   (y_1^2 - y_2^2 - \ldots - y_n^2)^{\beta- \! \smfrac{n}{2}}
   \:
   {\rm I}[y \in \scrc]
   \;,
 \label{eq.Rbeta}
\ee
which is initially defined for $\beta > (n-2)/2$ as a positive measure
[or for $\real \beta > \mbox{$(n-2)/2$}$ as a complex measure]
on $\R^n$
(and which is of course supported on $\scrc$),
can be analytically continued as a
tempered-distribution-valued entire function of $\beta$
\cite{Duistermaat_91} \cite[Theorem~VII.2.2]{Faraut_94}:
this is the {\em Riesz distribution}\/ $\scrr_\beta$
on the Euclidean Jordan algebra $\R \times \R^{n-1}$.\footnote{
   A slightly different normalization is used in \cite{Faraut_94},
   arising from the fact that the Jordan inner product on
   $\R \times \R^{n-1}$ is
   $\big( (x_0,{\bf x}) | (y_0,{\bf y}) \big) =
    2(x_0 y_0 + {\bf x} \cdot {\bf y})$:
   this has the consequence that
   $dx_{\rm Jordan} = 2^{n/2} \, dx_{\rm ordinary}$,
   and also the Laplace transform is written with an extra factor 2
   in the exponential.
   The change of sign from $x_0 y_0 - {\bf x} \cdot {\bf y}$
   to $x_0 y_0 + {\bf x} \cdot {\bf y}$ is irrelevant,
   because the Riesz distribution $\scrr_\beta(y)$ is invariant
   under the reflections $y_i \mapsto -y_i$ for $2 \le i \le n$.
}
The integral representation
\be
   (x_1^2 - x_2^2 - \ldots - x_n^2)^{-\beta}
   \;=\;
   \int
      e^{-(x_1 y_1 - x_2 y_2 - \ldots - x_n y_n)}  \,
      \scrr_\beta(y) \, dy
 \label{eq.laplace.Rbeta}
\ee
where $x \in \scrc$
then holds for all complex $\beta$, by analytic continuation.
However, the distribution $\scrr_\beta$ is a positive measure if and only if
either $\beta = 0$ or $\beta \ge (n-2)/2$.
This follows from general results of harmonic analysis on
Euclidean Jordan algebras
(i.e.\ Theorem~\ref{thm.riesz1}),
but we have given here a direct elementary proof.
Indeed, once one has in hand the fundamental properties of
the Riesz distribution $\scrr_\beta$,
one obtains an even simpler elementary proof
by observing that \reff{eq.Rbeta} is not locally integrable
near the boundary of the cone $\scrc$ when $\real\beta \le (n-2)/2$
and $\beta \neq (n-2)/2$,
hence \cite[Lemma~2.1]{Sokal_riesz} that the distribution $\scrr_\beta$
is not a locally finite complex measure in these cases.

\subsection{Explicit Laplace-transform formula for inverse powers
               of a quadratic form}
   \label{subsec.quadratic.explicit}

By a simple change of variables,
we can generalize \reff{eq.Riesz_integralrep}
to replace the Lorentz matrix $L_n = \diag(1,-1,\ldots,-1)$
by an arbitrary real symmetric matrix $A$
of inertia $(n_+,n_-,n_0) = (1,n-1,0)$.
This will provide, among other things,
an explicit Laplace-transform formula for $E_{2,n}(\x)^{-\beta}$
that generalizes the formulae \reff{eq.prop.K3} and \reff{eq.explicit.E24}
obtained previously for $n=3,4$, respectively.

So let $A$ be a real symmetric $n \times n$ matrix
with one positive eigenvalue, $n-1$ negative eigenvalues,
and no zero eigenvalues.
We first need a slight refinement of Lemma~\ref{lemma.quadratic.1}(b)
to take advantage of the fact that we now have $n_0 = 0$;
for simplicity we state it in the ``concrete'' situation $V = \R^n$.

\begin{lemma}
   \label{lemma.quadratic.A}
Fix $n \ge 2$, and let $A$ be a real symmetric $n \times n$ matrix
with one positive eigenvalue, $n-1$ negative eigenvalues,
and no zero eigenvalues.
Then there exists a nonempty open convex cone $\scrc \subset \R^n$
(which is uniquely determined modulo a sign) such that
\begin{subeqnarray}
    \scrc \cap (-\scrc)                        & = &   \emptyset    \\[1mm]
    \overline{\scrc} \cap (-\overline{\scrc})  & = &   \{0\}     \\[1mm]
    \scrchat \cap (-\scrchat)                  & = &   \emptyset    \\[1mm]
    \overline{\scrchat} \cap (-\overline{\scrchat})  & = &   \{0\}     \\[1mm]
    \{ \y\colon\, \y^{\rm T} A \y > 0\}        & = &  \scrc \cup (-\scrc)  \\[1mm]
    \{ \x\colon\, \x^{\rm T} A^{-1} \x > 0\}   & = &
        \scrchat \cup (-\scrchat)
\end{subeqnarray}
where
\be
   \scrchat
   \;=\;
   \{ \x\colon\; \x^{\rm T} \y > 0
      \hbox{ \rm for all } \y \in \overline{\scrc} \setminus \{0\} \, \}
\ee
is the open dual cone to $\scrc$.
\end{lemma}

\proof
We can write $A = S^{\rm T} L_n S$ where $S$ is a nonsingular real matrix.
Then the claims follow easily from the corresponding properties
of the Lorentz quadratic form.
\qed

\begin{proposition}
   \label{prop.quadratic.A}
Fix $n \ge 2$, and let $A$ be a real symmetric $n \times n$ matrix
with one positive eigenvalue, $n-1$ negative eigenvalues,
and no zero eigenvalues;
and let $\scrc$ be the open convex cone from Lemma~\ref{lemma.quadratic.A}.
Then for $\beta > (n-2)/2$ we have
\begin{eqnarray}
   & &
   (\x^{\rm T} A^{-1} \x)^{-\beta}
   \;=\;
   \left[ \pi^{(n-2)/2} \, 2^{2\beta-1} \, \Gamma(\beta)
               \, \Gamma\Bigl( \beta - {n-2 \over 2} \Bigr)
   \right] ^{\! -1}
   \: |\det A|^{1/2}
     \,\times \qquad\qquad   \nonumber \\
   & &
   \qquad\qquad\qquad\qquad\qquad\qquad
   \int\limits_{\scrc}
      e^{-\x^{\rm T} \y}  \,
      (\y^{\rm T} A \y)^{\beta- \! \smfrac{n}{2}} \, d\y
 \label{eq.Riesz_integralrep.A}
\end{eqnarray}
for $\x \in \scrchat$.
\end{proposition}

\proof
Write $A = S^{\rm T} L_n S$ where $S$ is a nonsingular real matrix,
and make the changes of variable
$\y = S \y'$ and $\x = L_n S^{-\rm T} \x'$
in \reff{eq.Riesz_integralrep}.
Then $d\y = |\det S| \, d\y'$ where $|\det S| = |\det A|^{1/2}$;
and the formula \reff{eq.Riesz_integralrep.A}
follows immediately after dropping primes.
\qed

Let us now specialize to matrices $A$ of the form
$A = \lambda E_n - \mu I_n$,
where $I_n$ is the $n \times n$ identity matrix
and $E_n$ is the $n \times n$ matrix with all entries 1.
Then $A$ has eigenvalues $n\lambda-\mu,-\mu,\ldots,-\mu$,
hence has inertia $(n_+,n_-,n_0) = (1,n-1,0)$
provided that $\mu > 0$ and $\lambda > \mu/n$;
and in that case we have $A^{-1} = \lambda' E_n - \mu' I_n$ where
\be
   \lambda'  \;=\;  {\lambda \over \mu(n\lambda-\mu)}
   \;,
   \qquad
   \mu'  \;=\;  {1 \over \mu}
   \;.
\ee
[The map $(\lambda,\mu) \mapsto (\lambda',\mu')$ is of course involutive.]
Furthermore, we have
\be
   \y^{\rm T} A \y
   \;=\;
   2\lambda E_{2,n}(\y) \,+\, (\lambda-\mu) \|\y\|^2
\ee
where $E_{2,n}(\y) = \sum\limits_{1 \le i < j \le n} \! y_i y_j$
and $\|\y\|^2 = \sum\limits_{i=1}^n y_i^2$,
and analogously for $\x^{\rm T} A^{-1} \x$.

\begin{corollary}
   \label{cor.quadratic.A}
Fix $n \ge 2$, $\mu > 0$ and $\lambda > \mu/n$.
Then for $\beta > (n-2)/2$ we have
\begin{eqnarray}
   & &
   \Biggl( {2\lambda \over \mu(n\lambda-\mu)} E_{2,n}(\x)
           \,-\,
           {(n-1)\lambda - \mu  \over \mu(n\lambda-\mu)} \|\x\|^2
   \Biggr)^{\! -\beta}
   \;=\;
        \nonumber \\[3mm]
   & & \qquad\qquad
   \left[ \pi^{(n-2)/2} \, 2^{2\beta-1} \, \Gamma(\beta)
               \, \Gamma\Bigl( \beta - {n-2 \over 2} \Bigr)
   \right] ^{\! -1}
   \: \mu^{(n-1)/2} (n\lambda-\mu)^{1/2}
     \,\times \qquad\qquad   \nonumber \\
   & &
   \qquad\qquad
   \int\limits_{\scrc}
      e^{-\x^{\rm T} \y}  \,
      \Bigl( 2\lambda E_{2,n}(\y) \,+\, (\lambda-\mu) \|\y\|^2
      \Bigr)^{\beta- \! \smfrac{n}{2}} \, d\y
  \label{eq.cor.quadratic.A}
\end{eqnarray}
for $\x \in \scrchat$.
\end{corollary}

Specializing further to the case $\lambda = 2/(n-1)$, $\mu = 2$
(corresponding to $\lambda'=\mu' = 1/2$), we obtain:

\begin{corollary}[Laplace-transform formula for $E_{2,n}^{-\beta}$]
   \label{cor.quadratic.E2n}
Fix $n \ge 2$.
Then for $\beta > (n-2)/2$ we have
\be
   E_{2,n}(\x)^{-\beta}
   \;=\;
   {(n-1)^{\smfrac{n-1}{2} - \beta}
    \over
    (2\pi)^{(n-2)/2} \, \Gamma(\beta)
               \, \displaystyle \Gamma\Bigl( \beta - {n-2 \over 2} \Bigr)
   }
   \:
   \int\limits_{\scrc}
      e^{-\x^{\rm T} \y}  \,
      \Bigl( E_{2,n}(\y) \,+\, {n-2 \over 2} \|\y\|^2
      \Bigr)^{\beta- \! \smfrac{n}{2}} \, d\y
  \label{eq.cor.quadratic.E2n}
\ee
for $\x \in \scrchat \supseteq (0,\infty)^n$.
\end{corollary}

\noindent
For $n=2$ this is elementary;
for $n=3$ it reduces to \reff{eq.prop.K3};
and for $n=4$ it reduces to \reff{eq.explicit.E24}.
The formula \reff{eq.cor.quadratic.E2n} provides,
in particular, an explicit elementary proof of
the direct (``if'') half of Corollary~\ref{cor.quadratic}.

\section{Positive-definite functions on cones}  \label{sec.posdef}

In this section we recall briefly the theory of
positive-definite functions (in the semigroup sense) on convex cones,
which closely parallels the theory of completely monotone functions
developed in Section~\ref{sec.compmono}
and indeed can be considered as a natural extension of it.
We then apply this theory to powers of the determinant
on a Euclidean Jordan algebra,
and derive (in Theorem~\ref{thm.posdef.cones.Jordan})
a strengthening of Theorem~\ref{thm1.det.cones.Jordan}.
As an application of this latter result,
we disprove (in Example~\ref{exam.gurau})
a recent conjecture of Gurau, Magnen and Rivasseau \cite{Gurau_08}.

This section is not required for the application to graphs and matroids
(Section~\ref{sec.graphs+matroids}).

\subsection{General theory}

Here we summarize the basic definitions and results
from the theory of positive-definite functions on convex cones
and, more generally, on convex sets.
A plethora of useful additional information concerning positive-definite
(and related) functions on semigroups can be found in the
monograph by Berg, Christensen and Ressel \cite{Berg_84}.

\begin{definition}
   \label{def.posdef.cones}
Let $V$ be a real vector space,
and let $C$ be a convex cone in $V$.
Then a function $f \colon\, C \to \R$
is termed {\em positive-definite in the semigroup sense}\/
if for all $n \ge 1$ and all $x_1, \ldots, x_n \in C$,
the matrix $\{ f(x_i + x_j) \}_{i,j=1}^n$ is positive-semidefinite;
or in other words,
if for all $n \ge 1$, all $x_1, \ldots, x_n \in C$
and all $c_1,\ldots,c_n \in \C$ we have
\be
   \sum_{i,j=1}^n \overline{c_i} c_j f(x_i + x_j)  \;\ge\;  0  \;.
\ee
Similarly, a function $f \colon\, C+iV \to \C$
is termed {\em positive-definite in the involutive-semigroup sense}\/
if for all $n \ge 1$ and all $x_1, \ldots, x_n \in C+iV$,
the matrix
\break
$\{ {f(x_i + \overline{x_j})} \}_{i,j=1}^n$ is positive-semidefinite.
\end{definition}

\begin{theorem}
   \label{thm.multidim.posdef}
Let $V$ be a finite-dimensional real vector space,
let $C$ be an {\em open}\/ convex cone in $V$,
and let $f \colon\, C \to \R$.
Then the following are equivalent:
\begin{itemize}
   \item[(a)] $f$ is continuous and positive-definite in the semigroup sense.
   \item[(b)] $f$ extends to an analytic function on the tube $C+iV$
       that is positive-definite in the involutive-semigroup sense.
   \item[(c)] There exists a positive measure $\mu$ on $V^*$ satisfying
\be
   f(x)  \;=\;  \int\limits e^{-\langle \ell,x \rangle} \, d\mu(\ell)
 \label{eq.thm.multidim.posdef}
\ee
for all $x \in C$.
\end{itemize}
Moreover, in this case the measure $\mu$ is unique,
and the analytic extension to $C+iV$ is given by~\reff{eq.thm.multidim.posdef}.
%
\end{theorem}

Please note that the completely monotone functions
(Theorem~\ref{thm.multidim.HBW}) correspond to the subset of
positive-definite functions that are bounded at infinity
(in the sense that $f$ is bounded on the set $x+C$ for each $x \in C$),
or equivalently decreasing (with respect to the order induced by the cone $C$),
or equivalently for which the measure $\mu$
is supported on the closed dual cone $C^*$
[rather than on the whole space $V^*$ as in
 Theorem~\ref{thm.multidim.posdef}(c)].
See \cite[Lemma~1, p.~579]{Nussbaum_55} for a direct proof
that complete monotonicity implies positive-definiteness.

We remark that the hypothesis of continuity
(or at least something weaker, such as measurability or local boundedness)
in Theorem~\ref{thm.multidim.posdef}(a) is essential,
even in the simplest case $V=\R$ and $C=(0,\infty)$.
Indeed, using the axiom of choice it can easily be shown
\cite[pp.~35--36, 39]{Aczel_66}
that there exist discontinuous solutions to the functional equation
$\rho(x+y) = \rho(x) \rho(y)$ for $x,y \in (0,\infty)$,
and any such function is automatically positive-definite
in the semigroup sense.
However, any such function is necessarily non-Lebesgue-measurable
and everywhere locally unbounded \cite[pp.~34--35, 37--39]{Aczel_66}.

Theorem~\ref{thm.multidim.posdef} is actually a special case
of a more general theorem for open convex {\em sets}\/
that need not be cones.  We begin with the relevant definition
\cite[p.~x]{Glockner_03}:

\begin{definition}
   \label{def.posdef.convexsets}
Let $V$ be a real vector space.
If $C \subseteq V$ is a convex set,
then a function $f \colon\, C \to \R$
is termed {\em positive-definite in the convex-set sense}\/
if for all $n \ge 1$ and all $x_1, \ldots, x_n \in C$,
the matrix $\bigl\{ f\bigl( {x_i + x_j \over 2} \bigr) \bigr\}_{i,j=1}^n$
is positive-semidefinite.
More generally,
if $C \subseteq V+iV$ is a conjugation-invariant convex set,
then a function $f \colon\, C \to \C$
is termed {\em positive-definite in the involutive-convex-set sense}\/
if for all $n \ge 1$ and all $x_1, \ldots, x_n \in C$,
the matrix
$\bigl\{ f\bigl( {x_i + \overline{x_j} \over 2} \bigr) \bigr\}_{i,j=1}^n$
is positive-semidefinite.
\end{definition}

Note that if $C$ is in fact a convex {\em cone}\/,
then a function $f \colon\, C \to \R$
is positive-definite in the convex-set sense
if and only if it is positive-definite in the semigroup sense.
So this concept is a genuine generalization of the preceding one.

\begin{theorem}
   \label{thm.convex.posdef}
Let $V$ be a finite-dimensional real vector space,
let $C \subseteq V$ be an {\em open}\/ convex set,
and let $f \colon\, C \to \R$.
Then the following are equivalent:
\begin{itemize}
   \item[(a)] $f$ is continuous and positive-definite in the convex-set sense.
   \item[(b)] $f$ extends to an analytic function on the tube $C+iV$
       that is positive-definite in the involutive-convex-set sense.
   \item[(c)] There exists a positive measure $\mu$ on $V^*$ satisfying
\be
   f(x)  \;=\;  \int\limits e^{-\langle \ell,x \rangle} \, d\mu(\ell)
 \label{eq.thm.convex.posdef}
\ee
for all $x \in C$.
\end{itemize}
Moreover, in this case the measure $\mu$ is unique,
and the analytic extension to $C+iV$ is given by~\reff{eq.thm.convex.posdef}.
%
%
\end{theorem}

Theorem~\ref{thm.convex.posdef} was first proven
by Devinatz \cite{Devinatz_55},
using the spectral theory of commuting unbounded self-adjoint operators
on Hilbert space
(he gives details for $\dim V = 2$ but states that the methods work
 in any finite dimension);
see also Akhiezer \cite[pp.~229--231]{Akhiezer_65}
for the special case in which $C$ is a Cartesian product of open intervals.
A detailed alternative proof,
based on studying positive-definiteness on convex sets
of {\em rational}\/ numbers as an intermediate step \cite{Atanasiu_95},
has been given by Gl\"ockner
\cite[Proposition~18.7 and Theorem~18.8]{Glockner_03},
who also gives generalizations to infinite-dimensional spaces $V$
and to operator-valued positive-definite functions.
See also Shucker \cite[Theorem~4 and Corollary]{Shucker_84}
and Gl\"ockner \cite[Theorem~18.8]{Glockner_03}
for the very interesting extension to convex sets $C$
that are not necessarily open (but have nonempty interior):
in this latter case the representation \reff{eq.thm.convex.posdef}
does {\em not}\/ imply the continuity of $f$ on $C$,
but only on line segments (or more generally, closed convex hulls
of finitely many points) within $C$.
But with this modification the equivalence (a${}'$) $\iff$ (c) holds.

Surprisingly, we have been unable to find in the literature
any complete proof of Theorem~\ref{thm.multidim.posdef}
except as a corollary of the more general Theorem~\ref{thm.convex.posdef}.
But see \cite[Theorem~16.6]{Glockner_03}
for a version of Theorem~\ref{thm.multidim.posdef}
for the subclass of positive-definite functions that are
$\alpha$-bounded with respect to a ``tame'' absolute value $\alpha$.

It would be interesting to try to find simpler proofs of
Theorems~\ref{thm.multidim.posdef} and \ref{thm.convex.posdef}.

\subsection{Powers of the determinant on a Euclidean Jordan algebra}

We can now deduce analogues of
Theorems~\ref{thm1.det.cones} and \ref{thm1.det.cones.Jordan}
in which complete monotonicity is replaced by
positive-definiteness in the semigroup sense.
For brevity we state only the abstract result in terms of
Euclidean Jordan algebras.
The ``converse'' half of this result
constitutes an interesting strengthening of
the corresponding half of Theorem~\ref{thm1.det.cones.Jordan};
we will apply it in Example~\ref{exam.gurau}.

\begin{theorem}
  \label{thm.posdef.cones.Jordan}
Let $V$ be a simple Euclidean Jordan algebra of dimension~$n$ and rank~$r$,
with $n = r + \frac{d}{2} r(r-1)$;
let $\Omega \subset V$ be the positive cone;
let $\Delta \colon\, V \to \R$ be the Jordan determinant;
and let $\beta \in \R$.
Then the following are equivalent:
\begin{itemize}
   \item[(a)] The map $x \mapsto \Delta(x)^{-\beta}$ is positive-definite
      in the semigroup sense on $\Omega$.
   \item[(b)] The map $x \mapsto \Delta(x)^{-\beta}$ is positive-definite
      in the semigroup sense on some nonempty open convex subcone
      $\Omega' \subseteq \Omega$.
   \item[(c)] $\beta \in \{0,\frac{d}{2},\ldots,(r-1)\frac{d}{2}\}
               \cup \bigl((r-1)\frac{d}{2},\infty\bigr)$.
\end{itemize}
\end{theorem}

Theorem~\ref{thm.posdef.cones.Jordan}
is an immediate consequence of facts about Riesz distributions
--- namely, the Laplace-transform formula \reff{eq.laplace.analcont}
and Theorem~\ref{thm.riesz1} --- 
together with Theorem~\ref{thm.multidim.posdef}.
Indeed, the proof of Theorem~\ref{thm.posdef.cones.Jordan}
is essentially the {\em identical}\/ to that of
Theorem~\ref{thm1.det.cones.Jordan},
but using Theorem~\ref{thm.multidim.posdef} in place of the
Bernstein--Hausdorff--Widder--Choquet theorem.
The point, quite simply, is that our proof of the ``converse'' half
of Theorem~\ref{thm1.det.cones.Jordan} used only the failure of positivity
of the Riesz distribution,
not any failure to be supported on the closed dual cone
(indeed, it is {\em always}\/ supported there);
so it proves Theorem~\ref{thm.posdef.cones.Jordan} as well.

\begin{example}
   \label{exam.gurau}
\rm
Let $V$ be the real vector space ${\rm Sym}(m,\R)$
of real symmetric $m \times m$ matrices,
let $\Pi_m(\R) \subset V$ be the cone of positive-definite matrices,
and let $C \subset \Pi_m(R)$ be the subcone consisting of matrices
that are also elementwise strictly positive.
It follows from Theorem~\ref{thm.posdef.cones.Jordan} that
the map $A \mapsto (\det A)^{-\beta}$
is positive-definite in the semigroup sense on $C$
$\iff$ it is positive-definite in the semigroup sense on $\Pi_m(\R)$
$\iff$ $\beta \in \{0,\smhalf,1,\smfrac{3}{2},\ldots\} \cup [(m-1)/2,\infty)$.
This disproves the conjecture of Gurau, Magnen and Rivasseau
\cite[Section~7, Conjecture~1]{Gurau_08}
that the map $A \mapsto (\det A)^{-\beta}$
would be positive-definite in the semigroup sense on $C$
for all $\beta \ge 0$.
\end{example}

\section{Application to graphs and matroids}   \label{sec.graphs+matroids}

\subsection{Graphs}   \label{sec.application.graphs}

Let $G=(V,E)$ be a finite undirected graph
with vertex set $V$ and edge set $E$;
in this paper all graphs are allowed to have
loops and multiple edges unless explicitly stated otherwise.
Now let ${\bf x} = \{x_e\}_{e \in E}$
be a family of indeterminates indexed by the edges of $G$.
If $G$ is a {\em connected}\/ graph, we denote by $T_G({\bf x})$
the generating polynomial of spanning trees in $G$, namely
\be
   T_G({\bf x})  \;=\;   \sum_{T \in \scrt(G)} \, \prod_{e \in T} x_e
 \label{def.TG}
\ee
where $\scrt(G)$ denotes the family of edge sets of spanning trees in $G$.
If $G$ is {\em disconnected}\/, we define $T_G({\bf x})$ to be
the product of the spanning-tree polynomials of its connected components.
Otherwise put, $T_G({\bf x})$ is in all cases the generating polynomial
of {\em maximal spanning forests}\/ in $G$.
This is a slightly nonstandard definition
(the usual definition would put $T_G \equiv 0$ if $G$ is disconnected),
but it is convenient for our purposes and is natural from a
matroidal point of view (see below).
In order to avoid any possible misunderstanding,
we have inserted in Theorems~\ref{thm1.serpar} and \ref{thm1.serpar}${}'$
and Corollary~\ref{cor1.det} the word ``connected'',
so that the claims made in the Introduction
will be true on either interpretation of $T_G$.
Please note that, in our definition, $T_G$ is always strictly positive
on $(0,\infty)^E$, because the set of maximal spanning forests is nonempty.
Note also that, on either definition of $T_G$,
loops in $G$ (if any) play no role in $T_G$.
And it goes without saying that $T_G$ is a {\em multiaffine}\/ polynomial,
i.e.\ of degree at most 1 in each $x_e$ separately.

If $e$ is an edge of $G$, the spanning-tree polynomial of $G$
can be related to that of the deletion $G \setminus e$
and the contraction $G/e$:
\be
   T_G({\bf x})  \;=\;
   \cases{ T_{G \setminus e}({\bf x}_{\neq e}) + x_e T_{G/e}({\bf x}_{\neq e})
              & if $e$ is neither a bridge nor a loop  \cr
           \noalign{\vskip 3pt}
           x_e T_{G \setminus e}({\bf x}_{\neq e}) \;=\; 
           x_e T_{G/e}({\bf x}_{\neq e})
              & if $e$ is a bridge \cr
           \noalign{\vskip 3pt}
           T_{G \setminus e}({\bf x}_{\neq e}) \;=\; 
           T_{G/e}({\bf x}_{\neq e})
              & if $e$ is a loop \cr
         }
\ee
where ${\bf x}_{\neq e}$ denotes $\{x_f\}_{f \in E \setminus \{e\}}$.
The fact that $T_{G \setminus e}$ equals $T_{G/e}$
(rather than equalling zero) when $e$ is a bridge
is a consequence of our peculiar definition of $T_G$.


Now let us take an analytic point of view,
so that the indeterminates $x_e$ will be interpreted
as real or complex variables.

\begin{definition}
   \label{def.gbeta}
For each $\beta > 0$,
we denote by $\scrg_\beta$ the class of graphs $G$
for which $T_G^{-\beta}$ is completely monotone on $(0,\infty)^E$.
\end{definition}

The naturality of the classes $\scrg_\beta$ is illustrated by the
following easy but fundamental result:

\begin{proposition}
   \label{prop.gbeta.minors}
Each class $\scrg_\beta$ is closed under taking minors
(i.e.\ under deletion and contraction of edges
 and deletion of isolated vertices),
under disjoint unions, and under gluing at a cut vertex.
\end{proposition}

\proof
Deletion of a non-bridge edge $e$ corresponds to taking $x_e \downarrow 0$.
Contraction of a non-loop edge $e$ corresponds to dividing by $x_e$
and taking $x_e \uparrow +\infty$.
Both of these operations preserve complete monotonicity.
Deletion of a bridge has the same effect as contracting it,
in our peculiar definition of $T_G$.
Contraction of a loop is equivalent to deleting it
(but loops play no role in $T_G$ anyway).
Isolated vertices play no role in $T_G$.
This proves closure under taking minors.

If $G$ is obtained from $G_1$ and $G_2$ either by disjoint union
or by gluing at a cut vertex,
then $T_G = T_{G_1} T_{G_2}$ (on disjoint sets of variables)
in our definition of $T_G$;
this again preserves complete monotonicity.
\qed

Proposition~\ref{prop.gbeta.minors} illustrates the principal reason
for allowing arbitrary constants ${\bf c} > \zero$
(rather than just ${\bf c} = \one$)
in Theorem~\ref{thm1.serpar} and subsequent results:
it leads to a minor-closed class of graphs.
This, in turn, allows for characterizations that
are necessary as well as sufficient.
A similar situation arises in studying the negative-correlation property
for a randomly chosen basis of a matroid.
If only the ``uniformly-at-random'' situation is considered
(i.e., element weights ${\bf x} = \one$),
then the resulting class of matroids is not minor-closed,
and closure under minors has to be added by hand,
leading to the class of so-called {\em balanced matroids}\/
\cite{Feder-Mihail}.
But it then turns out that the class of balanced matroids
is {\em not}\/ closed under taking 2-sums \cite{Choe_Rayleigh}.
If, by contrast, one demands negative correlation for {\em all}\/
choices of element weights ${\bf x} > \zero$,
then the resulting class ---
the so-called {\em Rayleigh matroids}\/ ---
is automatically closed under taking minors
(by the same $x_e \to 0$ and $x_e \to \infty$ argument
 as in Proposition~\ref{prop.gbeta.minors}).
Moreover, it turns out to be closed under 2-sums as well \cite{Choe_Rayleigh}.

The very important property of closure under 2-sums holds also in our context.
To see this, let us first recall the definitions of
parallel connection, series connection and 2-sum of graphs
\cite[Section~7.1]{Oxley_11},
and work out how $T_G$ transforms under these operations.

For $i=1,2$, let $G_i = (V_i,E_i)$ be a graph
and let $e_i$ be an edge of $G_i$;
it is convenient (though not absolutely necessary)
to assume that $e_i$ is neither a loop nor a bridge in $G_i$.
Let us furthermore choose an orientation
$\vec{e}_i = \overrightarrow{x_i y_i}$ for the edge $e_i$.
(To avoid notational ambiguity, it helps to assume that the sets
 $V_1,V_2,E_1,E_2$ are all disjoint.)
Then the {\em parallel connection}\/ of
$(G_1, \vec{e}_1)$ with $(G_2, \vec{e}_2)$
is the graph $(G_1, \vec{e}_1) \| (G_2, \vec{e}_2)$
obtained from the disjoint union $G_1 \cup G_2$
by identifying $x_1$ with $x_2$, $y_1$ with $y_2$, and $e_1$ with $e_2$.
[Equivalently, it is obtained from the
 disjoint union $(G_1 \setminus e_1) \cup (G_2 \setminus e_2)$
 by identifying $x_1$ with $x_2$ (call the new vertex $x$),
 $y_1$ with $y_2$ (call the new vertex $y$),
 and then adding a new edge $e$ from $x$ to $y$.]
The {\em series connection}\/ of
$(G_1, \vec{e}_1)$ with $(G_2, \vec{e}_2)$
is the graph $(G_1, \vec{e}_1) \bowtie (G_2, \vec{e}_2)$
obtained from the disjoint union $(G_1 \setminus e_1) \cup (G_2 \setminus e_2)$
by identifying $x_1$ with $x_2$
and adding a new edge $e$ from $y_1$ to $y_2$.
The {\em 2-sum}\/ of $(G_1, \vec{e}_1)$ with $(G_2, \vec{e}_2)$
is the graph $(G_1, \vec{e}_1) \oplus_2 (G_2, \vec{e}_2)$
obtained from the parallel connection $(G_1, \vec{e}_1) \| (G_2, \vec{e}_2)$
by deleting the edge $e$ that arose from identifying $e_1$ with $e_2$,
or equivalently from the series connection
$(G_1, \vec{e}_1) \bowtie (G_2, \vec{e}_2)$
by contracting the edge $e$.

To calculate the spanning-tree polynomial of a parallel connection,
series connection or 2-sum, it is convenient to change slightly the notation
and suppose that $e_1$ and $e_2$ have already been identified
(let us call this common edge $e$),
so that $E_1 \cap E_2 = \{e\}$.
It is then not difficult to see
\cite[Proposition~7.1.13]{Oxley_11}
that the spanning-tree polynomial of a parallel connection $G_1 \|_e G_2$ is
given by
\be
   T_{G_1 \|_e G_2}({\bf x})  \;=\;
       T_{G_1 \setminus e} T_{G_2/e}  \,+\,
       T_{G_1/e} T_{G_2 \setminus e}  \,+\,
       x_e T_{G_1/e} T_{G_2/e}
   \;,
 \label{eq.TG.para}
\ee
while that of a series connection $G_1 \bowtie_e G_2$ is
\be
   T_{G_1 \bowtie_e G_2}({\bf x})  \;=\;
       T_{G_1 \setminus e} T_{G_2 \setminus e}  \,+\,
       x_e T_{G_1 \setminus e} T_{G_2/e}  \,+\,
       x_e T_{G_1/e} T_{G_2 \setminus e}
   \;.
 \label{eq.TG.series}
\ee
(All the spanning-tree polynomials on the right-hand sides
 are of course evaluated at ${\bf x}_{\neq e}$.)
The spanning-tree polynomial of a 2-sum $G_1 \oplus_{2,e} G_2$ is therefore
\be
   T_{G_1 \oplus_{2,e} G_2}({\bf x})  \;=\;
       T_{G_1 \setminus e} T_{G_2/e}  \,+\,
       T_{G_1/e} T_{G_2 \setminus e}
   \;.
 \label{eq.TG.2-sum}
\ee

\begin{proposition}
   \label{prop.gbeta.parallel}
Each class $\scrg_\beta$ is closed under parallel connection
and under 2-sums.
\end{proposition}

\proof
Closure under parallel connection is an immediate consequence
of Proposition~\ref{prop.compmono.parallel}
and the formula~\reff{eq.TG.para} for parallel connection.
Since the 2-sum is obtained from the parallel connection by deletion,
closure under 2-sum then follows from Proposition~\ref{prop.gbeta.minors}.
\qed

\begin{proposition}
   \label{prop.gbeta.series}
The class $\scrg_\beta$ is closed under series connection
for $\beta \ge 1/2$ but not for $0 < \beta < 1/2$.
\end{proposition}

\proof
For $\beta \ge 1/2$, closure under series connection 
is an immediate consequence of Proposition~\ref{prop.compmono.series}
and the formula~\reff{eq.TG.series} for series connection.
For $0 < \beta < 1/2$, non-closure under series connection
follows immediately from the observation that
the 2-cycle $C_2 = K_2^{(2)}$
(a pair of vertices connected by two parallel edges)
lies in $\scrg_\beta$ for all $\beta > 0$,
but the series connection of a 2-cycle with another 2-cycle
is a 3-cycle, which lies in $\scrg_\beta$ only for $\beta \ge \smhalf$
(by Proposition~\ref{prop.K3}).
\qed

\medskip\noindent
{\bf Remarks.}
%
1.  Unlike the situation for the half-plane and Rayleigh properties,
the classes $\scrg_\beta$ are {\em not}\/ in general closed
under duality for planar graphs.
For instance, the graph $C_3^* = K_2^{(3)}$
(a pair of vertices connected by three parallel edges)
lies in $\scrg_\beta$ for all $\beta \ge 0$;
but its dual $C_3$ lies in $\scrg_\beta$ only for $\beta \ge \smhalf$
(by Proposition~\ref{prop.K3}).

However, the class $\scrg_\beta$ \emph{is} duality-invariant
for $\beta \in (\smhalf,1)$,
as it consists of all series-parallel graphs (Theorem~\ref{thm.ghalf1} below).
And since $\scrg_\beta$ = all graphs 
for $\beta \in \{\smhalf,1,\smfrac{3}{2},\ldots\}$,
these classes $\scrg_\beta \cap \scrp$ (where $\scrp$ = planar graphs)
are also duality-invariant.
We do not know whether $\scrg_\beta \cap \scrp$ is duality-invariant
for $\beta \in (1,\infty) \setminus \{\smfrac{3}{2},2,\smfrac{5}{2},\ldots\}$.
In any case, the duality question is most naturally posed for matroids
rather than for graphs.

2.  Since $\scrg_\beta$ is closed under 0-sums (disjoint unions),
1-sums (gluing at a cut vertex) and 2-sums (essentially gluing at an edge),
it is natural to ask whether it is also closed under
3-sums (gluing along triangles).
We do not know the answer.
In particular, $K_4 \in \scrg_\beta$ for all $\beta \ge 1$
by Corollary~\ref{cor.TG},
but as noted in the discussion after Problem~\ref{problem.Gbeta}${}'$,
we do not know whether $K_5 - e = K_4 \oplus_3 K_4$
belongs to $\scrg_\beta$ for $\beta \in (1, {3 \over 2})$.
\medskip

\bigskip

It is an well-known (and easy) result that any minor-closed class of graphs
is of the form
\be
   {\rm Ex}(\scrf)  \;=\;
   \{ G \colon\, G \hbox{ does not contain any minor from } \scrf\}
\ee
for some family $\scrf$ of ``excluded minors'';
indeed, the minimal choice of $\scrf$
consists of those graphs that do not belong to the class in question
but whose proper minors all do belong to the class.
(Here we consider isomorphic graphs to be identical,
 or alternatively take only one representative from each isomorphism class.)

In one of the deepest and most difficult theorems of graph theory,
Roberston and Seymour \cite{Robertson_04} sharpened this result by proving
that any minor-closed class of graphs
is of the form ${\rm Ex}(\scrf)$
for some {\em finite}\/ family $\scrf$.
Therefore, each of our classes $\scrg_\beta$
can be characterized by a finite family of excluded minors.
One of the goals of this paper
--- alas, incompletely achieved --- is to determine these excluded minors.

\subsection{Matroids}  \label{sec.application.matroids}

The foregoing considerations have an immediate generalization to matroids.
(Readers unfamiliar with matroids can skip this subsection without loss of
 logical continuity.)
Let $M$ be a matroid with ground set $E$,
and let ${\bf x} = \{x_e\}_{e \in E}$
be a family of indeterminates indexed by the elements of $M$.
We denote by $B_M({\bf x})$
the basis generating polynomial of $M$, namely
\be
   B_M({\bf x})  \;=\;   \sum_{B \in \scrb(M)} \, \prod_{e \in B} x_e
\ee
where $\scrb(M)$ denotes the family of bases of $M$.
Please note that loops in $M$ (if any) play no role in $B_M$.
Note also that if $M$ is the graphic matroid $M(G)$
associated to a graph $G$, we have $B_{M(G)}({\bf x}) = T_G({\bf x})$.
This identity would not hold for disconnected $G$
if we had taken the standard definition of $T_G$.

\begin{definition}
   \label{def.mbeta}
For each $\beta > 0$,
we denote by $\scrm_\beta$ the class of matroids $M$
for which $B_M^{-\beta}$ is completely monotone on $(0,\infty)^E$.
\end{definition}

Once again we have:

\begin{proposition}
   \label{prop.mbeta.minors}
Each class $\scrm_\beta$ is closed under taking minors
(i.e.\ deletion and contraction)
and direct sums.
\end{proposition}

\proof
The proof is identical to that of Proposition~\ref{prop.gbeta.minors}
if one substitutes ``element'' for ``edge'', ``coloop'' for ``bridge'',
and ``direct sum'' for either form of union.
\qed

We refer to \cite[Section~7.1]{Oxley_11}
for the definitions of parallel connection, series connection
and 2-sum of matroids, which generalize those for graphs.
The upshot \cite[Proposition~7.1.13]{Oxley_11}
is that the formulae \reff{eq.TG.para}--\reff{eq.TG.2-sum}
for the spanning-tree polynomials of graphs
extend unchanged to the basis generating polynomials of matroids.
We therefore have:

\begin{proposition}
   \label{prop.mbeta.parallel}
Each class $\scrm_\beta$ is closed under parallel connection
and under 2-sums.
\end{proposition}

\begin{proposition}
   \label{prop.mbeta.series}
The class $\scrm_\beta$ is closed under series connection
for $\beta \ge 1/2$ but not for $0 < \beta < 1/2$.
\end{proposition}

Since each $\scrm_\beta$ is a minor-closed class,
we can once again seek a characterization of $\scrm_\beta$
by excluded minors.  However, in this case no analogue of
the Robertson--Seymour theorem exists,
so we have no {\em a priori}\/ guarantee of finiteness of the
set of excluded minors.
Indeed, there exist minor-closed classes of matroids
having an infinite family of excluded minors
\cite[Exercise~6.5.5(g)]{Oxley_11};
and in fact, for any infinite field $F$,
the class of $F$-representable matroids
has infinitely many excluded minors
\cite[Theorem~6.5.17]{Oxley_11}.

We suspect that the classes $\scrm_\beta$ are {\em not}\/ closed
under duality in general.
For instance, Corollary~\ref{cor.quadratic} shows that
$U_{2,5} \in \scrm_\beta$ if and only if $\beta \ge 3/2$;
but we suspect (Conjecture~\ref{conj.Ern}) that
$U_{3,5} \in \scrm_\beta$ if and only if $\beta \ge 1$.
On the other hand, we shall show in Theorem~\ref{thm.mhalf1}
that $\scrm_\beta$ for $1/2 < \beta < 1$ consists precisely
of the graphic matroids of series-parallel graphs ---
a class that {\em is}\/ closed under duality.

\subsection{Partial converse to Corollary~\ref{cor1.det}}
   \label{subsec.matroids.converse}

It was remarked at the end of Section~\ref{subsec.det.major}
that if $\beta$ does not lie in the set described in
Theorem~\ref{thm1.det.cones}, then the map $A \mapsto (\det A)^{-\beta}$
is not completely monotone on any nonempty open convex subcone
of the cone of positive-definite matrices;
and in particular, if the matrices $A_1,\ldots,A_n$
span ${\rm Sym}(m,\R)$ or ${\rm Herm}(m,\C)$,
then the determinantal polynomial \reff{def.P.det}/\reff{def.P.det.bis2}
does not have $P^{-\beta}$ completely monotone on $(0,\infty)^n$.
The spanning-tree polynomial of the complete graph $K_{m+1}$
provides an example of this situation;
and it turns out that there are two other cases
arising from complex-unimodular matroids.
The following result thus provides a (very) partial converse
to Corollary~\ref{cor1.det},
where part~(a) concerns regular [= real-unimodular] matroids
and part~(b) concerns complex-unimodular matroids:


\begin{proposition}
   \label{cor2.det}
Let $M$ be a matroid on the ground set $E$,
and let $B_M({\bf x})$ be its basis generating polynomial.
\begin{itemize}
   \item[(a)] If $M=M(K_p)$
      [the graphic matroid of the complete graph $K_p$]
      and $\beta \notin \{0,\smhalf,1,\smfrac{3}{2},\ldots\} \cup
          [(p-2)/2,\infty)$,
      then $B_M^{-\beta}$ is {\em not}\/ completely monotone on
      any nonempty open convex subcone of $(0,\infty)^E$.
   \item[(b${}_{r=2}$)]  If $M=U_{2,4}$
      [the uniform matroid of rank 2 on four elements]
      and $\beta \notin \{0\} \cup [1,\infty)$,
      then $B_M^{-\beta}$ is {\em not}\/ completely monotone on
      any nonempty open convex subcone of $(0,\infty)^E$.
   \item[(b${}_{r=3}$)]  If $M=AG(2,3)$
      [the ternary affine plane]
      and $\beta \notin \{0,1\} \cup [2,\infty)$,
      then $B_M^{-\beta}$ is {\em not}\/ completely monotone on
      any nonempty open convex subcone of $(0,\infty)^E$.
\end{itemize}
\end{proposition}
   
\proof
When $M$ is a real-unimodular (resp.\ complex-unimodular) matroid
of rank~$r$ with $n$ elements, we let $B$ be an $r \times n$
real-unimodular (resp.\ complex-unimodular) matrix
that represents $M$;
we then define $A_i$ ($1 \le i \le n$) to be the outer product of the
$i$th column of $B$ with its complex conjugate.
We shall show that in the cases enumerated above,
we can choose $B$ so that the matrices $A_1,\ldots,A_n$
span ${\rm Sym}(r,\R)$ [resp.\ ${\rm Herm}(r,\C)$].
The result then follows from Theorem~\ref{thm1.det.cones}
and the observation made immediately after it.

(a)  Let $H_p$ be the directed vertex-edge incidence matrix
for $K_p$ with an arbitrarily chosen orientation of the edges;
it is of size $p \times {p \choose 2}$ and is real-unimodular.
Then $K_p$ is represented over $\R$ by the matrix $H'_p$
obtained from $H_p$ by deleting one of the rows.
But we can reorder the columns of $H'_p$ to get $H''_p = (I_{p-1}|H_{p-1})$
where $I_{p-1}$ is the $(p-1) \times (p-1)$ identity matrix
(and $H_{p-1}$ is defined using the orientation of $K_{p-1}$
 inherited from $K_p$).
The corresponding matrices $A_1,\ldots,A_{{p \choose 2}}$,
obtained by taking outer products of the columns of $H''_p$
with their complex conjugates,
are easily seen to span ${\rm Sym}(p-1,\R)$.

(b${}_{r=2}$)  The matrices $A_1,\ldots,A_4$ defined in \reff{eq.matrices.E24}
are easily seen to span ${\rm Herm}(2,\C)$.

(b${}_{r=3}$)  Write $\omega = e^{\pm 2\pi i/3}$;
then the matrix \cite[p.~597]{Whittle_97}
\be
   B  \;=\;
   \left( \! \begin{array}{ccccccccc}
                1 & 0 & 0 & 1 & 0 & 1 & 1 & 1 & 1 \\
                0 & 1 & 0 & 1 & 1 & 0 & 1+\omega & 1 & 1+\omega \\
                0 & 0 & 1 & 0 & 1 & \omega & \omega & 1+\omega & 1+\omega
             \end{array}
          \! \right)
\ee
is easily seen to be complex-unimodular and to represent $AG(2,3)$.
A tedious computation
(or an easy one using {\sc Mathematica} or {\sc Maple})
now shows that the matrices $A_1,\ldots,A_9$ are linearly independent,
hence span the 9-dimensional space ${\rm Herm}(3,\C)$.
\qed

Let us remark that the cases enumerated in Proposition~\ref{cor2.det}
{\em exhaust}\/ the list of regular or complex-unimodular matroids
(of rank $r \ge 2$)
for which the matrices $A_1,\ldots,A_n$
span ${\rm Sym}(r,\R)$ or ${\rm Herm}(r,\C)$, respectively.
Indeed, it is known that a simple rank-$r$ matroid that is {\em regular}\/
(or, more generally, is binary with no $F_7$ minor)
can have at most $r(r+1)/2$ elements;
furthermore, the unique matroid attaining this bound is $M(K_{r+1})$
\cite[Proposition~14.10.3]{Oxley_11}.
See also \cite{Baclawski_79} for an intriguing proof
that uses the matrices $A_1,\ldots,A_n$
(but over $GF(2)$ rather than $\C$).
Likewise, it is known \cite[Theorem~2.1]{Oxley_98}
that a simple rank-$r$ matroid that is {\em complex-unimodular}\/
can have at most $(r^2+3r-2)/2$ elements if $r \neq 3$,
or 9 elements if $r=3$;
furthermore, the unique matroid attaining this bound
is $T_r$ (defined in \cite{Oxley_98}) when $r \neq 3$,
or $AG(2,3)$ when $r=3$.
The only cases in which this size reaches
$\dim {\rm Herm}(r,\C) = r^2$ are thus $r=1,2,3$,
yielding $T_1 = U_{1,1}$, $T_2 = U_{2,4}$ and $AG(2,3)$, respectively.



\subsection{Series-parallel graphs (and matroids):
    Proof of Theorem~\ref{thm1.serpar}${}'$}
  \label{sec.serpar}

Before proving Theorem~\ref{thm1.serpar}${}'$,
let us prove a similar but simpler theorem
concerning the interval $0 < \beta < \smhalf$.

\begin{theorem}
  \label{thm.g0half}
Let $G$ be a graph, and let $\beta \in (0,\smhalf)$.
Then the following are equivalent:
\begin{itemize}
   \item[(a)]  $G \in \scrg_\beta$.
   \item[(b)]  $G$ can be obtained from a forest
      by parallel extensions of edges
      (i.e., replacing an edge by several parallel edges)
      and additions of loops.
   \item[(c)]  $G$ has no $K_3$ minor.
\end{itemize}
Moreover, these equivalent conditions imply
that $G \in \scrg_{\beta'}$ for all $\beta' > 0$.
\end{theorem}

\proof
The equivalence of (b) and (c) is an easy graph-theoretic exercise.

If $G$ is obtained from a forest by parallel extensions of edges
and additions of loops,
then $T_G({\bf x})$ is a product of factors of the form
$x_{e_1} + \ldots + x_{e_k}$
(where $e_1,\ldots,e_k$ are a set of parallel edges in $G$),
so that $T_G^{-\beta}$ is completely monotone on $(0,\infty)^E$
for all $\beta \ge 0$.
Therefore (b) $\implies$ (a).

Conversely, Proposition~\ref{prop.K3} tells us that
$K_3 \notin \scrg_\beta$ for $\beta \in (0,\smhalf)$.
Since $\scrg_\beta$ is a minor-closed class,
this proves that (a) $\implies$ (c).
\qed

Dave Wagner has pointed out to us
that Theorem~\ref{thm.g0half} extends easily to matroids:

\begin{theorem}
  \label{thm.m0half}
Let $M$ be a matroid, and let $\beta \in (0,\smhalf)$.
Then the following are equivalent:
\begin{itemize}
   \item[(a)]  $M \in \scrm_\beta$.
   \item[(b)]  $M$ is the graphic matroid $M(G)$
      for a graph $G$ that can be obtained from a forest
      by parallel extensions of edges
      and additions of loops.
   \item[(c)]  $G$ has no $M(K_3)$ or $U_{2,4}$ minor.
\end{itemize}
Moreover, these equivalent conditions imply
that $M \in \scrm_{\beta'}$ for all $\beta' > 0$.
\end{theorem}

\proof
Tutte has proven \cite[Theorem~10.3.1]{Oxley_11}
that a matroid is graphic
if and only if it has no minor isomorphic to
$U_{2,4}$, $F_7$, $F_7^*$, $M^*(K_5)$ or $M^*(K_{3,3})$.
Since $M(K_3)$ is a minor of the last four matroids on this list,
the equivalence of (b) and (c) follows from the graphic case.

(b) $\implies$ (a) has already been proven in Theorem~\ref{thm.g0half}.

Finally, Proposition~\ref{prop.K3} tells us that
$M(K_3) \notin \scrm_\beta$ for $\beta \in (0,\smhalf)$,
and Corollary~\ref{cor.E24} or \ref{cor.quadratic}
tells us that
$U_{2,4} \notin \scrm_\beta$ for $\beta \in (0,1)$.
Since $\scrm_\beta$ is a minor-closed class,
this proves that (a) $\implies$ (c).
\qed

Let us now prove the corresponding characterization
for $\smhalf < \beta < 1$,
which is a rephrasing of Theorem~\ref{thm1.serpar}${}'$
and concerns series-parallel graphs.
Unfortunately, there seems to be no completely standard
definition of ``series-parallel graph'';
a plethora of slightly different definitions can be found in the literature
\cite{Duffin_65,Colbourn_87,Oxley_86,Oxley_11,Brandstadt_99}.
So let us be completely precise about our own usage:
we shall call a loopless graph {\em series-parallel}\/
if it can be obtained from a forest by a finite sequence of
series and parallel extensions of edges
(i.e.\ replacing an edge by two edges in series or two edges in parallel).
We shall call a general graph (allowing loops) series-parallel
if its underlying loopless graph is series-parallel.\footnote{
   Some authors write ``obtained from a tree'', ``obtained from $K_2$''
   or ``obtained from $C_2$'' in place of ``obtained from a forest'';
   in our terminology these definitions yield, respectively,
   all {\em connected}\/ series-parallel graphs,
   all connected series-parallel graphs whose blocks form a path,
   or all {\em 2-connected}\/ series-parallel graphs.
   See \cite[Section 11.2]{Brandstadt_99} for a more extensive bibliography.
}

So we need to understand how the spanning-tree polynomial $T_G({\bf x})$
behaves under series and parallel extensions of edges.
Parallel extension is easy:
if $G'$ is obtained from $G$ by replacing the edge $e$
by a pair of edges $e_1$ and $e_2$ in parallel, then
\be
   T_{G'}({\bf x}_{\neq e}, x_{e_1}, x_{e_2})
   \;=\;
   T_G({\bf x}_{\neq e}, x_{e_1} + x_{e_2})
   \;.
 \label{eq.parallel.TG}
\ee
In other words, two parallel edges with weights $x_{e_1}$ and $x_{e_2}$
are equivalent to a single edge with weight $x_{e_1} + x_{e_2}$.
This is because the spanning trees of $G'$
are in correspondence with the spanning trees $T$ of $G$ as follows:
if $T$ does not contain $e$, then leave $T$ as is
(it is a spanning tree of $G'$);
if $T$ does contain $e$, then adjoin to $T \setminus e$
one but not both of the edges $e_1$ and $e_2$.

Series extension is only slightly more complicated:
if $G'$ is obtained from $G$ by replacing the edge $e$
by a pair of edges $e_1$ and $e_2$ in series, then
\be
   T_{G'}({\bf x}_{\neq e}, x_{e_1}, x_{e_2})
   \;=\;
   (x_{e_1} + x_{e_2}) \,
   T_G \biggl( {\bf x}_{\neq e}, {x_{e_1} x_{e_2} \over x_{e_1} + x_{e_2}}
       \biggr)
   \;.
 \label{eq.series.TG}
\ee
In other words, two series edges with weights $x_{e_1}$ and $x_{e_2}$
are equivalent to a single edge with weight
$x_{e_1} x_{e_2}/(x_{e_1} + x_{e_2})$
together with a prefactor that clears the resulting denominator.
This is because the spanning trees of $G'$
are in correspondence with the spanning trees $T$ of $G$ as follows:
if $T$ does not contain $e$, then adjoin to $T \setminus e$
one but not both of the edges $e_1$ and $e_2$;
if $T$ does contain $e$, then adjoin to $T \setminus e$
both of the edges $e_1$ and $e_2$.
Since $T_G({\bf x}) = T_{G \setminus e}({\bf x}_{\neq e}) +
                      x_e T_{G/e}({\bf x}_{\neq e})$
where $T_{G \setminus e}$ (resp.\ $T_{G/e}$)
counts the spanning trees of $G$ that do not (resp.\ do) contain $e$
(the latter without the factor $x_e$),
this proves \reff{eq.series.TG}.

Let us remark that the parallel and series laws for $T_G({\bf x})$
are precisely the laws for combining electrical conductances
in parallel or series.
This is no accident, because as Kirchhoff \cite{Kirchhoff_1847} showed
a century-and-a-half ago,
the theory of linear electrical circuits can be written in terms
of spanning-tree polynomials
(see e.g.\ \cite{Chen_76} for a modern treatment).
Let us also remark that the parallel and series laws for $T_G({\bf x})$
are limiting cases of the parallel and series laws for the
multivariate Tutte polynomial $Z_G(q,{\bf v})$,
obtained when $q \to 0$ and ${\bf v}$ is infinitesimal:
see \cite[Sections~4.4--4.7]{Sokal_bcc2005} for a detailed explanation.

We are now ready to state and prove the main result of this section:

\begin{theorem}
  \label{thm.ghalf1}
Let $G$ be a graph, and let $\beta \in (\smhalf,1)$.
Then the following are equivalent:
\begin{itemize}
   \item[(a)]  $G \in \scrg_\beta$.
   \item[(b)]  $G$ is series-parallel.
   \item[(c)]  $G$ has no $K_4$ minor.
\end{itemize}
Moreover, these equivalent conditions imply
that $G \in \scrg_{\beta'}$ for all $\beta' \ge \smhalf$.
\end{theorem}

\noindent
Please note that Theorems~\ref{thm.g0half} and \ref{thm.ghalf1}
together imply Theorem~\ref{thm1.serpar}${}'$.

\proofof{Theorem~\ref{thm.ghalf1}}
The equivalence of (b) and (c) is a well-known graph-theoretic result
\cite[Exercise 7.30 and Proposition 1.7.4]{Diestel_10}
(see also \cite{Duffin_65,Oxley_86}).

Now let $G$ be a series-parallel graph:
this means that $G$ can be obtained from a forest by
series and parallel extensions of edges
and additions of loops.
As shown in Theorem~\ref{thm.g0half},
if $G$ is a forest, then $T_G^{-\beta}$
is completely monotone for all $\beta \ge 0$.
Parallel extension \reff{eq.parallel.TG} obviously preserves
complete monotonicity.
By Lemma~\ref{lemma.serpar.compmono},
series extension \reff{eq.series.TG} preserves
complete monotonicity whenever $\beta \ge \smhalf$.
Finally, additions of loops do not affect $T_G$.
Therefore, every series-parallel graph $G$
belongs to the class $\scrg_\beta$ for all $\beta \ge \smhalf$.

Conversely, Proposition~\ref{cor2.det}(a)
tells us that $K_4 \notin \scrg_\beta$ for $\beta \in (\smhalf,1)$.
Since $\scrg_\beta$ is a minor-closed class,
this proves (a) $\implies$ (c).
\qed

\medskip\noindent
{\bf Remarks.}
1.  Instead of using series and parallel extension of edges
[eqns.~\reff{eq.parallel.TG}/\reff{eq.series.TG}
 and Lemma~\ref{lemma.serpar.compmono}],
we could equally well have written this proof in terms of
the more general concept of series and parallel connection of graphs
[eqns.~\reff{eq.TG.para}/\reff{eq.TG.series}
 and Propositions~\ref{prop.compmono.parallel} and \ref{prop.compmono.series}].

2.  It ought to be possible to give an ``elementary'' proof
of the fact that $K_4 \notin \scrg_\beta$ for $\beta \in (\smhalf,1)$
--- and more generally of the fact that
the derivatives of $T_{K_4}^{-\beta}$ at ${\bf c} = \one$
do not all have sign $(-1)^k$ ---
by asymptotic calculation of coefficients
\`a la Pemantle--Wilson--Baryshnikov
\cite{Pemantle_02,Pemantle_04,Pemantle_08,Pemantle_10,Baryshnikov_11}.
\bigskip

Once again, Dave Wagner has pointed out to us
that Theorem~\ref{thm.ghalf1} extends immediately to matroids:

\begin{theorem}
  \label{thm.mhalf1}
Let $M$ be a matroid, and let $\beta \in (\smhalf,1)$.
Then the following are equivalent:
\begin{itemize}
   \item[(a)]  $M \in \scrm_\beta$.
   \item[(b)]  $M$ is the graphic matroid $M(G)$
     [or equivalently the cographic matroid $M^*(G)$] of
     a series-parallel graph $G$.
   \item[(c)]  $G$ has no $M(K_4)$ or $U_{2,4}$ minor.
\end{itemize}
Moreover, these equivalent conditions imply
that $M \in \scrm_{\beta'}$ for all $\beta' \ge \smhalf$.
\end{theorem}

\noindent
The proof is completely analogous to that of Theorem~\ref{thm.m0half}
(see also \cite[Corollary~12.2.14]{Oxley_11} for an alternative proof of
 (b) $\iff$ (c) in this case).

\bigskip

{\bf Remark.}
It follows from Corollary~\ref{cor1.det}(a) that
$\scrg_0 = \scrg_{1/2} = \scrg_1 = \scrg_{3/2} = \ldots =$ all graphs.
But what is the story for matroids?
Does $\scrm_{1/2}$ contain {\em only}\/ regular matroids?
Does $\scrm_{1}$ contain {\em only}\/ complex-unimodular matroids?
We suspect that the answer to this last question is {\em no}\/,
since we suspect that $U_{n-2,n} \in \scrm_{1}$ for all $n \ge 2$
(Conjecture~\ref{conj.Ern}).

\subsection{Combining the determinantal method with constructions}
   \label{subsec.combining}

Let us now combine the {\em ab initio}\/ results from
the determinantal method
(Theorems~\ref{thm1.det}--\ref{thm1.det.cones.Jordan}
 and their corollaries)
with the constructions from
Sections~\ref{sec.constructions}, \ref{sec.application.graphs}
and \ref{sec.application.matroids}
(deletion, contraction, direct sum, parallel connection and series connection).
For graphs we have:

\proofof{Proposition~\ref{prop.TG.extended}}
For $p=2$ the result is trivial, as the graphs concerned
are precisely those that can be obtained from a forest
by parallel extensions of edges and additions of loops
(Theorem~\ref{thm.ghalf1}).
For $p \ge 3$ (hence $\beta \ge 1/2$),
the result is an immediate consequence of Corollary~\ref{cor.TG} and
Propositions~\ref{prop.gbeta.minors}, \ref{prop.gbeta.parallel}
and \ref{prop.gbeta.series}.
\qed

The matroid generalization of Proposition~\ref{prop.TG.extended} is:

\begin{proposition}
   \label{prop.BM.extended1}
Fix $r \ge 1$,
and let $M$ be any matroid (on the ground set $E$)
that can be obtained from regular matroids of rank at most $r$ by
parallel connection, series connection, direct sum, deletion and contraction.
Then $B_M^{-\beta}$ is completely monotone on $(0,\infty)^E$
for $\beta = 0,\smhalf,1,\smfrac{3}{2},\ldots$
and for all real $\beta \ge (r-1)/2$.
\end{proposition}

\noindent
The proof is essentially identical to the previous one,
but uses Corollary~\ref{cor1.det}(a) in place of Corollary~\ref{cor.TG}
and Propositions~\ref{prop.mbeta.minors}--\ref{prop.mbeta.series}
in place of Propositions~\ref{prop.gbeta.minors}--\ref{prop.gbeta.series}.

And we also have:

\begin{proposition}
   \label{prop.BM.extended2}
Fix $r \ge 1$,
and let $M$ be any matroid (on the ground set $E$)
that can be obtained from regular matroids of rank at most $2r-1$
and complex-unimodular matroids of rank at most $r$ by
parallel connection, series connection, direct sum, deletion and contraction.
Then $B_M^{-\beta}$ is completely monotone on $(0,\infty)^E$
for $\beta = 0,1,2,3,\ldots$
and for all real $\beta \ge r-1$.
\end{proposition}

\noindent
The proof is again identical,
but uses both parts of Corollary~\ref{cor1.det} instead of only part~(a).

\bigskip

Propositions~\ref{prop.TG.extended}, \ref{prop.BM.extended1}
and \ref{prop.BM.extended2} give a rather abstract characterization
of the class of graphs or matroids that they handle,
and so it is of interest to give a more explicit characterization.
Let us start with the graph case.
We say that a graph $G$ is {\em minimally 3-connected}\/
if $G$ is 3-connected but, for all edges $e$ of $G$,
the graph $G\setminus e$ is not 3-connected.
We then have:

\begin{proposition}
   \label{prop.min3conn.graphs}
Let $\mG_p$ be the class of graphs obtained from $K_p$ by
disjoint union, gluing at a cut vertex, series and parallel connection,
deletion and contraction of edges, and deletion of vertices.
Then, for $p\ge3$, $\mG_p$ is minor-closed, and the
minimal excluded minors for $\mG_p$ are the minimally 3-connected graphs
on $p+1$ vertices.
\end{proposition}

\proof
Let $\mH_p$ be the class of graphs with no minimally 3-connected minor on
$p+1$ vertices.  It is clear that $\mG_p$ is minor-closed, so our aim is
prove that $\mG_p=\mH_p$.

We first show that $\mG_p\subseteq\mH_p$.  It is clear that $\mH_p$ is
minor-closed, and is closed under disjoint union (``0-sum'')
and gluing at a cut vertex (``1-sum'').
It is easily checked that if a graph $G$ is obtained by
parallel connection of graphs $G_1$ and $G_2$,
then any 3-connected minor of $G$ is a minor of either $G_1$ or $G_2$;
it follows that $\mH_p$ is closed under parallel connection.
Since series connection can be obtained by combining the other operations
(exploiting $K_3 \in \mH_p$),
we conclude that $\mH_p$ is also closed under series connection.
Finally, we note that $K_p\in\mH_p$,
and as $\mG_p$ is the closure of $\{K_p\}$ under these operations,
we see that $\mG_p\subseteq\mH_p$.

We now show that $\mH_p\subseteq\mG_p$.  For suppose otherwise, and choose
$G\in\mH_p\setminus\mG_p$ with the minimal number of vertices.
Clearly $G$ is 2-connected and has at least $p+1$ vertices.
If $G$ is not 3-connected, then $G$ has a cutset $\{x,y\}$
and there is a decomposition $G=G_1\cup G_2$
where $G_1$ and $G_2$ are connected graphs with at least 3 vertices
and $V(G_1)\cap V(G_2)=\{x,y\}$.  Now let $G_1'$, $G_2'$ be the graphs
obtained from $G_1$, $G_2$ by adding the edge $xy$ if not present.
Then $G_1'$ and $G_2'$ are both minors of $G$
(obtained by contracting the other side to a single edge).
As $\mH_p$ is minor-closed, we have $G_1', G_2'\in\mH_p$.
Therefore, by minimality of $G$ we have $G_1',G_2'\in\mG_p$.
But $G$ can be obtained by taking a parallel connection of $G_1'$ and $G_2'$
along $xy$ and then deleting the edge $xy$ if necessary,
yielding $G \in \mG_p$, contrary to hypothesis.
We conclude that $G$ must be 3-connected.

We now use the fact that
every 3-connected graph other than $K_4$ has an edge
that can be contracted to produce another 3-connected graph
\cite[Lemma~3.2.4]{Diestel_10}.
Contracting suitable edges, we see that $G$ has a 3-connected minor on
$p+1$ vertices, which in turn (deleting edges if necessary) contains a
minimally 3-connected minor on $p+1$ vertices.
But this contradicts $G\in\mH_p$.
\qed

The matroid case is analogous.  We refer to \cite[Chapter~8]{Oxley_11}
for the definitions of 3-connectedness and minimal 3-connectedness
for matroids.
The following result and its proof are due to Oxley \cite{Oxley_private}:

\begin{proposition}
   \label{prop.min3conn.matroids}
Let $\mF$ be a class of matroids that is closed under minors,
direct sums and 2-sums.
Let $r \ge 1$,
let $\mF_r = \{ M \in \mF \colon\, {\rm rank}(M) \le r \}$,
and let $\mF'_r$ denote the class of matroids that can be obtained
from $\mF_r$ by direct sums and 2-sums.
Then $\mF'_r$ is closed under minors
(and of course also under direct sums and 2-sums),
and it consists of all matroids in $\mF$
that have no $\mF^\circ_{r+1}$ minor, where
\be
   \mF^\circ_{r+1}
   \;=\;
   \{ M \in \mF \colon\:  M \hbox{ is minimally 3-connected and }
        {\rm rank}(M) = r+1 \}
   \;.
\ee
\end{proposition}

\proof
Consider first the case $r=1$.
The unique minimally 3-connected matroid of rank 2 is $U_{2,3}$.
If $U_{2,3} \in \mF$, then the result clearly holds;
if $U_{2,3} \notin \mF$, then $\mF'_1 = \mF$
and the result again holds.

Now assume $r \ge 2$, and let $\mH_r$ denote the class of matroids in $\mF$
that have no $\mF^\circ_{r+1}$ minor.
Clearly $\mH_r$ is minor-closed, and our goal is to show that $\mF'_r = \mH_r$.

We first show that $\mF'_r \subseteq \mH_r$.
It is easy to see \cite[Propositions~4.2.20 and 8.3.5]{Oxley_11}
that if a matroid $M$ is either a direct sum or a 2-sum
of matroids $M_1$ and $M_2$, then any 3-connected minor of $M$
is isomorphic to a minor of either $M_1$ or $M_2$;
it follows that $\mH_r$ is closed under direct sum and 2-sum.
Since $\mF_r \subseteq \mH_r$, and $\mF'_r$ is the closure of $\mF_r$
under direct sum and 2-sum, we see that $\mF'_r \subseteq \mH_r$.

We now show that $\mH_r \subseteq \mF'_r$.  For suppose otherwise,
and choose $M \in \mH_r \setminus \mF'_r$ with the minimal number of elements.
It is not hard to see that $M$ must be 3-connected.\footnote{
   If $M$ were a direct sum or a 2-sum of matroids $M_1$ and $M_2$,
   each having at least one element, then $M_1$ and $M_2$ would be
   minors of $M$ \cite[Proposition~7.1.21]{Oxley_11},
   hence $M_1, M_2 \in \mH_r$ because $\mH_r$ is minor-closed.
   But $M_1$ and $M_2$ cannot both belong to $\mF'_r$ because
   $\mF'_r$ is closed under direct sum and 2-sum and $M \notin \mF'_r$.
   Therefore $M_1$ or $M_2$ would be a counterexample to the
   minimality of $M$.
}
Note also that ${\rm rank}(M) > r$ since $M \in \mF \setminus \mF'_r$.

We now use Tutte's Wheels and Whirls Theorem
\cite[Theorem~8.8.4]{Oxley_11},
which says that every 3-connected matroid $N$
that is not a wheel or a whirl has an element $e$ such that
either $N \setminus e$ or $N/e$ is 3-connected (or both).
On the other hand, if $N$ is a wheel or a whirl,
then by contracting a rim element and deleting one of the spokes
adjacent to that rim element, we obtain a 3-connected minor $N'$ of $N$
such that ${\rm rank}(N') = {\rm rank}(N) - 1$.

So we apply this argument repeatedly to $M$ until we obtain a 3-connected minor
$M'$ of $M$ with ${\rm rank}(M') = r+1$.  We then delete elements from $M'$
while maintaining 3-connectedness until we arrive at a minimally
3-connected matroid $M''$.  Therefore $M$ has a minor
$M'' \in \mF^\circ_{r+1}$, contradicting the hypothesis that $M \in \mH_r$.
%
%
%
\qed

Proposition~\ref{prop.min3conn.matroids} with $\mF = $ regular matroids
gives an excluded-minor characterization of the class of matroids
handled by Proposition~\ref{prop.BM.extended1}.
We leave it as a problem for readers more expert in matroid theory
than ourselves to provide an analogous characterization
for Proposition~\ref{prop.BM.extended2}.

%
%



\appendix

\section{The Moore determinant for quaternionic hermitian matrices}
  \label{sec.moore}

In this appendix we review the definition and properties
of the Moore determinant for quaternionic hermitian matrices.
For the most part we have followed Aslaksen \cite{Aslaksen_96}
and Alesker \cite{Alesker_03,Alesker_05}.

The {\em (real) quaternions}\/ $\HH$
are the associative algebra over $\R$
consisting of objects of the form
\be
   q \;=\; x_0 1 + x_1 i + x_2 j + x_3 k  \qquad (x_0,x_1,x_2,x_3 \in \R) \,,
\ee
equipped with the multiplication law where $1$ is the identity element and
\be
   i^2 = j^2 = k^2 = -1 ,\quad
   ij = -ji = k ,\quad
   jk = -kj = i ,\quad
   ki = -ik = j \;.
\ee
The {\em conjugate}\/ of a quaternion $q = x_0 1 + x_1 i + x_2 j + x_3 k$
is $\bar{q} = x_0 1 - x_1 i - x_2 j - x_3 k$.
Note that conjugation is an involutive antiautomorphism of~$\HH$,
i.e.\ $\bar{\bar{q}} = q$ and $\overline{uv} = \bar{v} \bar{u}$.
We also write $|q|^2 = q \bar{q} = \bar{q} q = x_0^2 + x_1^2 + x_2^2 + x_3^2$.
The complex numbers $a+bi$ can be identified
with the subalgebra of $\HH$ consisting of quaternions $a1 + bi + 0j + 0k$.

The {\em hermitian conjugate}\/ of a matrix $M \in \HH^{n \times n}$
is $M^* = (\bar{M})^{\rm T}$.
The map $M \mapsto M^*$ is an involutive antiautomorphism of
the algebra of $n \times n$ quaternionic matrices.
A matrix $M$ is called {\em hermitian}\/ if $M = M^*$.
We denote by ${\rm Herm}(n,\HH)$ the set of $n \times n$
hermitian quaternionic matrices.

The map $\psi \colon\, \HH \to \C^{2 \times 2}$
defined by
\be
   \psi(x_0 1 + x_1 i + x_2 j + x_3 k)
   \;=\;
   \left(\! \begin{array}{cc}
                  x_0 + i x_1   &  x_2 + i x_3 \\
                  -x_2 + i x_3  &  x_0 - i x_1
            \end{array}
   \!\right)
\ee
is an injective *-homomorphism of the algebra of quaternions
into the algebra of $2 \times 2$ complex matrices.
It can equivalently be written as
\be
   \psi(a+bj)  \;=\;  \left(\! \begin{array}{cc}
                                   a  &  b \\
                                   -\bar{b} & \bar{a}
                               \end{array}
                      \!\right)
   \quad\hbox{for } a,b \in \C \;.
\ee
More generally, the map $\Psi \colon\, \HH^{n \times n} \to \C^{2n \times 2n}$
defined by
\be
   \Psi(A+Bj)  \;=\;  \left(\! \begin{array}{cc}
                                   A  &  B \\
                                   -\bar{B} & \bar{A}
                               \end{array}
                      \!\right)
   \quad\hbox{for } A,B \in \C^{n \times n}
 \label{def.Psi}
\ee
is an injective *-homomorphism
of the algebra of $n \times n$ quaternionic matrices
into the algebra of $2n \times 2n$ complex matrices.
Its image is
\be
   \Psi(\HH^{n \times n})
   \;=\;
   \{Z \in \C^{2n \times 2n} \colon\:  JZ = \bar{Z} J \}
\ee
where
\be
   J  \;=\;  \left(\! \begin{array}{cc}
                         0    &  I_n \\
                         -I_n &  0
                      \end{array}
             \!\right)
   \:.
\ee
Note that $\Psi(M^*) = \Psi(M)^* = -J \Psi(M)^{\rm T} J$;
in particular, $M$ is hermitian if and only if $\Psi(M)$ is.

It is convenient now to define $\Phi(M) = J \Psi(M)$, i.e.
\be
   \Phi(A+Bj)  \;=\;  \left(\! \begin{array}{cc}
                                   -\bar{B}  &  \bar{A} \\
                                   -A        &  -B
                               \end{array}
                      \!\right)
   \quad\hbox{for } A,B \in \C^{n \times n}
   \;.
 \label{def.Phi}
\ee
The image of $\Phi$ is the same as that of $\Psi$, i.e.
\be
   \Phi(\HH^{n \times n})
   \;=\;
   \{Y \in \C^{2n \times 2n} \colon\:  JY = \bar{Y} J \}
   \;.
\ee
Note that $\Phi(I) = J$ and $\Phi(M^*) = -\Phi(M)^{\rm T}$.
In particular, $M$ is hermitian if and only if $\Phi(M)$ is antisymmetric,
and we have
\be
   \Phi({\rm Herm}(n,\HH))
   \;=\;
   \{Y \in \C^{2n \times 2n} \colon\:
             Y = -Y^{\rm T} \hbox{ and } JY = \bar{Y} J \}
   \;.
\ee

Recall now that on the space ${\rm Antisym}(2n,\C)$
of $2n \times 2n$ complex antisymmetric matrices
there is a polynomial $\pf$, called the {\em pfaffian}\/,
satisfying $(\pf A)^2 = \det A$ and $\pf J = 1$.
Indeed, the pfaffian is uniquely defined by these conditions.
Here are some of its properties:

\begin{lemma}[Properties of the pfaffian]
   \label{lemma.properties.pfaffians}
The map $\pf \colon\, {\rm Antisym}(2n,\C) \to \C$
has the following properties:
\begin{itemize}
   \item[(a)] $\pf$ is a homogeneous polynomial of degree $n$,
      with integer coefficients.\footnote{
   More precisely, if $A =(a_{ij}) \in {\rm Antisym}(2n,\C)$,
   then $\pf A$ is a homogeneous polynomial with integer coefficients
   in the variables $\{a_{ij}\}_{1 \le i < j \le n}$.
   Indeed, all the coefficients lie in $\{-1,0,1\}$.
}
   \item[(b)] $\pf J = 1$.
   \item[(c)]  $(\pf A)^2 = \det A$.
   \item[(d)]  $\pf(X A X^{\rm T}) = (\det X) (\pf A)$
      for any $2n \times 2n$ matrix $X$.
   \item[(e)]  {\bf (minor summation formula for pfaffians
        \cite{Ishikawa_95,Ishikawa_06}\footnote{
    See also \cite[Theorem~A.15 and remarks after it]{CSS_cayley}
    for an alternative proof using Grassmann--Berezin integration.
})\ }
      More generally, we have
\be
    \pf(X A X^{\rm T})   \;=\;
    \sum\limits_{\begin{scarray}
                    I \subseteq [2n] \\
                    |I| = 2m
                 \end{scarray}}
    (\det X_{\star I}) \, (\pf A_{II})
 \label{eq.app.pfaffian.minorsummation}
\ee
 for any $2m \times 2n$ matrix $X$ ($m \le n$).
 Here $X_{\star I}$ denotes the submatrix of $X$ with columns $I$
 (and all its rows),
 and $A_{II}$ denotes the submatrix of $A$ with rows and columns $I$.
\end{itemize}
\end{lemma}

\noindent
See \cite{Stembridge_90,Knuth_96,Fulton_98,Hamel_01,Lang_02,Ishikawa_06,%
Fulmek_10,CSS_cayley}
for further information on pfaffians.

Let us also observe the following useful fact about pfaffians:
If $B$ is an arbitrary $n \times n$ complex matrix,
then
\be
   \pf \left(\!\! \begin{array}{cc}
                     0  &  B \\
                    -B^{\rm T} &  0
                \end{array}
      \!\!\right)
   \;=\;
   \pf \left[ \left(\!\! \begin{array}{cc}
                           0  &  -B \\
                           I  &  0
                       \end{array}
              \!\!\right)
              \left(\!\! \begin{array}{cc}
                           0  &  I \\
                          -I  &  0
                       \end{array}
              \!\!\right)
              \left(\!\! \begin{array}{cc}
                           0  &  I \\
                           -B^{\rm T}  &  0
                       \end{array}
              \!\!\right)
        \right]
   \;=\;
   \det \left(\!\! \begin{array}{cc}
                           0  &  -B \\
                           I  &  0
                       \end{array}
              \!\!\right)
   \;=\;
   \det B
 \label{eq.pf.block}
\ee
where the second equality used Lemma~\ref{lemma.properties.pfaffians}(b,d)
and the last equality used row (or column) interchanges.

We now define a ``determinant'' for {\em hermitian}\/
quaternionic matrices only:

\begin{definition}[Moore determinant of a hermitian quaternionic matrix]
\hfill\break
The {\em Moore determinant} of a matrix $M \in {\rm Herm}(n,\HH)$
is defined by
\be
   \Mdet(M)  \;=\;  \pf(\Phi(M))  \;=\;  \pf(J \Psi(M))
   \;.
 \label{def.Mdet}
\ee
\end{definition}

\noindent
Moore's original definition \cite{Moore_22} was very different from this one,
but the two definitions can be proven to be equivalent
\cite[pp.~63--64]{Aslaksen_96}
\cite[pp.~141--152, especially Theorem~8.9.4]{Mehta_89}.
It also turns out that the Moore determinant
is the same as the Jordan determinant on the Euclidean Jordan algebra
${\rm Herm}(n,\HH)$:
see \cite[Exercise~II.7 (pp.~39--40), Exercise~III.1 (p.~58),
and pp.~84 and 88]{Faraut_94}.
This equivalence is what allows us to formulate
   Theorems~\ref{thm1.det}(c) and \ref{thm1.det.cones}(c)
   in terms of the Moore determinant,
   while the key Theorem~\ref{thm1.det.cones.Jordan}
   is stated and proven in terms of the Jordan determinant.


\begin{proposition}[Properties of the Moore determinant]
   \label{prop.properties.Mdet}
\hfill\break
The map $\Mdet \colon\, {\rm Herm}(n,\HH) \to \C$
has the following properties:
\begin{itemize}
   \item[(a)] $\Mdet$ is a homogeneous polynomial of degree $n$,
      with integer coefficients.\footnote{
   More precisely, if for $M =(m_{ij}) \in {\rm Herm}(n,\HH)$ we write
   $m_{ij} = \alpha_{ij} 1 + \beta_{ij} i + \gamma_{ij} j + \delta_{ij} k$
   with $\alpha_{ij}, \beta_{ij}, \gamma_{ij}, \delta_{ij} \in \R$,
   then $\Mdet(M)$ is a homogeneous polynomial with integer coefficients
   in the variables $\{\alpha_{ii}\}_{1 \le i \le n}$
   and
   $\{\alpha_{ij}, \beta_{ij}, \gamma_{ij}, \delta_{ij} \}_{1 \le i < j \le n}$.
}
   \item[(b)] $\Mdet I = 1$.
   \item[(c)] $\Mdet$ is real-valued.
   \item[(d)] If $M$ is a {\em complex} hermitian matrix,
      then its Moore determinant equals its ordinary determinant.
   \item[(e)] If $X \in \HH^{n \times n}$ (not necessarily hermitian),
      then $\Mdet(XX^*) = \det(\Phi(X)) = \det(\Psi(X)) \ge 0$.
   \item[(f)] If $M \in {\rm Herm}(n,\HH)$ and $X \in \HH^{n \times n}$,
      then $\Mdet(XMX^*) = \Mdet(XX^*) \, \Mdet(M)$.
   \item[(g)] {\bf (restricted Cauchy--Binet formula for the Moore determinant)\ }
      If $X \in \HH^{m \times n}$ with $m \le n$, 
      and $D \in \R^{n \times n}$ is a real diagonal matrix, then
\be
    \Mdet(XDX^*)  \;=\;
    \sum\limits_{\begin{scarray}
                    I \subseteq [n] \\
                    |I| = m
                 \end{scarray}}
    \Mdet [X_{\star I} (X_{\star I})^*]  \, \det(D_{II})  \;.
 \label{eq.cauchy-binet}
\ee
   \item[(h)] {\bf (restricted multiplicativity)\ }
      If $M,N \in {\rm Herm}(n,\HH)$ with $MN=NM$, then
      $\Mdet(MN)  =  \Mdet(M) \, \Mdet(N)$.
\end{itemize}
\end{proposition}

\proof
(a) is an immediate consequence of Lemma~\ref{lemma.properties.pfaffians}(a)
together with the definitions \reff{def.Mdet}
and \reff{def.Psi}/\reff{def.Phi}.

(b) is likewise immediate from Lemma~\ref{lemma.properties.pfaffians}(b).

(c) Note that for any $Y \in \Phi({\rm Herm}(n,\HH))$
we have $JYJ^{\rm T} = \bar{Y}$, hence
\be
   \overline{\pf(Y)}  \;=\;  \pf(\bar{Y})
   \;=\; \pf(JYJ^{\rm T})  \;=\;  (\det J) \, \pf(Y)  \;=\;  \pf(Y)
\ee
by Lemma~\ref{lemma.properties.pfaffians}(d).

(d)
If $M$ is a complex hermitian matrix, then
$\Phi(M) = \displaystyle \left(\! \begin{array}{cc}
                                   0  &  M^{\rm T} \\
                                   -M &  0
                               \end{array}
                      \!\right)$
and $\Mdet(M) = \pf \displaystyle \left(\! \begin{array}{cc}
                                   0  &  M^{\rm T} \\
                                   -M &  0
                               \end{array}
                      \!\right) = \det M$
by \reff{eq.pf.block}.


(e,f)  We have
$\Mdet(XMX^*) = \pf[J \Psi(XMX^*)]  = \pf[J \Psi(X) \Psi(M) \Psi(X^*)]
   =$
\break
$\pf[\Phi(X) J \Phi(M) J \Phi(X^*)]
   = \pf[- \Phi(X) J \Phi(M) J \Phi(X)^{\rm T}]
   = \pf[(\Phi(X) J) \, \Phi(M) \, (\Phi(X) J)^{\rm T}]
   = \det(\Phi(X) J) \, \pf(\Phi(M))
   = \det(\Phi(X)) \, \Mdet(M)$.
Applying this with $M=I$,
we see that $\Mdet(XX^*) = \det(\Phi(X)) = \det(\Psi(X))$.
Therefore, $\Mdet(XMX^*) = \Mdet(XX^*) \, \Mdet(M)$.

We will prove in Lemma~\ref{lemma.funnypositivity}(c) below
that $\det(\Psi(X)) \ge 0$.

(g) Define $\Psi_{mn} \colon\, \HH^{m \times n} \to \C^{2m \times 2n}$
by the same formula \reff{def.Psi};
then it is a homomorphism in the sense that if
$M \in \HH^{m \times n}$ and $N \in \HH^{n \times p}$,
then
\be
   \Psi_{mp}(MN)  \;=\;  \Psi_{mn}(M) \, \Psi_{np}(N)  \;.
\ee
We also have, for $M \in \HH^{m \times n}$,
\be
   \Psi_{nm}(M^*)  \;=\;  -J_n \Psi_{mn}(M)^{\rm T} J_m
   \;.
\ee
Then we have
\begin{eqnarray}
   \Mdet(XDX^*)
   & = &
   \pf[J_m \Psi_{mm}(XDX^*)]
        \nonumber \\[1mm]
   & = &
   \pf[J_m \Psi_{mn}(X) \Psi_{nn}(D) \Psi_{nm}(X^*)]
        \nonumber \\[1mm]
   & = &
   \pf[J_m \Psi_{mn}(X) [-J_n \Phi(D)] [-J_n \Psi_{mn}(X)^{\rm T} J_m]]
        \nonumber \\[1mm]
   & = &
   \pf[(J_m \Psi_{mn}(X) J_n) \, \Phi(D) \, (J_m \Psi_{mn}(X) J_n)^{\rm T}]
        \nonumber \\[1mm]
   & = &
   \!\!
    \sum\limits_{\begin{scarray}
                    I \subseteq [2n] \\
                    |I| = 2m
                 \end{scarray}}
    \det[(J_m \Psi_{mn}(X) J_n)_{\star I}] \: (\pf \Phi(D)_{II})
   \;.
\end{eqnarray}
Now write $I = K \cup L$ with
$K \subseteq [n]$, $L \subseteq [2n] \setminus [n]$ and $|K| + |L| = 2m$.
We claim that $\pf \Phi(D)_{II} = 0$ unless $L = K+n$.
To see this, note first that for any complex matrix~$D$ (diagonal or not)
we have
\be
   \Phi(D)  \;=\;  J \Psi(D)
   \;=\;
       \left(\! \begin{array}{cc}
                         0    &  \bar{D} \\
                         -D   & 0
                      \end{array}
             \!\right)
   \,.
\ee
If in addition $D$ is diagonal and $L \neq K+n$,
then $\Phi(D)_{II}$ has a zero row (and a zero column),
so that $\det \Phi(D)_{II} = 0$.
If furthermore $D$ is real, then $\Phi(D)_{II}$ is antisymmetric
and $(\pf \Phi(D)_{II})^2 = \det \Phi(D)_{II} = 0$.
So only the terms $L = K+n$ survive in the sum over $I$.

Then in this case we have
\be
   \pf \Phi(D)_{II}
   \;=\;
   \pf \left(\! \begin{array}{cc}
                     0      &  D_{KK} \\
                    -D_{KK} &  0
                \end{array}
       \!\right)
    \;=\;
    \det D_{KK}
\ee
by \reff{eq.pf.block}.
Moreover, a simple calculation shows that
\be
   (J_m \Psi_{mn}(X) J_n)_{\star I}
   \;=\;
   J_m \Psi_{mm}(X_{\star K}) J_m
   \;;
\ee
and writing $Y = X_{\star K}$ we see easily that
$\det(J \Psi(Y) J) = \det(\Phi(Y) J) = \det(\Phi(Y)) = \Mdet(YY^*)$
by part (e), hence
\be
   \det[(J_m \Psi_{mn}(X) J_n)_{\star I}]
   \;=\;
   \Mdet [X_{\star K} (X_{\star K})^*]
   \;.
\ee

(h) If $M,N \in {\rm Herm}(n,\HH)$ with $MN=NM$,
we have $MN \in {\rm Herm}(n,\HH)$, so its Moore determinant is well-defined;
then
$[\Mdet(MN)]^2 = [\pf(\Phi(MN))]^2 = \det(\Phi(MN)) = \det(\Psi(MN))
 = \det(\Psi(M) \Psi(N)) = \det(\Psi(M)) \, \det(\Psi(N)) =$
\linebreak  
$\det(\Phi(M)) \, \det(\Phi(N)) = [\pf(\Phi(M))]^2 [\pf(\Phi(N))]^2
 = [\Mdet(M)]^2 [\Mdet(N)]^2$.
Now apply this with $M,N$ replaced by
$M_\lambda = (1-\lambda) I + \lambda M$,
$N_\lambda = (1-\lambda) I + \lambda N$
where $\lambda \in \R$.
Then $P(\lambda) = \Mdet(M_\lambda N_\lambda)$,
$Q(\lambda) = \Mdet(M_\lambda)$ and $R(\lambda) = \Mdet(N_\lambda)$
are polynomials satisfying
$P(\lambda)^2 = Q(\lambda)^2 R(\lambda)^2$ and $P(0) = Q(0) R(0)$.
It follows that $P(\lambda) = Q(\lambda) = R(\lambda) = 1$.
Now evaluate at $\lambda=1$.
\qed

\begin{lemma}   
   \label{lemma.funnypositivity}
\begin{itemize}
   \item[(a)]  Let $M \in \C^{n \times n}$.
      Then the eigenvalues of $M \bar{M}$ are real
      or come in complex-conjugate pairs;
      and the negative real eigenvalues have even algebraic multiplicity.
   \item[(b)]  If $M \in \C^{n \times n}$ and $\kappa \ge 0$,
      then $\det(\kappa I + M \bar{M}) \ge 0$.
   \item[(c)]  If $A,B \in \C^{n \times n}$, then
\be
       \det \left(\! \begin{array}{cc}
                                   A  &  B \\
                                   -\bar{B} & \bar{A}
                               \end{array}
                      \!\right)
       \;\ge\;  0  \;.
\ee
\end{itemize}
\end{lemma}

\proof
(a)  It is well known that $AB$ and $BA$ have the same
characteristic polynomial.
Therefore, $M \bar{M}$ and $\bar{M} M$ have the same characteristic
polynomial, i.e.\ the coefficients of this characteristic polynomial
are real.  It follows that the eigenvalues of $M \bar{M}$
are real or come in complex-conjugate pairs.

If $M \bar{M}$ has $n$ distinct eigenvalues, then we can argue as follows:
Suppose that $M \bar{M} x = \lambda x$ with $x \neq 0$.
Then by taking complex conjugates we have
$\bar{M} M \bar{x} = \bar{\lambda} \bar{x}$;
and left-multiplying by $M$ we obtain
$M \bar{M} (M \bar{x}) = \bar{\lambda} (M \bar{x})$.
If $\lambda$ is real, the distinct-eigenvalues hypothesis implies
that $M \bar{x} = \mu x$ for some $\mu \in \C$.
By complex-conjugating we get $\bar{M} x = \bar{\mu} \bar{x}$,
and left-multiplying by $M$ we obtain
$M \bar{M} x = \bar{\mu} M \bar{x} = |\mu|^2 x$.
Since $x \neq 0$ it follows that $\lambda = |\mu|^2 \ge 0$.
Therefore all real eigenvalues of $M \bar{M}$ are nonnegative.

Let us now observe that the matrices $M$ for which
$M \bar{M}$ has $n$ distinct eigenvalues are dense in $\C^{n \times n}$.
To see this, let us consider the real and imaginary parts of
the entries of $M$ to be $2n^2$ distinct indeterminates,
and let us form the discriminant of the characteristic polynomial
of $M \bar{M}$.
This is a polynomial (with real coefficients, though we don't need this fact)
in the $2n^2$ indeterminates,
and its zero set in $\R^{2n^2}$ corresponds precisely to the matrices
$M \in \C^{n \times n}$
for which $M \bar{M}$ does not have $n$ distinct eigenvalues.
Now the zero set of a polynomial in $\R^N$
is either all of $\R^N$ or else a proper subvariety
(which in particular has empty interior).
Since there do exist matrices $M \in \C^{n \times n}$ for which
$M \bar{M}$ has $n$ distinct eigenvalues,
it follows that such matrices form a dense open set in $\C^{n \times n}$.

Since the eigenvalues depend continuously on the matrix,
it follows by density and continuity that,
for arbitrary $M \in \C^{n \times n}$,
the negative real eigenvalues of $M \bar{M}$ have even algebraic multiplicity.

(b) follows immediately from (a), using the fact that the
determinant of a matrix is the product of its eigenvalues
taken with their algebraic multiplicity.

(c) If $A$ is invertible, then we have (writing $C = A^{-1} B$)
\be
        \left(\! \begin{array}{cc}
                     A  &  B \\
                     -\bar{B} & \bar{A}
                 \end{array}
        \!\right)
        \;=\;
        \left(\! \begin{array}{cc}
                     A  &  0 \\
                     0 & \bar{A}
                 \end{array}
        \!\right)
        \left(\! \begin{array}{cc}
                     I  &  C \\
                     -\bar{C} & I
                 \end{array}
        \!\right)
        \;=\;
        \left(\! \begin{array}{cc}
                     A  &  0 \\
                     0 & \bar{A}
                 \end{array}
        \!\right)
        \left(\! \begin{array}{cc}
                     I  & C \\
                     0  & I
                 \end{array}
        \!\right)
        \left(\! \begin{array}{cc}
                     I+C\bar{C}  &  0 \\
                     -\bar{C}    & I
                 \end{array}
        \!\right)
        \,,
\ee
from which it follows that
\be
    \det\left(\! \begin{array}{cc}
                     A  &  B \\
                     -\bar{B} & \bar{A}
                 \end{array}
        \!\right)
    \;=\;
    |\det A|^2 \: \det(I+C\bar{C})
    \;,
\ee
which is nonnegative by part (b).
Since the invertible matrices are dense in $\C^{n \times n}$,
the general result follows by continuity.
\qed

\bigskip

{\bf Remarks on Proposition~\ref{prop.properties.Mdet}.} \quad
1. The assertion in Proposition~\ref{prop.properties.Mdet}(e)
that $\Mdet(XX^*) \ge 0$ can alternatively be proven as follows:
Observe that $XX^*$ is positive-semidefinite
and hence has nonnegative real eigenvalues
\cite[Corollary~5.3 and Proposition~5.2]{Farenick_03};
so by the spectral theorem
\cite[Theorem~3.3 and Proposition~3.8]{Farenick_03}
we can write $XX^* = UDU^*$ with $UU^* = U^* U = I$
and $D$ diagonal with nonnegative real entries;
it then follows from part (f)
that $\Mdet(XX^*) = \Mdet(UU^*) \Mdet(D) = \Mdet(D)$,
which is $\ge 0$ by part (d).

2.  The restriction in part (g) to {\em real diagonal}\/ matrices
is essential:  indeed, already in the first nontrivial case $m=1$, $n=2$
one sees that
\be
   \left( r \;\: s \right)
   \left( \! \begin{array}{cc}
                  a & q \\
                  \bar{q} & b
             \end{array}
   \!\right)
   \left( \! \begin{array}{c}
                  \bar{r} \\
                  \bar{s}
             \end{array}
   \!\right)
   \;=\;
   r a \bar{r} + r q \bar{s} + s \bar{q} \bar{r} + s b \bar{s}
\ee
while the right side of \reff{eq.cauchy-binet} is
$(r \bar{r}) a + (s \bar{s}) b$.

To see why part (g) states the appropriate Cauchy--Binet formula for the
Moore determinant, recall first the ordinary Cauchy--Binet formula for
matrices over a {\em commutative}\/ ring $R$:
for $X \in R^{m \times n}$, $M \in R^{n \times n}$ and $Y \in R^{n \times m}$, 
we have
\be
    \det(XMY)  \;=\!
    \sum\limits_{\begin{scarray}
                    I,J \subseteq [n] \\
                    |I| = |J| = m
                 \end{scarray}}
    \!
    (\det X_{\star I}) \, (\det M_{IJ}) \, (\det Y_{J \star})
   \;.
\ee
But in our case the determinant is defined only for hermitian matrices,
so we need to
\begin{quote}
\begin{itemize}
   \item[(i)] require $Y=X^*$ and $M = M^*$ so that $XMY$ is hermitian,
   \item[(ii)] restrict to $I=J$ so that $M_{IJ}$ is hermitian, and
   \item[(iii)] rewrite $(\det X_{\star I}) (\det (X^*)_{I \star})$
       as $\det [X_{\star I} (X_{\star I})^*]$.
\end{itemize}
\end{quote}
But restricting the sum to $I=J$ is correct in the commutative case
(e.g.\ if all the quaternionic matrix elements happen to belong to $\R$)
only if $\det M_{IJ} = 0$ whenever $I \neq J$,
and this happens if and only if $M$ is diagonal
(the ``only if'' is seen by considering $m=1$).
So we really need to require that $M$ be a real diagonal matrix,
at least when $m=1$.
On the other hand, in the square case $m=n$ we have $I=J=[n]$ automatically,
so this case does {\em not}\/ require $M$ to be diagonal
[cf.\ part (f)].

The special case $D=I$ of \reff{eq.cauchy-binet}
was proven by Liebend\"orfer \cite[Theorem~1]{Liebendorfer_05}.
We do not know whether the general result is new.

3.  Regarding part (h), note that
the product $MN$ is hermitian if and only if $M$ and $N$ commute,
so this hypothesis is required for $\Mdet(MN)$ to be well-defined.

Part (h) can alternatively be proven by invoking the
spectral theorem for commuting hermitian quaternionic matrices
\cite[Theorem~3.3 and Propositions~3.5, 3.6 and 3.8]{Farenick_03}:
if $M,N \in {\rm Herm}(n,\HH)$ with $MN=NM$,
then there exist a unitary matrix $U \in \HH^{n \times n}$
(i.e.\ $UU^* = U^* U = I$)
and real diagonal matrices $D,E$ such that
$U^* D U = M$ and $U^* E U = N$.
The result is then an easy consequence of parts (d) and (f).
\qed

\medskip

{\bf Remarks on Lemma~\ref{lemma.funnypositivity}.} \quad
1. Lemma~\ref{lemma.funnypositivity}(a)
goes back at least to Asano and Nakayama \cite[S\"atze~19--21]{Asano_38},
who show that $M \bar{M}$ is similar to the square of a real matrix.
It is part of more general theorem of Youla \cite[Theorem~1]{Youla_61}
that gives a normal form for matrices $M \in \C^{n \times n}$
under unitary congruence, i.e.\ $M = U^{\rm T} \Sigma U$
with $U$ unitary and $\Sigma$ of a special form.
See also \cite[p.~147]{Hong_88}, \cite{Fassbender_07},
\cite[Sections~4.4 and 4.6, especially pp.~252--253]{Horn_00}
and
\cite[Sections~4.4 and 4.6, especially Theorem~4.4.9 and
   Corollaries~4.4.13 and 4.6.16]{Horn_12}.
Lemma~\ref{lemma.funnypositivity}(c) goes back at least to
Zhang \cite[Proposition~4.2]{Zhang_97}
--- with the same proof as given here ---
but is probably much older.

2. In Lemma~\ref{lemma.funnypositivity}(a),
it {\em is}\/ possible for $M \bar{M}$ to have negative real eigenvalues:
consider, for instance, the $90^\circ$ rotation matrix
$M = \displaystyle{ \left(\! \begin{array}{cc}
                                   0   & 1  \\
                                   -1  & 0
                               \end{array}
                    \!\right)}$.
It is also possible for $M \bar{M}$ to be nondiagonalizable:
consider, for instance,
$M = \displaystyle{ \left(\! \begin{array}{cc}
                                   1  & \epsilon  \\
                                   0  & 1
                               \end{array}
                    \!\right)}$
with $\epsilon \in \R \setminus \{0\}$.

3. The assertions in Lemma~\ref{lemma.funnypositivity}(a)
about the eigenvalues of $M \bar{M}$ are best possible,
in the sense that every $n$-tuple $(\lambda_1,\ldots,\lambda_n) \in \C^n$
satisfying these conditions is the set of eigenvalues for
some matrix $M \bar{M}$ (indeed, one that is real and diagonalizable).
To see this, it suffices to consider $1 \times 1$ matrices $M = (\lambda)$
with $\lambda \ge 0$,
and $2 \times 2$ matrices
$M = \lambda \displaystyle{ \left(\! \begin{array}{cc}
                                   \cos\theta   & \sin\theta  \\
                                   -\sin\theta  & \cos\theta
                               \end{array}
                    \!\right)}$
with $\lambda \ge 0$ and $\theta \in \R$,
and then form direct sums.

4. If we consider the real and imaginary parts of the entries of
$M \in \C^{n \times n}$ to be $2n^2$ distinct indeterminates,
then Lemma~\ref{lemma.funnypositivity}(b) asserts that $\det(I + M \bar{M})$
is a nonnegative polynomial on $\R^{2n^2}$.
It would be nice to have a more direct proof of this fact.
For instance, is this nonnegative polynomial actually a sum of squares?
And likewise for Lemma~\ref{lemma.funnypositivity}(c).
\qed

\bigskip

We also have a {\em restricted}\/ formula for the Moore determinant of the sum
of two matrices:

\begin{proposition} {$\!\!\!$ \bf \protect\cite[Proposition~1.1.11]{Alesker_03} \ }
   \label{prop.detsum.quaternions}
If $M \in {\rm Herm}(n,\HH)$
and $D \in \R^{n \times n}$ is a real diagonal matrix, then
\be
   \Mdet(D+M)  \;=\;  \sum_{I \subseteq [n]}  (\det D_{II}) (\Mdet M_{I^c I^c})
   \;.
\ee
\end{proposition}

{\bf Remark.}
For matrices over a {\em commutative}\/ ring we have an {\em unrestricted}\/
formula for $\det(A+B)$, which involves a double sum over subsets $I,J$
(Lemma~\ref{lemma.detpoly1}).
But in the quaternionic case we must restrict to $I=J$
in order to ensure that $A_{IJ}$ and $B_{I^c J^c}$ are hermitian;
and to justify this restriction we must require that at least one
of the matrices $A,B$ be diagonal.
Therefore only the restricted formula (Corollary~\ref{cor.detpoly1})
generalizes to quaternions.
\qed

\medskip

For matrices over a noncommutative division ring
(such as the quaternions),
one can in principle define four ranks:
the left and right row-ranks and the left and right column-ranks.
In general the left row-rank equals the right column-rank,
and the right row-rank equals the left column-rank,
but these two numbers need not be equal:
for instance, the matrix
$\Biggl(\!\! \begin{array}{cc}
                      1 & i \\
                      j & k
                 \end{array}
            \!\!\Biggr)$
has left row-rank and right column-rank 2,
but right row-rank and left column-rank 1.
Because matrix multiplication always forms
right linear combinations of columns
and left linear combinations of rows,
it is left row-rank = right column-rank that is well-behaved
with respect to multiplication;
we therefore follow Jacobson \cite[pp.~22, 51]{Jacobson_53}
in focussing attention on these two,
which we call the ``rank'' {\em tout court}\/.

The rank of a general quaternionic matrix can be characterized as follows:

\begin{proposition} {$\!\!\!$ \bf \protect\cite[Theorem~8]{Wolf_36}
                                  \protect\cite[Theorem~7.3]{Zhang_97}
                 \ }
 \label{prop.rank}
$M \in \HH^{n \times n}$ has rank $r$
if and only if $\Psi(M)$ [or equivalently $\Phi(M)$] has rank $2r$.
\end{proposition}

For {\em hermitian}\/ quaternionic matrices we conjecture an alternate
characterization in terms of principal minors,
which if true would extend
a well-known result for complex hermitian matrices\footnote{
   See e.g.\ \cite[Theorem~7.12.4 and Corollary~7.12.5]{Hohn_73}.
}:

\begin{conjecture}
Let $M \in {\rm Herm}(n,\HH)$ have rank $r$.
Then:
\begin{itemize}
   \item[(a)]  All principal minors of $M$ of order $r$
(i.e.\ the quantities $\Mdet M_{II}$ with $|I|=r$) have the same sign,
and at least one of them is nonzero.
   \item[(b)]  All principal minors of $M$ of order $s > r$ vanish.
\end{itemize}
\end{conjecture}

We also have an analogue of Proposition~\ref{prop.detpoly2}
for hermitian quaternionic matrices:

\begin{proposition}
   \label{prop.detpoly2.quaternions}
Let $A_1,\ldots,A_k \in {\rm Herm}(n,\HH)$, and define
\be
   P(x_1,\ldots,x_k)  \;=\; \Mdet\!\left( \sum_{i=1}^k x_i A_i \right)
\ee
for $x_1,\ldots,x_k \in \R$.
Then $P$ is a homogeneous polynomial of degree $n$ with real coefficients.
Furthermore, the degree of $P$ in the variable $x_i$ is $\le {\rm rank}(A_i)$.
\end{proposition}

\proof
It follows immediately from the definition of $\Mdet$
that $P$ is a homogeneous polynomial of degree $n$
that is real-valued when $x_1,\ldots,x_k$ are real.
Moreover we have
\be
   P(x_1,\ldots,x_k)^2
   \;=\;
   \det \Phi\!\left( \sum_{i=1}^k x_i A_i \right)
   \;=\;
   \det \!\left( \sum_{i=1}^k x_i \Phi(A_i) \right)
   \;.
\ee
It then follows from Proposition~\ref{prop.detpoly2}
that the degree of $P^2$ in the variable $x_i$
is at most ${\rm rank}(\Phi(A_i))$,
which equals twice ${\rm rank}(A_i)$ by Proposition~\ref{prop.rank}.
\qed

\section{Elementary proof of Gindikin's criterion for the positivity of
   Riesz distributions}  \label{sec.gindikin}

In this appendix we present an elementary proof
of Gindikin's \cite{Gindikin_75} necessary and sufficient condition
for the positivity of the one-parameter Riesz distributions $\scrr_\alpha$
on a simple Euclidean Jordan algebra (Theorem~\ref{thm.riesz1}).
This proof is due to Shanbhag \cite{Shanbhag_88}
and Casalis and Letac \cite{Casalis_94},
as reworked and simplified by one of us \cite[Appendix]{Sokal_riesz};
it uses only the basic properties of the Riesz distributions as set forth in
Theorem~\ref{thm.basicproperties} and Proposition~\ref{prop.positive}.

\proofof{Theorem~\ref{thm.riesz1}}
By Proposition~\ref{prop.positive},
if $\alpha \in \big\{0,\frac{d}{2},\ldots,(r-1)\frac{d}{2} \big\}
    \cup {\big((r-1)\frac{d}{2},\infty\big)}$,
then $\scrr_\alpha$ is a positive measure.
We shall prove the converse in several steps:

1) If $\alpha \in \C \setminus \R$, then $\scrr_\alpha$
is not even real\footnote{
   This can be seen by restricting the definition \reff{def.scrr}
   to test functions $\varphi$ that have compact support contained in $\Omega$,
   and then analytically continuing both sides in $\alpha$;
   or it can be seen alternatively from the
   Laplace transform \reff{eq.laplace.analcont}.
},
so it is surely not a positive measure.

2) Consider next $(r-2) {d \over 2} < \alpha <  (r-1) {d \over 2}$.
Since $\supp \scrr_\alpha \subseteq \overline{\Omega}$
and $\Delta(x)$ is nonnegative on $\overline{\Omega}$,
if $\scrr_\alpha$ were a positive measure, then so would be
$\Delta(x) \, \scrr_\alpha$,
which by \reff{eq.riesz.product} equals $C_\alpha \scrr_{\alpha+1}$,
where
\be
   C_\alpha \;=\; \prod\limits_{j=0}^{r-1} \Bigl( \alpha - j \frac{d}{2} \Bigr)
   \;<\; 0  \;.
\ee
It follows that $\scrr_{\alpha+1}$ would be a negative measure.
But in fact {\em no}\/ Riesz distribution $\scrr_\beta$ can be a
negative measure, because its Laplace transform \reff{eq.laplace.analcont}
is strictly positive (when $\beta \in \R$)
or nonreal (when $\beta \in \C \setminus \R$).
We conclude that $\scrr_\alpha$ is not a positive measure
for $\alpha \in \big( (r-2) {d \over 2},\, (r-1) {d \over 2} \big)$.

3) By Proposition~\ref{prop.positive}(a), $\scrr_{d/2}$ is a positive measure;
therefore, whenever $\scrr_\alpha$ is a positive measure,
so is $\scrr_{\alpha + d/2} = \scrr_\alpha * \scrr_{d/2}$
[by \reff{eq.riesz.convolution}].
By induction it follows that $\scrr_\alpha$ is not a positive measure for
$\alpha \in \bigcup\limits_{k=1}^\infty
 \big( (r-k-1) {d \over 2},\, (r-k) {d \over 2} \big)$.

4) The only remaining values of $\alpha$ are negative multiples of $d/2$.
But $\scrr_\alpha$ cannot be a positive measure for any $\alpha < 0$:
for if it were, then its Laplace transform would be a decreasing function
on the cone $\Omega$;
but for $\alpha < 0$ the Laplace transform \reff{eq.laplace.analcont}
is in fact an {\em increasing}\/ (to $+\infty$) function on $\Omega$.
\qed


%

\section*{Acknowledgments}

We wish to thank Christian Berg, Petter Br\"and\'en,
Jacques Faraut, James Oxley, Robin Pemantle, Dave Wagner and Geoff Whittle
for helpful conversations and/or correspondence.
In particular, we thank James Oxley for giving us permission
to include his proof of Proposition~\ref{prop.min3conn.matroids}.
We also wish to thank Christian Berg
for drawing our attention to the work of Hirsch \cite{Hirsch_72},
Malek Abdesselam for drawing our attention to the work of
Gurau, Magnen and Rivasseau \cite{Gurau_08},
Roger Horn for valuable information concerning the history of
Lemma~\ref{lemma.funnypositivity},
and Muriel Casalis for providing us with a copy of \cite{Bonnefoy-Casalis_90}.
Finally, we thank an anonymous referee for valuable suggestions
that helped us to improve the organization of this article.

This research was supported in part
by U.S.\ National Science Foundation grant PHY--0424082.



\end{document}